\documentclass[openany,11pt,oneside,a4paper]{book}

\usepackage[T1]{fontenc}
\usepackage[brazil]{babel}
\usepackage[latin1]{inputenc}
\usepackage[pdftex]{graphicx}           
\usepackage{setspace}                   
\usepackage{indentfirst}                
\usepackage{makeidx}                    
\usepackage{amsmath}                    
\usepackage{amsfonts}
\usepackage{amssymb}
\usepackage{amsmath, amsthm}
\usepackage[nottoc]{tocbibind}          
\usepackage{courier}                    
\usepackage{type1cm}                    
\usepackage{listings}                   
\usepackage{titletoc}
\usepackage[fixlanguage]{babelbib}
\usepackage[font=small,format=plain,labelfont=bf,up,textfont=it,up]{caption}
\usepackage[usenames,svgnames,dvipsnames]{xcolor}
\usepackage[a4paper,top=2.54cm,bottom=2.0cm,left=2.27cm,right=2.27cm]{geometry} 
\usepackage[pdftex,plainpages=false,pdfpagelabels,pagebackref,colorlinks=true,citecolor=DarkGreen,linkcolor=NavyBlue,urlcolor=DarkRed,filecolor=green,bookmarksopen=true]{hyperref} 
\usepackage[all]{hypcap}                
\fontsize{60}{62}\usefont{OT1}{cmr}{m}{n}{\selectfont}
\usepackage{subfigure}
\usepackage{psfrag}

\usepackage{fancyhdr}
\renewcommand{\chaptermark}[1]{\markboth{\MakeUppercase{#1}}{}}

\graphicspath{{./figuras/}}             
\makeindex                              
\fontsize{60}{62}\usefont{OT1}{cmr}{m}{n}{\selectfont}
\cleardoublepage
\normalsize

\newtheorem{definicao}{Definição}
\newtheorem{teorema}{Teorema}
\newtheorem{proposicao}{Proposição}
\newtheorem{lema}{Lema}
\newtheorem{exemplo}{Exemplo}
\newtheorem{corolario}{Corolário}
\newtheorem{suposicao}{Suposição}

\newcommand{\R}{\mathbb{R}}
\newcommand{\um}{\mathbf{1}}
\newcommand{\Mat}{\operatorname{Mat}}

\begin{document}
\frontmatter 
\fancyhead[RO]{{\footnotesize\rightmark}\hspace{2em}\thepage}
\setcounter{tocdepth}{2}
\fancyhead[LE]{\thepage\hspace{2em}\footnotesize{\leftmark}}
\fancyhead[RE,LO]{}
\fancyhead[RO]{{\footnotesize\rightmark}\hspace{2em}\thepage}

\onehalfspacing  

\thispagestyle{empty}
\begin{center}
    \vspace*{2.3cm}
    \textbf{\Large{Sincronização em Redes Complexas:\\
     Estabilidade e Persistência}}\\
    
    \vspace*{1.2cm}
    \Large{Marcos Daniel Nogueira Maia}
    
    \vskip 2cm
    \textsc{
    Dissertação apresentada\\[-0.25cm] 
    à\\[-0.25cm]
    Universidade Federal do ABC\\[-0.25cm]
    para\\[-0.25cm]
    obtenção do título\\[-0.25cm]
    de\\[-0.25cm]
    Mestre em Matemática Aplicada}
    
    \vskip 1.5cm
    Programa de Pós-Graduação em Matemática Aplicada\\
    Orientador: Prof. Dr. Tiago Pereira da Silva\\
    Coorientador: Prof. Dr. Rafael de Mattos Grisi

   	\vskip 1cm
    
    \vskip 0.5cm
    \normalsize{Santo André, março de 2013}
\end{center}

%
\newpage
\thispagestyle{empty}
    \begin{center}
        \vspace*{2.3 cm}
        \textbf{\Large{Sincronização em Redes Complexas:\\
        Estabilidade e Persistência}}\\
        \vspace*{2 cm}
    \end{center}

    \vskip 2cm

    \begin{flushright}
	Esta dissertação trata-se da versão original\\
	do aluno Marcos Daniel Nogueira Maia.
    \end{flushright}

\pagebreak

%
%




    
%

\pagenumbering{roman}     

\chapter*{Agradecimentos}
Agradeço primeiramente à Deus, por tornar os meus sonhos possíveis e sempre
me dar força para alcançá-los.

Em especial, agradeço à minha esposa, Lidiana, por seu amor, apoio e compreensão. Por sempre
entender que em vários momentos, tive que abdicar de sua companhia para me dedicar
ao estudo e a consequente produção deste trabalho.

Agradeço aos meus orientadores, o Prof. Tiago Pereira e o Prof. Rafael Grisi,
por sua dedicação em sempre estar disponível para me ajudar. Ao Prof. Tiago Pereira,
em especial agradeço, por me receber como orientando, mesmo sem me conhecer, a princípio.
Pela paciência e disposição
que tiveram para me ensinar, desde o mais simples ao mais complexo. Pelo incentivo
ao estudo disciplinado e à pesquisa, e pelos conselhos de vida, agradeço.

Aos professores da pós-graduação,
por me fornecerem o suporte matemático tão precioso e necessário
para a minha formação, agradeço.

Aos amigos da pós-mat, em especial ao José, Rafael A., Rafael B., Moisés, Renato, Sue Ellen, pelas
horas de estudos e disposição em ajudar aos demais colegas, agradeço. Agradeço,
também aos amigos do grupo de estudos em Sistemas Dinâmicos, Prof. Tiago Pereira, 
Élcio, Fernando e Jihard, pelos conhecimentos que compartilhamos.

Aos amigos e professores de graduação do UNASP-SP. Em particular, ao Prof. Ivanildo
Prado, por sempre animar os seus estudantes na continuação dos estudos. Dentre os meus
amigos da graduação, e que também são amigos pessoais, ao Décio, pelo seu companheirismo,
por sempre me receber em sua casa quando eu precisei, agradeço.

À minha família, em especial aos meus pais, José e Maria, por seu apoio moral, 
pelo suporte financeiro e por entender que a distância de casa, muitas
vezes é necessária para alcançarmos os nosso sonhos, agradeço. 

E por fim, porém, não menos importante, agradeço à UFABC, pelo apoio financeiro
durante todo o curso de mestrado e por possibilitar
a oportunidade dessa conquista.
\chapter*{Resumo}

Investigamos o surgimento do comportamento coletivo global em redes
de osciladores idênticos difusivamente acoplados, que no modelo
estabelecido é uma variedade invariante pelas equações do movimento. 
A interação é modelada através da teoria
de grafos e de sistemas dinâmicos. Utilizamos a teoria de contrações uniformes
em equações diferenciais lineares não-autônomas para estabelecer os critérios sobre
o parâmetro global de acoplamento que por sua vez define o estado síncrono e sua
respectiva estabilidade sob perturbações lineares e não-lineares. 
O parâmetro crítico global de interação é dado em função
somente da dinâmica individual dos osciladores, de propriedades espectrais
da função de acoplamento e do segundo autovalor do laplaciano da rede.

\noindent \textbf{Palavras-chave:} Sincronização, Redes, Estabilidade.

\chapter*{Abstract}
We investigate emergence of the global collective behavior in networks of
diffusively coupled identical oscillators, which in the established model is an 
invariant manifold of the motion equations. The interaction is modeled with the 
graph theory and dynamical systems theory. We use the uniform contractions
theory in non-autonomous linear differential equations to address the criterion under
the global coupling parameter, which it turns defines the synchronized motion and it
stability under small linear and non-linear perturbations. The critical 
global interaction parameter
is given only by the isolated dynamics, by the spectral properties of the coupling
function and the second eigenvalue network laplacian.

\noindent \textbf{Keywords:} Synchronization, Complex networks, Stability.

\tableofcontents    


\chapter{Lista de Símbolos}
\begin{tabular}{ll}
$\approx$ & Aproximadamente\\
$\operatorname{Mat}(\mathbb{R},n)$ & Conjunto das matrizes quadradas de dimensão $n$ com entradas reais\\
$\mathbb{R}_+ $ & Conjunto $[0,+\infty)$\\
$\delta_{ij}$ & Delta de Kronecker\\
$\mbox{det}(\cdot)$ & Determinante\\
$\setminus$ & Diferença entre conjuntos\\
$\exists$ & Existe\\
$\exp(\cdot)$ & Exponencial (o mesmo que número de Euler)\\
$\blacksquare$ & Fim da prova\\
$\sim$ & $i\sim j$ ($i$ é vizinho de $j$)\\
$\mbox{diag}(\cdot)$ & Matriz diagonal\\
$Df(x)$ & Matriz Jacobiana do campo de vetores $f$ calculada em $x$\\
$\max$ & Máximo\\
$\Vert\cdot\Vert$ & Norma \\
$e$ & Número de Euler ($e=2,718281828...$)\\
$\circ$ & Operação composição\\
$T(t,s)$ & Operador de evolução\\
$\nabla$ & Operador gradiente \\
$\alpha_c$ & Parâmetro crítico de acoplamento\\
$\alpha$ & Parâmetro global de acoplamento \\
$\colon{=}$ & Por definição igual a\\
$\otimes$ & Produto de Kronecker\\
$\left \langle \cdot,\cdot\right \rangle$ & Produto interno \\
$\lambda_2$ & Segundo autovalor do laplaciano do grafo\\
$\mathcal{O}(\cdot)$ & Símbolo de Landau\\
$\oplus$ & Soma direta\\
$\sup$ & Supremo\\
$\dagger$ & Transposto conjugado\\
$\vert\cdot\vert$ & Valor absoluto\\
$\mathbf{0}\in\mathbb{R}^n$ & Vetor em $\mathbb{R}^n$ com todas as entradas iguais a $0\in\mathbb{R}$\\
\end{tabular}

\listoffigures            

\mainmatter

\fancyhead[RE,LO]{\thesection}

\singlespacing              

\chapter{Introdução}
\label{cap:introducao}

O termo sincronização, que do grego significa ``ocorrência ao mesmo tempo'',
está relacionado a uma gama de fenômenos presentes em muitos ramos
das ciências naturais, engenharia e vida social \cite{amj01},\cite{strogatz03},\cite{tbk}.
Este fenômeno está enraizado, por exemplo, na vida humana, ocorrendo desde
processos metabólicos em nossas células às mais altas tarefas cognitivas \cite{kurths}.
Um dos primeiros indivíduos a estudar esse tipo de fenômeno foi o físico
alemão  Christiaan Huygens (1629 - 1694). Ele descobriu que relógios de pêndulo, quando
pendurados em um mesmo suporte, tendem à sincronização (movimento dos pêndulos).
De forma geral, esse estado síncrono surge da competição e colaboração entre os elementos 
de uma rede. Exemplos vão desde a sincronização em relógios de pêndulos à doenças neurais como 
Mal de Parkinson \cite{rosenblum98} e Epilepsia, onde
está última acontece quando um grupo específico de neurônios mantem um estado síncrono \cite{jung10}.

Modelamos então essa interação entre elementos através de uma rede, em termos matemáticos,
através da teoria dos grafos e da teoria de sistemas dinâmicos. Passamos a distinguir a estrutura 
da rede, a natureza de interação entre os elementos e o comportamento dinâmico individual dos
osciladores (elementos da rede). Uma rede é dita complexa quando a mesma não possui
uma estrutura regular de conectividade, e no mundo real a análise de sistemas interagentes
é feita a partir da modelagem em redes com estrutura complexa.
Existem vários tipos de redes que se enquadram
nesta característica, como por exemplo as redes aleatórias, redes com estrutura de pequeno
mundo, redes \textit{scale-free}, entre outras \cite{chung2006complex}.

Estaremos interessados em sincronização que ocorre de forma global. E
global no que se refere a totalidade dos elementos da rede e ao tempo, ou seja,
uma vez que o estado síncrono ocorre, o mesmo permanece por todo o tempo futuro. 
Um dos principais ingredientes para o surgimento da sincronização global na rede é que os
sistemas interagentes sejam idênticos, ou seja, que todos eles sejam descritos exatamente
pela mesma equação de movimento \cite{rosenblum03}.
Estaremos interessados apenas em redes de sistemas idênticos.


Recentemente, Pecora e Carroll \cite{pc98} estudando modelos de osciladores
em redes complexas com interação semelhante à difusão, demonstraram que
tais redes de sistemas dinâmicos não-lineares, porém idênticos, podem sincronizar globalmente, 
mesmo manifestando uma dinâmica individual complicada. Utilizamos
um modelo de rede semelhante ao utilizado por Pecora e Carroll em \cite{pc98}.

Considere que um grafo
$G$ modela uma rede de $n$ vértices. Sobre cada vértice da 
grafo introduzimos uma cópia da equação diferencial $\dot{x}=f(x)$, 
$f:\mathbb{R}^m\rightarrow\mathbb{R}^m$, onde
requeremos que $f\in C^r$, $r\geq2$, de forma que,	
essencialmente, $f$ representa uma dinâmica não-linear, possivelmente caótica,
e que as condições de existência e unicidade sejam satisfeitas. 
Para o acoplamento entre os osciladores, consideramos
o caso de interação difusiva, isto é, o acoplamento depende da diferença de uma função
dos estados dos vértices que estão conectados. Então, dado um vértice $i$ em $G$, a sua dinâmica passa a ser descrita
na forma
\begin{equation}\label{eq:modelointro}
\dot{x}_i  = f(x_i) + \alpha \sum_{j=1}^n{A_{ij}[H(x_j)-H(x_i)]},
\end{equation}
onde os $A_{ij}$'s são as entradas da matriz de adjacência, que codifica a informação
topológica da rede dizendo quais elementos possuem interação direta, 
$\alpha$ é o parâmetro global de acoplamento e $H:\mathbb{R}^m\rightarrow\mathbb{R}^m$ 
é uma função de acoplamento, que sem perder a generalidade podemos considerar que $H$,
é uma matriz. Trabalharemos apenas com o caso em que $H$ é uma matriz positiva-definida.

Note que um importante detalhe deste modelo é que quando $x_i = x_j$ para todo $i,j$, então
o termo de acoplamento desaparece identicamente de forma que $x_i(t) = x_j(t)$ ficam iguas
por todo o tempo, caracterizando assim um subespaço invariante. Passamos então
a chamar o estado síncrono de variedade de sincronização. A variedade de sincronização
é uma variedade invariante pelas equações do movimento e a análise da
sincronização global mostra que tal movimento síncrono e sua consequente
estabilidade é definida pela intensidade do acoplamento $\alpha$. Objetivamos então
caracterizar a estabilidade desse subespaço, que possui uma geometria trivial.

Nos últimos 20 anos houve uma evolução significativa no campo da sincronização em redes.
A análise da sincronização beneficia-se da análise estrutural de
tais redes \cite{bullmore09}, isto é, a estrutura
de interação entre os elementos influência significativamente
na sincronização \cite{pc98},\cite{wu}. Por exemplo, sabe-se que redes com estrutura de pequeno mundo, ou seja, redes
bem conectadas, possuem uma maior propensão à sincronização do que as redes regulares que
as geram, e consequente menos propensas à sincronização
do que as redes complexas aleatórias advindas da construção da rede pequeno mundo 
\cite{wattsandstrogatz}. Por outro lado,
a heterogeneidade na rede, isto é, redes em que a maioria dos elementos possuem poucos
vizinhos e apenas alguns elementos com muitos vizinhos, dificulta a sincronização
global da rede \cite{motter05}, no entanto facilita a sincronização entre elementos 
que possuem uma grande quantidade de vizinhos \cite{tiago}.

O modelo \eqref{eq:modelointro} tem sido muito usado como ponto de partida para se entender a emergência da coletividade.
Nijmeijer e Pogromsky em \cite{henk} estudaram o modelo do ponto de vista de controle e eles 
provam que as soluções
existem globalmente. Hasler e colaboradores em \cite{hasler} estudam
o problema construindo funções de Lyapunov para provar a estabilidade 
da variedade de sincronização. Macau e colaboradores em \cite{elbert} também estudam esse
problema utilizando a existência funções de Lyapunov para garantir a estabilidade 
do comportamento síncrono.
Josic em \cite{josic2} estudou esse problema do ponto de vista das variedades
invariantes. Pereira em \cite{lecturenotes} estudou este problema do ponto de vista
de dicotomias e nossa abordagem é baseada na mesma perspectiva.
%
%

%

\section{Principais Resultados}
\label{sec:contribucoes}

Os estudos aqui desenvolvidos  provém uma compreensão rigorosa do comportamento
coletivo em redes e assim abre frentes para novos resultados na área.

Uma abordagem muito utilizada para se estudar a estabilidade da variedade de sincronização 
em redes é baseada
na teoria dos expoentes de Lyapunov \cite{pc98}, \cite{pecora09}, \cite{tiago}. 
No entanto, sabe-se que
neste contexto, perturbações arbitrariamente pequenas podem causar mudanças abruptas na
estabilidade do movimento síncrono \cite{barreira2002lyapunov}.
Abordamos então o fenômeno da sincronização através de uma nova perspectiva, isto é,
através da teoria de contrações uniformes em equações diferenciais lineares não-autônomas, de forma que
poderemos garantir que a variedade de sincronização seja resistente à pequenas perturbações.
Destacamos as principais contribuições abordadas neste trabalho:
\\

1. Existência global das soluções de \eqref{eq:modelointro}:
Construímos uma função de Lyapunov para a dinâmica isolada e extendemos para a dinâmica coletiva. 
Dessa forma, independentemente
da rede considerada, as soluções de \eqref{eq:modelointro} são limitas e portanto 
existem globalmente. Veja o Teorema \ref{teo:egscaposc}.
\\

2. Estabilidade da variedade de sincronização:
A partir da teoria de contrações uniformes em equações diferenciais lineares não-autônomas 
poderemos determinar um parâmetro crítico de acoplamento $\alpha_c$ tal que
para $\alpha>\alpha_c$ em \eqref{eq:modelointro} a variedade de sincronização é estável.
Precisamente, temos o Teorema \ref{teo:principal1}.
\\

3. Persistência da sincronização:
Introduzimos um modelo de perturbação que age
na função de acoplamento de \eqref{eq:modelointro} e utilizamos um resultado de persistência
de \cite{coppel} para estabelecer a magnitude das perturbações que não destrói a estabilidade
da variedade de sincronização. 
Veja o Teorema \ref{teo:persistenciaglobal}.
\\

Agora, gostaríamos de expor algumas considerações sobre os itens acima.
No item 1 contruíremos uma função de Lyapunov para a coletividade supondo
que existe uma função de Lyapunov quadrática para a dinâmica isolada.
O item 2 é um resultado fundamental e
vários autores já o provaram \cite{josic2}, \cite{henk}, \cite{lecturenotes}.
A novidade então, é que abordamos na prova a teoria de contrações uniformes 
em equações diferenciais lineares
não-autônomas \cite{coppel}. Essa teoria não é usualmente empregada para estudar redes,
porém foi recentemente usada para provar o item 2 em \cite{lecturenotes}. O que fazemos
aqui é expandir e abordar com mais detalhes a prova. Portanto,
não há nenhuma técnica nova mas no contexto de sincronização o resultado é novo.

O item 3 é o ponto novo e
segue praticamente como um corolário do resultado apresentado no
item 2. O caso onde $f$ em \eqref{eq:modelointro} não é igual foi tratado
por Pereira e colaboradores em \cite{tiagocoherence}, porém o caso em que
$H$ não é igual é novo na literatura. No caso da pertuação na função de acoplamento,
a estrutura da rede tem um papel muito importante para a conservação da estabilidade.

\section{Organização do Trabalho}
\label{sec:organizacao_trabalho}

Queremos deixar claro ao leitor que, do Capítulo \ref{cap:conceitos} ao Capítulo 
\ref{cap:grafos} temos por objetivo fundamentar toda
a teoria matemática que será utilizada para modelar o nosso problema, a saber, o problema
de estudar a o surgimento da sincronização e sua estabilidade em redes. Portanto, o leitor deve
saber que esse é um tema interdisciplinar e que por essa característica, faz-se necessário
fundamentar e então interligar algumas diferentes áreas da matemática, para então dar o
devido tratamento rigoroso sobre o mesmo. Isso faz com que, a princípio,
pareça que estamos tratando de temas que estão a parte do título apresentado pela dissertação,
porém os mesmos são de fundamental importância para o tratamento do problema proposto.

No Capítulo \ref{cap:conceitos}, apresentamos os conceitos preliminares, e as notações matemáticas
que utilizaremos durante o desenvolvimento do trabalho.
No Capítulo \ref{cap:dinamicanaolinear} abordaremos alguns tópicos da teoria qualitativa
de equações diferenciais ordinárias, à qual 
é usada para modelar e entender o
comportamento individual dos elementos da rede. 
Abordaremos 
em especial
a construção da função de Lyapunov. Em seguida, no Capítulo \ref{cap:edolinear} trabalharemos 
sobre a estabilidade e suas
persistências, da solução trivial em equações diferenciais lineares não-autônomas.
No Capítulo \ref{cap:grafos} abordaremos os fundamentos da teoria espectral dos grafos, que por sua vez
é utilizada para modelar topologicamente uma rede qualquer. Seguido pelo Capítulo
\ref{cap:osciladores} o qual introduziremos formalmente o modelo de acoplamento difusivo em redes,
e enunciamos os principais resultados da presente dissertação.
No Capítulo \ref{cap:doisosciladoresacoplados} trabalharemos vários exemplos para a teoria abordada,
principalmente com o exemplo mais trivial de rede, a saber, 
o de dois osciladores difusivamente acoplados, do qual extrairemos valorosos resultados. E
finalmente, no Capítulo \ref{cap:estabilidadedassolucoessincronas} apresentaremos as provas dos
principais resultados enunciados no Capítulo \ref{cap:osciladores}.
Por fim, no Capítulo \ref{cap:conclusoes} discutiremos algumas conclusões e 
abordaremos ideias sobre uma possível continuação do respectivo trabalho. 
Nos Apêndices \ref{ape:algebralinearA} e \ref{ape:edoB}, trataremos de resultados gerais
e específicos sobre respectivamente Álgebra Linear e Equações Diferenciais Ordinárias, os quais
são utilizados ao longo da dissertação.

\chapter{Conceitos Preliminares}
\label{cap:conceitos}

O presente capítulo é utilizado, principalmente, para
introduzir a notação  
que servirá como base para todo o texto.  Estão disponíveis
nos Apêndices \ref{ape:algebralinearA} e \ref{ape:edoB} outras informações
basilares para o entendimento da problemática exposta. Trabalharemos essencialmente
com espaços que possuem norma e produto interno e matrizes quadradas.

\section{Espaços Normados e Espaços com Produto Interno}
\label{sec:fundamentosdealgelin}

Considere $\R$ o corpo dos número reais.
Indicamos por $\mathbb{R}_+$ o conjunto $[0,+\infty)$.
Seja $E$ um espaço vetorial sobre $\R$. Uma \textbf{norma} em $E$ é a função
$\Vert \cdot \Vert : E \rightarrow \mathbb{R}_+$ satisfazendo as
seguintes propriedades:

\begin{enumerate}
\item[(i)] Para todo $x \in E \setminus \{\mathbf{0}\}$ tem-se que $\Vert x \Vert >0$.
\item[(ii)] $\Vert x \Vert = 0$ se e somente se $x=\mathbf{0}$ (vetor nulo).
\item[(iii)] Dados $\alpha \in \mathbb{R}$ e $x \in E$, tem-se que 
$\Vert \alpha x\Vert = \vert\alpha\vert\Vert x\Vert$, onde $\vert \alpha\vert$ é o valor
absoluto de $\alpha$.
\item[(iv)] ({Desigualdade Triangular)} Para todo $x,y \in E$ vale que
$\Vert x+y\Vert\leq\Vert x\Vert +\Vert y\Vert$.
\end{enumerate}

Dado $x=(x_1,\cdots,x_n)\in\mathbb{R}^n$, importantes normas que surgem são
\begin{enumerate}
\item ($p$-norma) $\Vert x\Vert_p = \left(\sum_{i=1}^n \vert x_i\vert^p\right)^{1/p}$; Quando
$p=2$ temos a norma Euclidiana.
\item (Norma do máximo) $\Vert x\Vert_\infty = \max_i{\vert x_i\vert}$;
\item Norma da soma ($p=1$).
\end{enumerate}

Dadas duas normas $\Vert \cdot\Vert_u$ e $\Vert\cdot\Vert_v$ em $E$. Dizemos que
tais normas são equivalentes se existem números positivos $c_1$ e $c_2$ tais que
$$
c_1\Vert x\Vert_u\leq\Vert x\Vert_v\leq c_2\Vert x\Vert_u,
$$
para todo $x \in E$.

Devido a forte garantia do teorema a seguir, por todo o presente texto, não fazemos
uso de uma norma específica na maioria das situações em que trabalhamos com tal conceito.

\begin{teorema}
Duas normas quaisquer em $\mathbb{R}^n$ são equivalentes.
\end{teorema}

\noindent
\textbf{Prova:} Ver \cite{lima1987curso}, p. 19.

De forma mais geral, em espaços de dimensão finita, todas as normas são equivalentes 
\cite{o2006metric}.

Uma norma em um espaço vetorial $E$ dá origem à noção de distância em $E$. Para todo
$x,y\in E$, a {distância} de $x$ a $y$ é definida por
$$
d(x,y) = \Vert x-y\Vert.
$$
Dessa forma, como em um espaço vetorial normado tem-se a noção de distância, ou métrica, 
o mesmo pode ser chamado de \textbf{espaço métrico}.
%
Um espaço métrico $E$ é dito ser \textbf{completo} se toda sequência de Cauchy em $E$
é convergente em $E$ \cite{o2006metric}. 

\begin{definicao}\label{def:uniflim}
Seja $E$ um espaço métrico completo e $f:U\subset\mathbb{R} 
\rightarrow E$ uma função contínua. Dizemos que $f$ é \textbf{uniformemente limitada} em $t$, 
$t\in U$, se existe alguma constante $c$ tal que
$$
\sup_{t\in U} \Vert f(t)\Vert =c.
$$
\end{definicao}

Seja $E$ um espaço vetorial sobre o corpo $\R$. 
O \textbf{produto interno} é a função real $\langle \cdot,\cdot\rangle:E\times E\rightarrow\mathbb{R}$, que associa
cada par de vetores $u,v\in E$ ao número $\langle u,v \rangle$, de modo que,
são válidas as propriedades de simetria, bilinearidade e positividade \cite{lima1987curso}.






\section{Matrizes}

Indicamos por $\Mat(\R,n)$ o conjunto das matrizes quadradas de dimensão $n\times n$ com entradas 
em $\R$.
\begin{definicao}
Dizemos que uma matriz $H \in Mat(\R,n)$ é \textbf{simétrica} se $H=H^\dagger$, onde $\dagger$
representa o transposto. 
\end{definicao}
A matriz identidade,
indicada por $I_n\in\Mat(\mathbb{R},n)$, é uma matriz simétrica.
Qualquer matriz $A\in\Mat(\R,n)$ pode ser decomposta na forma
$$
A=\frac{A+A^\dagger}{2}+\frac{A-A^\dagger}{2},
$$
onde a parcela ${(A+A^\dagger)}/{2}$ é chamada de parte simétrica de $A$.

Reservamos a notação ``$\mbox{diag}$'' para representar as matrizes diagonais em $\Mat(\R,n)$. Ou seja,
as matrizes cujos os valores fora da diagonal são nulos. Por exemplo, $I_3 = \mbox{diag}(1,1,1)$.

\begin{definicao}
Seja $x\in\mathbb{R}^n\setminus\{\mathbf{0}\}$. Dizemos que uma matriz 
$A\in\Mat(\mathbb{R},n)$, simétrica, é
positiva-definida $[$positiva semi-definida$]$ se $\langle x,Ax \rangle >0$ 
$[ \langle x,Ax \rangle \geq0]$ e negativa-definida $[$negativa semi-definida$]$ se 
$\langle x,Ax \rangle <0$ $[\langle x,Ax \rangle \leq0]$.
\end{definicao}

Note que $\langle x,Ax\rangle = x^\dagger Ax$, de forma que em muitos contextos poderemos usar
essa equivalência.

\subsection{Teoria da Representação}




\begin{definicao}
Diremos que uma matriz $Q\in\Mat(\R,n)$ é ortogonal se $Q^\dagger =Q^{-1}$, onde $Q^{-1}$ é
a inversa de $Q$.
\end{definicao}

\begin{definicao}
Uma matriz real $A$ é ortogonalmente diagonalizável se existe uma matriz ortogonal $Q$ tal que
$Q^\dagger A Q =D$ é uma matriz diagonal.
\end{definicao}

\begin{teorema}[Teorema Espectral]\label{teo:teoremaespectral1}
Uma matriz $A\in\Mat(\mathbb{R},n)$ é ortogonalmente diagonalizável se e 
somente se é simétrica.
\end{teorema}

\noindent
\textbf{Prova:} Ver \cite{lima2008algebra}, página 167.




\begin{corolario}\label{cor:teoremaespectral}
Seja A uma matriz real e simétrica. Então existe uma matriz ortogonal Q tal que
\begin{enumerate}
 \item $D = Q^\dagger A Q$ é uma matriz diagonal onde os elementos da diagonal são
os autovalores de $A$ e
\item As colunas de $Q$ são formadas pelos autovetores de $A$.
\end{enumerate}
\end{corolario}

\begin{proposicao}
Seja $A \in\Mat(\R,n)$ uma matriz simétrica. Então para todo $\lambda$, autovalor de
$A$, tem-se que
$A$ é positiva-definida $[$positiva semi-definida$]$ se e somente se $\lambda >0$ 
$[\lambda \geq0]$ e negativa-definida $[$negativa semi-definida$]$ se e somente se 
$\lambda <0$ $[\lambda \leq0]$.

\end{proposicao}

Como $A$ é simétrica então podemos formar uma base de $\R^n$ com os autovetores de $A$ os quais
podem ser tomados ortonormalizados. Seja $\beta = \{v_1,\cdots,v_n\}$ tal base, então
$x=\sum_{i=1}^n c_i v_i$, $x\in \R^n$, o que produz
\begin{align*}
\langle x,Ax \rangle & = \left\langle \sum_i c_i v_i,A \sum_i c_i v_i\right\rangle\\
& = \left\langle \sum_i c_i v_i, \sum_i c_i \lambda_i v_i \right\rangle
= \sum_i {c_i}^2 \lambda_i.
\end{align*}
Sendo $\lambda_i>0$ para todo $i$ tem-se que $\langle x,Ax \rangle>0$. Os demais resultados
seguem alterando-se as hipóteses sobre o sinal de cada $\lambda_i$.


\subsection{Norma Induzida de Matriz}\label{seccao:normamatriz}

Se $\Vert\cdot\Vert$ é uma norma em $\R^n$ e $A\in\Mat(\R,n)$, definimos a norma de
$A$ por 
$$
\Vert A\Vert = \sup_{x\in \R^n, \Vert x\Vert\leq1}{\Vert Ax\Vert}.
$$
Essa norma é conhecida como norma induzida \cite{lancaster}. 
Por toda a dissertação, fazemos uso dessa norma
ao se trabalhar com matrizes. Uma notação mais precisa seria $\Vert A\Vert_p$, onde o índice $p$ indica
qual a norma sobre $\R^n$ estamos induzindo a norma em $\Mat(\R,n)$. Omitiremos tal índice, 
indicando-o apenas quando for necessário. 

Observação: a noção de função uniformemente limitada dada pela Definição \ref{def:uniflim} pode
também ser extendida à matrizes considerando a norma induzida.



Uma das propriedades mais importantes dessa norma é a propriedade sub-multiplicativa. Mais
especificamente temos o
\begin{lema}\label{lema:matrixpropeties}
Dados $A,B\in\Mat(\R,n)$ temos
\begin{enumerate}
 \item $\Vert Ax\Vert \leq \Vert A\Vert \Vert x\Vert$, para cada $x\in \R^n$;
\item $\Vert AB\Vert \leq \Vert A\Vert \Vert B \Vert$ (propriedade sub-multiplicativa).
\end{enumerate}
\end{lema}

\noindent
\textbf{Prova:} Ver \cite{claus}, p. 351.


\begin{exemplo}
Sejam $A\in\Mat(\R,n)$, $A=[A_{ij}]$, e $x\in \R^n$. Considere a norma do máximo
definida em $\R^n$ por $\Vert x\Vert_{\infty} = \max_i \vert x_i\vert$, onde
$x=(x_1,\cdots,x_n)$. Então
\begin{align*}
\Vert A \Vert_{\infty} & = \sup_{\Vert x\Vert_{\infty}=1} \Vert Ax\Vert_\infty\\
& = \max_i \left\vert \sum_j A_{ij}x_j \right\vert\\
& \leq \max_i \sum_j \vert A_{ij}\vert \vert x_j\vert
\leq \Vert x\Vert_\infty \max_i \sum_j\vert A_{ij}\vert = \max_i \sum_j\vert A_{ij}\vert.
\end{align*}
Por outro lado, podemos mostrar também que $\Vert A\Vert_\infty\geq \max_i \sum_j\vert A_{ij}\vert$.
Certamente, para algum $1\leq k\leq n$ tem-se $\max_i \sum_j\vert A_{ij}\vert = \sum_j \vert A_{kj}\vert$.
Então, existe $x_0=(\alpha_1,\cdots,\alpha_n)\in C^n$ com $\Vert x_0\Vert_\infty =1$ escrito da forma
$$
\alpha_j = \left\{ \begin{array}{rl}
 \frac{\vert A_{kj}\vert}{A_{kj}} &\mbox{ se $A_{kj}\neq0$} \\
  0 &\mbox{ se $A_{kj}=0$ }
       \end{array} \right.
$$
que cumpre
\begin{align*}
\Vert Ax_0\Vert_\infty & = \max_i \left\vert \sum_j A_{ij}\alpha_j\right\vert\\
& = \max_i \left\vert \sum_j A_{ij} \frac{\vert A_{kj}\vert}{A_{kj}}\right\vert
 = \sum_j \vert A_{kj}\vert
\geq \max_i \sum_j \vert A_{ij}\vert.
\end{align*}
Portanto,
$$
\Vert A\Vert_\infty = \max_i \sum_j \vert A_{ij}\vert,
$$
que simplesmente é a soma máxima dos valores absolutos das linhas de $A$.
\end{exemplo}


\subsection{Produto de Kronecker}

Sejam $A=[A_{ij}]\in\Mat(\R, n)$ e
$B=[B_{ij}]\in\Mat(\R,m)$, então o produto de Kronecker de $A$ e $B$, escrito como
$A\otimes B\in\Mat(\R,n\times m)$, é a matriz formada pelas submatrizes $A_{ij}B$ escrita
da forma
$$
A\otimes B =\left[
\begin{matrix}
A_{11}B & A_{12}B & \cdots & A_{1n}B\\
A_{21}B & A_{22}B  & \cdots & A_{2n}B\\
\vdots & \vdots & \quad & \vdots\\
A_{n1}B & A_{m2}B & \cdots & A_{nn}B
\end{matrix}\right].
$$



\begin{teorema}\label{teo:kroneckerdis}
Sejam $A\in \Mat(\R,n)$, $B\in \Mat(\R,m)$, $U\in \Mat(\R,n)$ e
$V\in \Mat(\R,m)$. Então
$$
(A\otimes B) (U\otimes V) = AU\otimes BV.
$$
\end{teorema}

\noindent
\textbf{Prova:} Ver \cite{lancaster}, p. 408.

\begin{teorema}\label{teo:kroneckerdagger}
Sejam $A\in \Mat(\R,n)$ e $B\in \Mat(\R,m)$, então
$$
(A\otimes B)^\dagger = A^\dagger \otimes B^\dagger.
$$
\end{teorema}

\noindent
\textbf{Prova:} Por definição tem-se $(A\otimes B)^\dagger = (A_{ij}[B])^\dagger = A_{ji}[B^\dagger]
= A^\dagger \otimes B^\dagger$. $\blacksquare$

\begin{teorema}\label{teo:positivedifinitekronecker}
Sejam $A\in\Mat(\mathbb{R},n)$ e $B\in\Mat(\mathbb{R},m)$.
\begin{enumerate}
\item Se $A$ e $B$ são positiva-definidas então $A\otimes B$ é positiva-definida.
\item Se $A$ é positiva-definida e $B$ é positiva semi-definida então $A\otimes B$ é positiva semi-definida.
\item Se $A$ é positiva semi-definida e $B$ é positiva definida então $A\otimes B$ é positiva semi-definida.
\end{enumerate}
\end{teorema}

\noindent
\textbf{Prova:} Considere o item \textit{1}. 
Sejam $x\in\mathbb{R}^n$ e $y\in\mathbb{R}^m$ não-nulos. Então por hipótese
tem-se que $x^\dagger Ax>0$ e $y^\dagger By>0$. Assim  
$x^\dagger Ax \otimes y^\dagger By = (x \otimes y)^\dagger (A\otimes B) (x\otimes y)>0$ (Teoremas
\ref{teo:kroneckerdis} e \ref{teo:kroneckerdagger}). 
Os demais itens seguem de forma análoga. $\blacksquare$

\begin{teorema}
Sejam $x\in\mathbb{R}^n$, $y\in\mathbb{R}^m$ e $\Vert\cdot\Vert_p$ uma $p$-norma. Então
$$
\Vert x\otimes y\Vert_p = \Vert x\Vert_p \Vert y\Vert_p.
$$
\end{teorema}
\noindent
\textbf{Prova:} $\Vert x\otimes y\Vert_p^p = \sum_{k=1}^n\sum_{l=1}^m\vert x_k y_l\vert^p =
\sum_{k=1}^n\vert x_k\vert^p \sum_{l=1}^m\vert y_l\vert^p = \Vert x\Vert_p^p \Vert y\Vert_p^p$.
$\blacksquare$

O resultado é estensível para matrizes condiderando a norma induzida.

\section{Conjuntos}

A \textit{bola aberta} de centro em $x_0\in\mathbb{R}^n$ e raio $\delta>0$ é o conjunto
$
B_{\delta}(x_0) = \{ x\in\mathbb{R}^n : d(x,x_0) <\delta\}.
$
Analogamente, a \textit{bola fechada} de centro $x_0$ e raio $\delta> 0$ é o conjunto
$
\overline{B}_{\delta}(x_0) = \{ x\in\mathbb{R}^n : d(x,x_0) \leq\delta\}.
$
\begin{definicao}
Dizemos que um conjunto $E\subset\mathbb{R}^n$ é \textbf{limitado} quando o mesmo está inteiramente 
contido em alguma bola fechada.
\end{definicao}

Dado um ponto $a\in E$, dizemos que $a$ é \textit{ponto interior} a $E$ quando, para algum $\delta>0$,
tem-se que $B_{\delta}(a)\subset E$. 
\begin{definicao}
Um  conjunto $E$ é \textbf{aberto} quando todo ponto de $E$ é
ponto interior a $E$.
\end{definicao}


\begin{definicao}
Em relação à $\mathbb{R}^n$, dizemos que um conjunto $E$ é \textbf{fechado} se e somente se
o seu complementar, $\mathbb{R}^n\setminus E$ é aberto.
\end{definicao}

\begin{definicao}
Um conjunto $E\subset\mathbb{R}^n$ chama-se \textbf{compacto} quando o mesmo é limitado e fechado.
\end{definicao}

\begin{teorema}[Weierstrass]
Seja $E$ um conjunto compacto de $\mathbb{R}^n$. Então, toda função contínua $f:E\rightarrow\mathbb{R}$
assume um valor máximo e um valor mínimo. Ou seja, existem $x_0, x_1\in E$ tais que
$f(x_0)\leq f(x)\leq f(x_1)$ para todo $x\in E$. 
\label{weierstrass}
\end{teorema}

\noindent
\textbf{Prova:} Ver \cite{lima1987curso}, p. 45. O Teorema de Weierstrass também é válido
para espaços compactos quaisquer \cite{munkres2000topology}.


\chapter{EDO's - Tópicos em Teoria Qualitativa}\label{cap:dinamicanaolinear}

Para estudar sincronização e sua consequente estabilidade, é necessário impor que os elementos
da rede estudada possuam dinâmica. Em particular, essa dinâmica pode ser uma dinâmica
não-linear (o modelo de rede que iremos estudar será formalmente apresentado no Capítulo
\ref{cap:osciladores}). 

Encontrar as soluções de uma equação diferencial não-linear pode não ser uma tarefa fácil. Na
maioria dos casos (principalmente com sistemas em altas dimensões) pode não ser possível
de se computar analiticamente as soluções \cite{brauer}.
Na verdade, existem apenas alguns poucos casos em que somos capazes de encontrar analiticamente
a solução [soluções] de uma equação diferencial [sistema de equações diferenciais], como por
exemplo em equações diferenciais lineares cujos coeficientes são constantes. Portanto, neste capítulo
trabalharemos alguns pontos da teoria qualitativa de equações diferencias ordinárias (EDO's).

Seja $D$ um aberto de $\mathbb{R}^m$,
$f: D \rightarrow \mathbb{R}^m$, $m\geq1$, um campo vetorial de classe $C^r$, $r\geq2$, e considere a 
equação diferencial
\begin{equation}\label{x=ftx}
\dot{x} = f(x).
\end{equation}
Para um dado $x(s) \in D$, o {problema de valor inicial} é encontrar
um intervalo $I \ni s$, da reta real, onde $s$ é ponto interior e $x(t)$
está definida em $I$ satisfazendo a equação \eqref{x=ftx} e respeitando a condição
inicial $x(s)=x_s$. 
Como $f$ é de classe $C^r$, $r\geq2$, garante-se
o teorema de existência e unicidade das soluções (Teorema \ref{teo:picardlindelof}).
Na maioria dos casos, vamos considerar $I\subset\mathbb{R}_+$ e $s=0$.

Por vezes, nos referimos à equação diferencial \eqref{x=ftx} pela palavra sistema e $f$ pela
palavra dinâmica.
Muitos nomes são dados à solução de uma equação diferencial,
fazemos a escolha de usar a palavra trajetória ou simplesmente solução.




\begin{definicao}
Seja $\phi : I\times D\rightarrow\mathbb{R}^m$ a aplicação $\phi (t,x_0) = x(t)$. Dizemos que
$\phi$ é o \textbf{fluxo} do campo de vetores $f$ onde $x(t) : I\rightarrow\mathbb{R}^m$ é a trajetória
de $f$ por $x_0$ em cada ponto $t\in I$.
\end{definicao}
Podemos utilizar também a notação $\phi_t(x_0)$ para o fluxo. De forma mais simplificada, o fluxo é
uma aplicação que mapeia a condição inicial na solução. Escrevendo $\phi_t(x(s))=\phi_{t,s}(x)$
e considerando a operação composição, nota-se que o fluxo é um grupo de difeomorfismos a um
parâmetro \cite{arnold1992ordinary}, onde são válidas as propriedades
\label{pag:fluxo}
\begin{enumerate}
 \item $\phi_{t,s} \circ \phi_{s,u} (x) = \phi_{t,u} (x)$ para todo $s,u,t\in I$ e $x\in D$.
\item $\phi_{t,t} (x) = x$ para todo $t\in I$ e $x\in D$.
\end{enumerate}

\begin{exemplo}
Dado o sistema linear $\dot{x} = Ax$, onde $A \in\Mat(\mathbb{R},m)$, temos que para 
cada ponto $x_0\in\mathbb{R}^m$ dado,
$$
x(t) = \exp({tA})x_0
$$
é a solução de $\dot{x} = Ax$ com a condição inicial $x(0)=x_0$, de modo que o fluxo
$\phi : \mathbb{R}\times\mathbb{R}^m \rightarrow \mathbb{R}^m$ é dado por
$$
\phi (t,x_0) = \phi_t(x_0) = \exp({tA})x_0.
$$
Para $t=0$ o fluxo resume-se ao operador identidade e para cada $t$, $\phi_t = \exp({tA})$ 
é um isomorfismo linear cujo inverso é o fluxo $\phi_t = \exp({-tA})$.
\end{exemplo}

\section{Sistemas Dissipativos}

\begin{definicao}\label{def:invariantset}
Dizemos que um conjunto $\Omega\subset D$ é \textbf{invariante} pelo fluxo
$\phi$ do campo de vetores ${f}$ se $\phi_t(\Omega)\subseteq\Omega$
para todo $t\in\mathbb{R}$.
\end{definicao}

Equivalentemente, dizemos que $\Omega$ é \textbf{positivamente invariante} se
$\phi_t(\Omega)\subseteq\Omega$ para todo $t\geq 0$. E negativamente invariante
quando $\phi_t(\Omega)\subseteq\Omega$ para $t\leq0$.
Em outras palavras se $x(s)\in \Omega$ então $x(t)\in\Omega$ para todo $t\geq s>0$,
no caso de $\Omega$ ser positivamente invariante. De forma
mais intuitiva, dizemos que o conjunto $\Omega$ é positivamente invariante se 
as trajetórias entrarem em $\Omega$ lá permanecendo por todo o tempo futuro.
Assim, temos a

\begin{definicao}
Dizemos que o sistema \eqref{x=ftx} é \textbf{dissipativo} se suas soluções entram, a tempo finito,
em um conjunto positivamente invariante $\Omega\subset D$.
\end{definicao}


$\Omega$ será chamado de \textbf{domínio absorvente} do sistema $\dot{x}=f(x)$. A existência de um domínio
absorvente garante que as soluções do mesmo são limitadas e portanto existam globalmente
em decorrência do teorema de extensão (Teorema \ref{teo:extensao}). A seção
a seguir tratará de condições para garantir a existência do domínio absorvente.

\subsection{Segundo Método de Lyapunov}

Esta seção é de valorosa importância pois a técnica que traz a ideia
da função de Lyapunov será usada principalmente
para garantir a existência global das soluções e encontrar
regiões de estabilidade assintótica. 

Descoberta por
Aleksandr Mikhailovich Lyapunov (1857 - 1918) no final do século 19, a técnica é conhecida por segundo método
de Lyapunov ou método direto \cite{brauer}, pois a mesma pode ser aplicada diretamente
às equações diferenciais sem ter nenhum conhecimento das soluções. 
A função de Lyapunov é bastante utilizada para estudar a estabilidade de pontos fixos.
Porém, esse não será o nosso propósito no uso desta ferramenta. A estabilidade de pontos fixos, 
será estudada no Capítulo \ref{cap:edolinear}.
A dificuldade desta técnica é que não existe uma fórmula fechada
para encontrar tal função para o sistema de equações estudado.

Seja $\Omega$ um subconjunto conexo não-vazio de $\mathbb{R}^m$. Seja $V: \R^m\rightarrow\mathbb{R}$
uma função escalar continuamente diferenciável. Definimos
as noções de definitude\footnote{Classe dos conceitos de função positiva-definida
e negativa-definida} de $V$:

\begin{definicao}
Dizemos que a função escalar $V:\mathbb{R}^m\rightarrow\mathbb{R}$ é \textbf{positiva-definida}
com relação ao conjunto $\omega\subset\mathbb{R}^m$, contendo a origem, se $V(x)>0$ para todo
$x\in \mathbb{R}^m\setminus \omega$ e $V(\mathbf{0}) = 0$.
\end{definicao}

\begin{definicao}
Dizemos que $V:\mathbb{R}^m\rightarrow\mathbb{R}$ é \textbf{negativa-definida} em relação
ao conjunto $B\in\mathbb{R}^m$ se $-V$ é positiva-definida com respeito à esse conjunto.
\end{definicao}




A derivada de $V$ com respeito ao sistema \eqref{x=ftx} é o produto interno
\begin{equation}
V'(x) = \left \langle \nabla V(x),f(x) \right \rangle,
\end{equation}
onde $\nabla V(x) = \left( \frac{\partial V}{\partial x_1}(x),\cdots,
\frac{\partial V}{\partial x_m}(x)  \right)$.
De fato, visto que, se $x(t)$ é solução da equação \eqref{x=ftx}, então
pela regra da cadeia tem-se que
\begin{align*}
\frac{d V(x)}{dt} & = \left \langle \nabla V(x),\frac{dx(t)}{dt} \right \rangle \\
& = \left \langle \nabla V(x),f(x) \right \rangle.
\end{align*}

\begin{definicao}
Sejam $f:D\rightarrow \mathbb{R}^m$ um campo vetorial no aberto $D\subseteq\mathbb{R}^m$
e $V:\mathbb{R}^m\rightarrow\mathbb{R}$ uma função contínua. Dizemos que a função $V$
é uma \textbf{função de Lyapunov} para $f$ com respeito ao conjunto
$\Omega\subset D$ se
\begin{enumerate}
 \item[(i)] $V$ é positiva-definida em relação ao conjunto $\Omega$ e 
 \item[(ii)] $V'$ é negativa-definida em relação ao mesmo conjunto $\Omega$.
\end{enumerate}
\end{definicao}

\begin{definicao}
Dizemos que $V$
é radialmente ilimitada se
$$
\lim_{\Vert x\Vert\rightarrow\infty}{V(x)} = \infty.
$$
\end{definicao}

\begin{definicao}
Dizemos que a função escalar contínua $V:\mathbb{R}^m\rightarrow\mathbb{R}$ é própria se
dado um conjunto compacto $K\in\mathbb{R}$ então a pré-imagem $V^{-1}(K)$ é compacto de $\mathbb{R}^m$.
\end{definicao}

\begin{lema}
Uma função contínua $V:\mathbb{R}^m\rightarrow\mathbb{R}$ é própria se e somente se é radialmente
ilimitada.
\end{lema}

\noindent
\textbf{Prova:} Ver \cite{terrell2009stability}, p. 180.

\begin{corolario}
Se $V$ é radialmente ilimitada então os conjuntos de nível $V(x) = c$ são compactos.
\end{corolario}

\begin{teorema}[Lyapunov]\label{teo:lyapunov}
Seja $V$ uma função de Lyapunov para $f$ em \eqref{x=ftx} com respeito ao conjunto $\Omega\subset D$
contendo a origem, e além disso suponha que $V$ é radialmente ilimitada.
Sendo assim, as trajetórias da equação \eqref{x=ftx} entram, a tempo finito, em
$\Omega$, e lá permanecem por todo o tempo futuro.
\end{teorema}

Dizer que as trajetórias entram a tempo finito no conjunto absorvente $\Omega$ significa
dizer que o sistema \eqref{x=ftx} é dissipativo. 
%
Certamente o Teorema \ref{teo:lyapunov} depende da existência da função escalar $V$ com
as propriedades citadas. Porém tal resultado não nos mostra como construir a referida função, e esta
é a principal limitação deste método. Então, não existe um procedimento específico para a construção
de $V$, mas em boa parte dos casos, uma forma quadrática é uma boa candidata. Iremos então
considerar apenas as formas quadráticas para a função de Lyapunov e 
assumir que a mesma é dada considerando a

\begin{figure}[!ht]
  \centering
  \includegraphics[width=0.8\textwidth]{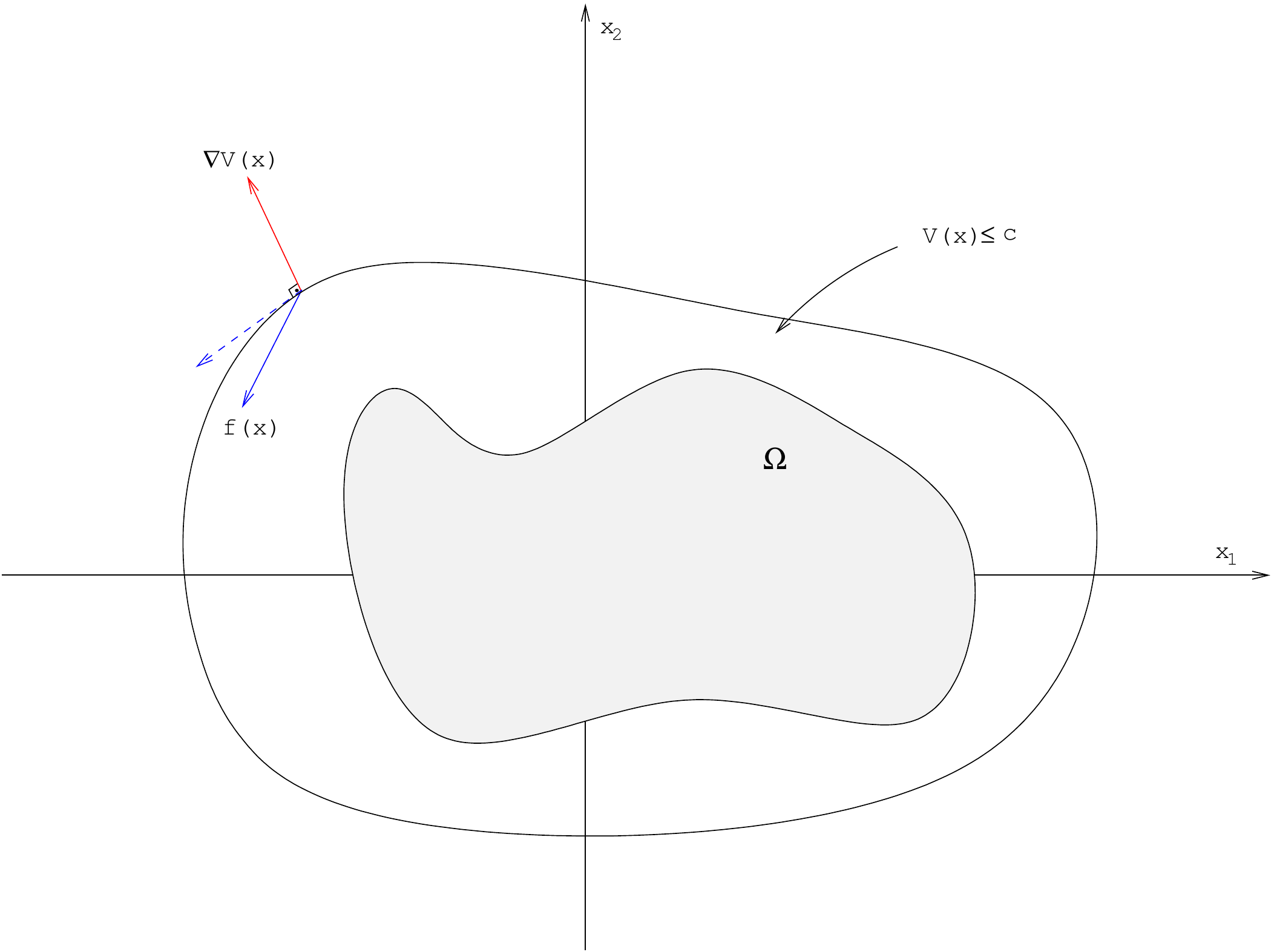} 
  \caption{Ilustração geométrica para o Teorema \ref{teo:lyapunov} com $m=2$.}
  \label{fig:lyapunov} 
\end{figure}

\begin{suposicao}\label{sup:atrator}
Existe uma matriz positiva-definida $Q$ tal que
\begin{equation}
V(x)= 1/2\langle x-a,Q(x-a) \rangle
\end{equation}
onde $a\in\mathbb{R}^m$ é fixo,  $V$ é positiva-definida em relação ao conjunto
\begin{equation}
\Omega \colon{=} \{x\in\mathbb{R}^m :  V(x)\leq \rho \},
\end{equation}
para algum $\rho\in\mathbb{R}_+$ e
$ V'(x)$ é negativa-definida em relação à $\Omega$.
\end{suposicao}
O papel de $a\in\mathbb{R}^m$ é realizar uma translação de eixos, se necessária.
A existência de um domínio absorvente para o sistema \eqref{x=ftx},
garante que as soluções são limitadas e existem globalmente, pois como construído,
o domínio absorvente $\Omega$ é um cojunto compacto de forma que
o resultado segue pelo Teorema de Extensão (\ref{teo:extensao}).
Vamos abordar apenas heuristicamente a prova do Teorema \ref{teo:lyapunov}. Uma prova formal
poderá ser encontrada em \cite{brauer}.
Vamos olhar para esse teorema geometricamente
(veja a Figura \ref{fig:lyapunov} para o caso $m=2$). De forma específica, vamos discutir
a condição $V'(x)<0$, onde $V$ é positiva-definida com respeito a uma região $\Omega$ de $\mathbb{R}^m$.

Seja $c$ uma constante positiva e considere a equação $V(x)=c$. Essa equação define superfícies
de nível em $\mathbb{R}^m$, as quais indicamos da forma
$$
S_c = \{ x\in\mathbb{R}^m : V(x)\leq c\}.
$$
Como $V$ é radialmente ilimitada, tais conjuntos são compactos e portanto positivamente invariantes.
Considere que a trajetória $x(t)$, solução de \eqref{x=ftx}, está fora de $\Omega$ e que
$S_c \supset \Omega$.
Por definição, $V'(x) = \left \langle \nabla V(x),f(x) \right \rangle$ onde o vetor
$\nabla V(x)$ é um vetor normal a superfície de nível $V(x) = c$. Então, a hipótese
de que $V'(x)<0$ significa que $f(x)$ deve apontar para o interior da
superfície de nível $S_c$. Por outro lado, o vetor $f(x)$ é um vetor tangente a
trajetória do sistema $\dot{x} = f(x)$ em cada ponto $x(t)$. Portanto, a trajetória cruza do exterior
para o interior da região limitada pela superfície de nível $V(x)=c$ para todo $c$. Isso faz 
com que quando $c$ tenda ao bordo de $\Omega$ as trajetórias entrem em $\Omega$. Uma vez em $\Omega$,
as trajetórias não deixam $\Omega$. Suponha que $x(t)\in\partial\Omega$ para algum $t=s>0$ e que $x(p)
\notin\Omega$ para $p$ ligeiramente maior que $s$. 
Esta hipótese contradiz o fato de que $V(x(t))$ é descrescente ao longo da trajetória
\footnote{Note que $\int_s^t V'(x(u))du = V(x(t))-V(x(s))<0$ para
todo $t>s$, portanto $V(x(t))$ é descrescente ao longo da trajetória $x(t)$.}.


\subsection{Sistema de Lorenz}\label{section:lorenz}

Um exemplo de dinâmica não-linear é a conhecida dinâmica de Lorenz. 
Utilizaremos tal dinâmica para efetuar as simulações computacionais dos principais resultados
que serão apresentados ao longo da presente dissertação.

Edward Lorenz (1917 - 2008) foi um meteorologista do MIT interessado em previsões do clima à tempos
longos. Com o advento do computador, muitas pessoas direcionavam suas tentativas de previsão do clima
a partir da análise numérica de equações que governam a dinâmica da atmosfera. Algo que era muito usado
eram as aproximações estatísticas, especialmente a regressão linear. Lorenz, porém, acreditava que
tais métodos eram essencialmente falhos, pois as reais equações de evolução estavam longe de serem lineares.
Depois de experimentar vários exemplos para testar suas ideias, Lorenz ao estudar um trabalho por
B. Saltzman \cite{bsaltz}, concernente a convecção de fluídos térmicos, percebeu que o sistema de equações estudado
por Saltzman era o ideal para provar o seu ponto de vista \cite{viana}. 
A maioria dos modelos de previsão
para o clima envolve equações diferenciais parciais, porém Lorenz encontrou um modelo muito mais
simples de ser analisado \cite{smale}.

O sistema de Lorenz, $\dot{\mathbf{x}}=f(\mathbf{x})$, é tal que
\begin{equation}\label{eq:lorenzsystem}
{f}(\mathbf{x}) =
\left(
\begin{matrix}
\sigma (y-x) \\
x(r-z) -y \\
-bz+xy
\end{matrix}
\right) ,
\end{equation}
com $\mathbf{x} = (x,y,z)\in\mathbb{R}^3$; $\sigma$, $r$ e $b$ constantes positivas, 
conhecidas respectivamente
por número de Prandtl, número de Rayleigh, e uma proporção física, e seus valores clássicos são:
$\sigma = 10$, $r=28$, e $b=8/3$. Esses são os valores que iremos considerar sempre que trabalharmos
com o sistema de Lorenz.

Lorenz, ao utilizar condições iniciais distintas, porém muito próximas, para este sistema, verificou que
as trajetórias divergiam gradativamente, de forma a produzir resultados bem diferentes, ou seja,
ele percebeu que tal sistema apresenta uma dependência sensível das condições iniciais,
e essa é a essência do comportamento caótico. Dessa forma,
assumindo que o clima tem um comportamento parecido com este modelo, Lorenz concluiu ser impossível 
fazer previsões à longo prazo.
Resumidamente, o sistema de Lorenz exibe uma dinâmica caótica \cite{viana}.

\subsubsection{Propriedades do Sistema de Lorenz}

\begin{figure}[!ht]
  \centering
  \includegraphics[width=0.5\textwidth]{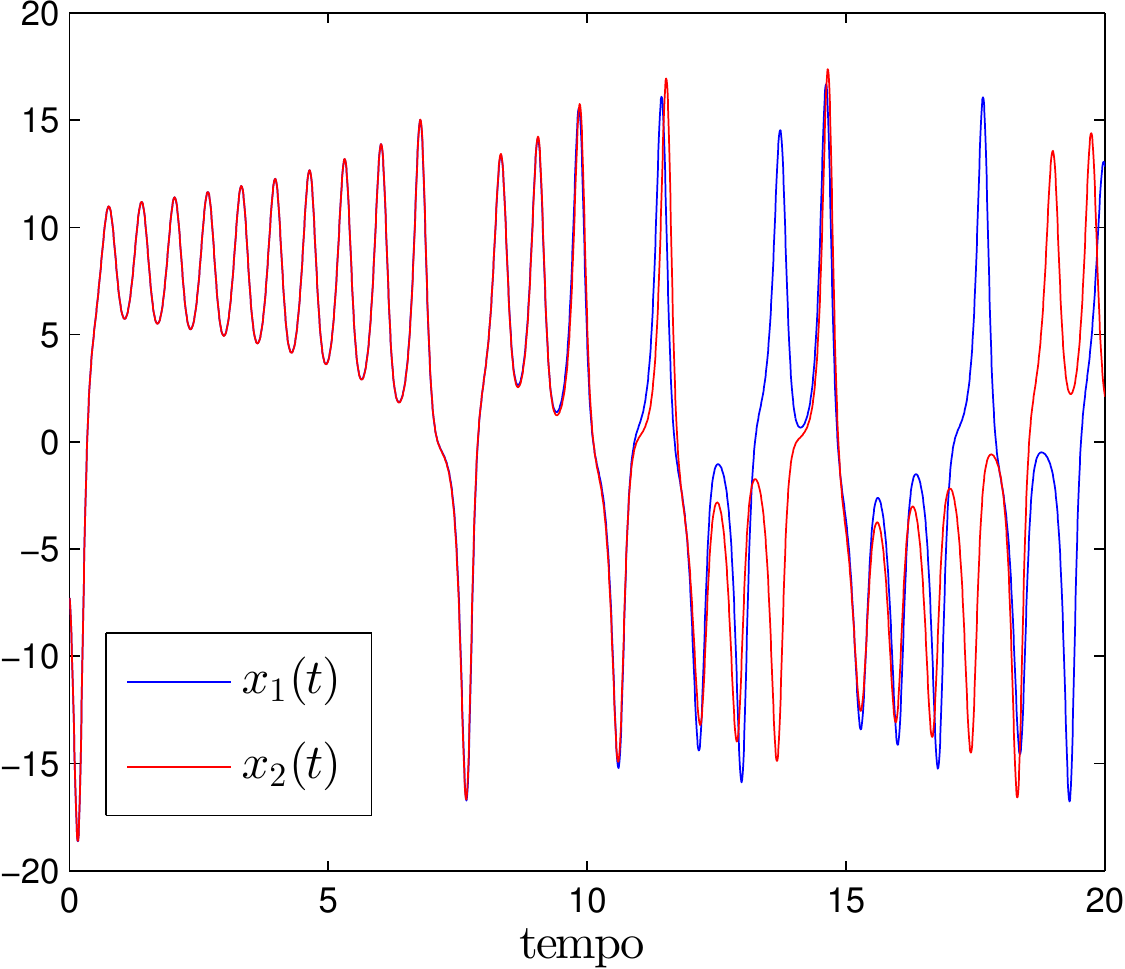} 
  \caption{O eixo horizontal representa o tempo e o vertical
as componentes ${x}_{1}(t)$ e ${x}_{2}(t)$ das trajetórias cuja diferença das condições iniciais,
considerando a norma euclidiana, é aproximadamente $0,014$.}
  \label{fig:x1x2} 
\end{figure}

A não-linearidade do sistema de Lorenz está associada aos termos $xy$ e $xz$.
Como o sistema de Lorenz é um sistema caótico, as trajetórias, para pequenas diferenças nas condições
iniciais, certamente irão divergir ao longo do tempo. 
Considere duas equações diferenciais, ambas com o sistema de Lorenz, a saber,
$\dot{\mathbf{x}}_1 = f(\mathbf{x}_1)$ e $\dot{\mathbf{x}}_2 = f(\mathbf{x}_2)$, onde 
$\mathbf{x}_1(t) = (x_{1}(t),y_{1}(t),z_{1}(t))$
e $\mathbf{x}_2(t) = (x_{2}(t),y_{2}(t),z_{2}(t))$ com condições inciais próximas, porém distintas.
A partir de integrações numéricas, utilizando
o método Runge-Kutta de quarta ordem, podemos verificar o fenômeno que aqui expomos.
Utilizando as seguintes condições iniciais para o sistema: $(-7,10,5)$ e $(-7.01,10.01,5)$,
exibimos as séries temporais ${x}_{1}(t)$ e ${x}_{2}(t)$ que representam as primeiras 
componentes das
trajetórias para o sistema considerado com as respectivas condições iniciais indicadas. O resultado é
apresentado na Figura \ref{fig:x1x2}.

\begin{figure}[!ht]
\centering
\includegraphics[width=0.6\textwidth]{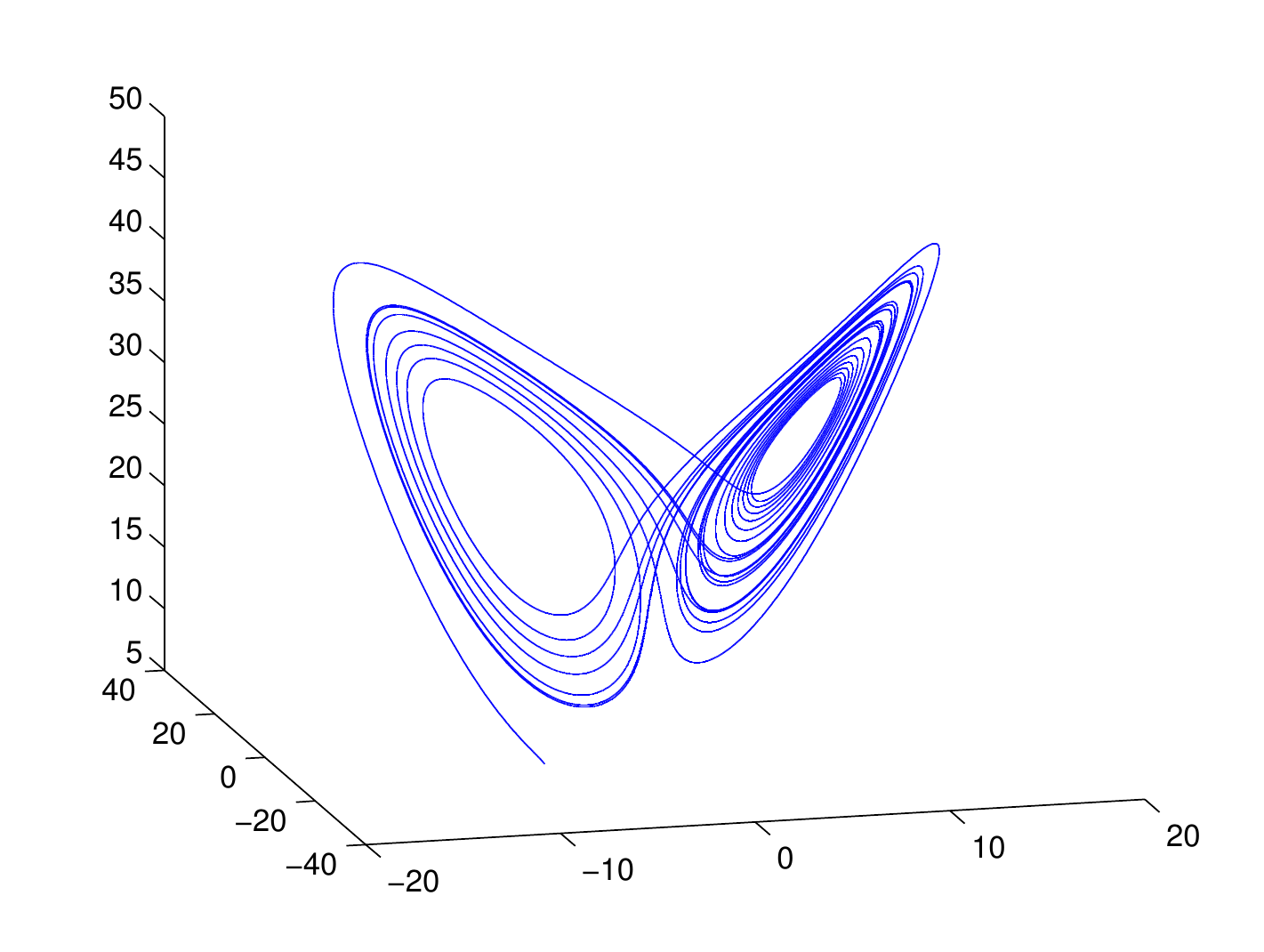}
\caption{Típica trajetória para o sistema de Lorenz considerando os parâmetros
$\sigma = 10$, $r=28$ e $b=8/3$.}
\label{fig:borboleta3d}
\end{figure}

Temos o mesmo sistema com diferença nas condições iniciais bem próximas, 
entretanto, essa diferença cresce com o tempo porém a mesma é limitada.
Mesmo considerando condições iniciais bem distintas, por exemplo, $(0,2,0)$ e 
$(0,-2,0)$, ambas as trajetórias divergem e se assemelham à borboletas confinadas em uma mesma região
do $\mathbb{R}^3$ (veja Figura \ref{fig:borboleta3d}). Isto nos elucida um importante fato sobre o sistema
de Lorenz: todas as soluções, que não são pontos de equilíbrio, tendem ao mesmo conjunto. Tal
conjunto é chamado de {atrator estranho de Lorenz} \cite{viana}. 

\subsubsection{Função de Lyapunov Para o Sistema de Lorenz}\label{subsection:flsl}

Vamos discorrer, nesta subseção, sobre a construção da função de Lyapunov para a dinâmica de Lorenz, 
e por consequência conseguimos também o domínio absorvente para a mesma.

Considere a função
$$
V(\mathbf{x}) = 1/2\langle \mathbf{x}-a,Q(\mathbf{x}-a) \rangle
$$
onde $\mathbf{x} = (x,y,z)$, $a=(0,0,2r)$ e 
$Q=\mbox{diag}(r,\sigma,\sigma)$, a qual deve ser positiva-definida
com relação a algum $\Omega\subset\mathbb{R}^3$. De forma explícita a função $V$ se lê
\begin{align*}
V(x,y,z) & = 1/2\langle (x,y,z-2r),Q(x,y,z-2r) \rangle\\
& = 1/2\langle (x,y,z-2r),(rx,\sigma y,\sigma (z-2r)) \rangle\\
& = 1/2[rx^2+\sigma y^2 +\sigma (z-2r)^2].
\end{align*}
Note que a condição $V(x,y,z) = c>0$ define um elipsoide em $\mathbb{R}^3$ com centro em $(0,0,2r)$.
Vamos então mostrar que as soluções de \eqref{eq:lorenzsystem} 
entram, e permanecem confinadas, dentro da região limitada pelo elipsoide 
$rx^2+\sigma y^2 +\sigma (z-2r)^2 = 2c$ para algum $c>0$.
Defina 
$$\Omega = \{(x,y,z)\in\mathbb{R}^3 : rx^2+\sigma y^2+\sigma(z-2r)^2\leq 2c\}.$$
A derivada da função $V$ se lê 
\begin{align*}
V'(x,y,z) & =  rx\frac{dx}{dt}+\sigma y\frac{dy}{dt} +\sigma (z-2r)\frac{dz}{dt}\\
& = r\sigma xy - r\sigma x^2+r\sigma xy-\sigma xyz-\sigma y^2-b\sigma z^2+2rb\sigma z
+\sigma xyz - 2r\sigma xy\\
& = -r\sigma x^2 -\sigma y^2 - b\sigma(z^2-2rz) = -\sigma \left[ rx^2+y^2+b(z^2-2rz)\right],
\end{align*}
mas $b(z^2-2rz)$ também pode ser escrito como
$b\left[ (z-r)^2-r^2\right]$ de forma que temos
\begin{equation}
V'(x,y,z) = -\sigma \left[ rx^2+y^2+b(z-r)^2-br^2\right].
\end{equation}
A equação $rx^2+y^2+b(z-r)^2 = \rho$ também define um elipsoide quando $\rho>0$.
Neste caso, teremos $V'<0$ quando $\rho> br^2$.
Além disso, para que sejam satisfeitas simultaneamente as duas
condições para que $V$ seja uma função de Lyapunov devemos escolher $c>0$ tal que o elipsoide $rx^2+y^2+b(z-r)^2=br^2$
esteja inteiramente contido no elipsoide
$rx^2+\sigma y^2+\sigma(z-2r)^2=2c$. Portanto, temos dois elipsoides
$$
E_1 : \frac{x^2}{br}+\frac{y^2}{br^2}+\frac{(z-r)^2}{r^2}=1 \quad \mbox{e}\quad
E_2 : \frac{x^2}{(2c/r)}+\frac{y^2}{(2c/\sigma)}+\frac{(z-2r)^2}{(2c/\sigma)}=1.
$$
Então, devemos procurar o valor mínimo
de $c$ para que se cumpra a condição $E_1\subset E_2$. Isso pode ser feito
através do método dos multiplicadores de Lagrange sobre dada
a restrição $V'<0$. Considerando os valores clássicos dos parâmetros do Lorenz, ou seja
$\sigma=10$, $r=28$, $b=8/3$, o valor de $c$ que devemos tomar é
$$
c=\frac{b^2r^2}{2(b-1)} .
$$
Para maiores detalhes veja \cite{sparrow1982lorenz}.

Assim, podemos redefinir o conjunto $\Omega$ da forma
\begin{equation}\label{eq:lorenzatractor}
\Omega = \{(x,y,z)\in\mathbb{R}^3 : rx^2+\sigma y^2+\sigma(z-2r)^2\leq {b^2r^2}/({b-1})\},
\end{equation}
o qual é uma região limitada e fechada de $\mathbb{R}^3$, portanto um conjunto compacto.
Assim, como são satisfeitas as condições do Teroema de Lyapunov (\ref{teo:lyapunov}),
qualquer solução do sistema de Lorenz que tem condições iniciais fora de $\Omega$
convergirá à $\Omega$ e permanecerá dentro de $\Omega$ por todo o tempo futuro.

Podemos então, enunciar a seguinte

\begin{proposicao}
O sistema de Lorenz é um sistema dissipativo e suas 
trajetórias certamente entram, a tempo finito, no domínio
absorvente \eqref{eq:lorenzatractor}.
\end{proposicao}

\chapter{Estabilidade em EDO's Lineares Não-Autônomas}
\label{cap:edolinear}

Vamos estudar algumas noções acerca de equações diferenciais lineares não-autônomas, em especial,
a estabilidade de seu ponto fixo trivial. Nosso objetivo é
transformar o problema de estudar a sincronização e sua estabilidade num problema de estudar
a estabilidade da solução trivial de uma equação diferencial linear não-autônoma. 

Considere $A(t):J\rightarrow\Mat(\mathbb{R},n)$ uma função matricial contínua e 
uniformemente limitada
que toma valores $t\in J\subseteq\mathbb{R}_+$ e a equação diferencial linear não-autônoma
\begin{equation}\label{edo}
\dot{x} = A(t)x.
\end{equation}
A equação \eqref{edo} 
satisfaz as condições de existência e unicidade dada a condição inicial 
$x(s)=x_s$ e mais do que isso, as soluções estão definidas para todo $t\in\mathbb{R}_+$ \cite{teschl}.

Note que $x(t) = \mathbf{0}$ para todo $t \in \R_+$ é uma solução de \eqref{edo} para 
um dado problema de valor inicial. Chamamos tal solução de \textbf{solução trivial}.
Por todo este capítulo focaremos o estudo da estabilidade da solução trivial de \eqref{edo}
com condição inicial $x(s)\in\mathbb{R}^m$.

Para equações diferenciais lineares autônomas, a estabilidade de $x=\mathbf{0}$ pode 
ser caracterizada em termos dos autovalores da matriz constante $A(t)=A$. 
Quando, porém, consideramos que as equações diferenciais lineares são não-autônomas,
os autovalores já não servem mais como critério de caracterização da estabilidade \cite{coppel}.
De fato, podemos ter casos em que a 
matriz $A(t)$ tem todos os autovalores com parte real negativa e ainda assim, a solução 
trivial seja instável, ver \cite{coppel} para exemplos. Ou ainda, casos em que as soluções 
não-triviais de \eqref{edo} convergem ao vetor nulo, porém uma pequena pertubação pode destruir 
sua estabilidade, ver \cite{barreira2} e \cite{josic} para exemplos.

\section{Alguns Tipos de Estabilidade}\label{section:tiposdeestabilidade}
Para tratar os principais resultados aqui expostos, necessitamos de
algumas definições.
A teoria de estabilidade é ampla e vários conceitos sobre estabilidade são encontrados
\cite{teschl}, \cite{claus}. Abordamos alguns desses conceitos que servirão ao nosso propósito.
\begin{definicao}\label{def:estabilidadelyapunov}
Dizemos que o ponto  fixo $x=\mathbf{0}$ de \eqref{edo} é \textbf{estável no sentido de Lyapunov} 
em $t=s$ se para todo $\epsilon>0$ existe $\delta = \delta(s,\epsilon)$ tais que
$$
\Vert x(s)\Vert <\delta\quad \Rightarrow\quad \Vert x(t)\Vert <\epsilon \quad \forall \mbox{ }t\geq s.
$$
\end{definicao}

Assim, a estabilidade segundo Lyapunov é
definida sobre um tempo inicial e, em particular, esse tipo de estabilidade não requer que as
trajetórias que começam próximas da solução trivial convirjam assintoticamente à mesma.
Outro importante conceito de estabilidade é a \textbf{estabilidade uniforme}. É um caso
especial de estabilidade segundo Lyapunov, que garante que o ponto
fixo não esteja perdendo a estabilidade. Assim, a diferença de tal conceito
com relação ao definido em \ref{def:estabilidadelyapunov} é que $\delta$ não depende do tempo inicial
$s$.

\begin{definicao}
Dizemos que o ponto fixo $x=\mathbf{0}$ de \eqref{edo} é
\textbf{assintoticamente estável} em $t=s$ se o mesmo é estável no sentido de Lyapunov e se é 
localmente atrator, isto é, existe uma vizinhança 
$U(\mathbf{0})\ni x(s)$ tal que $x(t)$ converge assintoticamente à $\mathbf{0}$ quando 
$t\rightarrow\infty$, com $x(t) \in U(\mathbf{0})$ para cada $t$. Ou seja
$$
\lim_{t\rightarrow\infty} x(t) = \mathbf{0}.
$$
\end{definicao}

Estamos, porém, interessados na estabilidade uniformemente assintótica.
O interesse
nesse tipo de estabilidade está baseada no fato de que a mesma é preservada quando o sistema
sofre pequenas pertubações caracterizando assim a propriedade de persistência.

\begin{definicao}\label{def:estbua}
O ponto fixo $x=\mathbf{0}$ de \eqref{edo} é 
\textbf{uniformemente assintoticamente estável} se existe $\delta>0$, independente de $t$, tal que,
para todo $\epsilon >0$, existe $\tau = \tau(\epsilon)$
satisfazendo $\Vert x(s) \Vert < \delta$ e $\Vert x(t) \Vert < \epsilon$ para todo $t\geq s+\tau$.
\end{definicao}


\section{Operador de Evolução}
Para o nosso propósito, precisamos de uma forma fechada para expressar as soluções de \eqref{edo}. 
A partir da teoria de 
equações
diferenciais (\cite{teschl}) nós podemos escrever tal solução sob a forma
\begin{equation}
x(t) = T(t,s)x(s)
\end{equation}
onde $T(t,s)$ é chamado de \textbf{operador de evolução} e $x(s)$
é a condição inicial. O operador de evolução satisfaz as seguintes propriedades:
\begin{enumerate}
 \item $T(t,s)T(s,u)=T(t,u)$
 \item $T(t,s)T(s,t) = I_m$
\end{enumerate}

A propriedade 1 é imediata tendo-se em vista que podemos escrever 
$T(t,s)x(s) = \phi_t(x(s)) = \phi_{t,s} (x)$ e as propriedades sobre o fluxo (ver página 
\pageref{pag:fluxo}). Considerando a propriedade 1, podemos sempre escolher $t=u$, produzindo 
$T(u,s)T(s,u) = T(u,u)=I_m$, donde segue que $T(t,s)$ é um isomorfismo com inversa dada por 
$T^{-1}(t,s)= T(s,t)$, valendo-se assim a propriedade 2.

\begin{exemplo}
Quando o sistema \eqref{edo} é autônomo, ou seja, quando $A(t) = A$, o operador
de evolução do mesmo se lê
$$
T(t,s) = \exp{[(t-s)A]}.
$$
\end{exemplo}


Existem processos para se encontrar explicitamente o operador de evolução (ver \cite{teschl} por exemplo),
porém não estamos interessados em obter analiticamente tal operador mas em estimar sua magnitude.
A seguir, daremos uma definição necessária à abordagem dos principais resultados aqui expostos.

\begin{definicao}\label{def:contracaouniforme}
 Seja $T(t,s)$ o operador de evolução associado ao sistema \eqref{edo}. 
Dizemos que tal operador possui \textbf{contração uniforme} se
$$
\Vert T(t,s) \Vert \leq k e^{-\eta(t-s)}
$$
com $k$ e $\eta$ positivos.
\end{definicao}

O termo contração uniforme, associado ao operador de evolução, significa que tal operador de
fato possui contração, como expressado pela norma, pois a mesma converge exponencialmente
rápido à zero, e que as constantes
$\eta$ e $k$ não possuem dependência do tempo inicial $s$.

O próximo teorema à ser abordado, relaciona intimamente as definições \ref{def:estbua} e
\ref{def:contracaouniforme}. 

\begin{teorema}\label{teo:contracaouniformexestabilidade}
 Seja $T(t,s)$ o operador de evolução da equação \eqref{edo} e considere a solução trivial da 
mesma $x(t) = \mathbf{0}$. Tal 
solução
é uniformemente assintoticamente estável se e somente se $T(t,s)$ possui contração uniforme.
\end{teorema}

\noindent
\textbf{Prova:} Faremos primeiro a volta. Escrevamos a solução de \eqref{edo} sob a forma 
$x(t) = T(t,s)x(s)$ e então 
tomamos a norma
em ambos os lados, assim
\begin{align*}
 \Vert x(t)\Vert & = \Vert T(t,s)x(s)\Vert \\
& \leq \Vert T(t,s)\Vert \Vert x(s) \Vert \\
& \leq ke^{-\eta(t-s)}\Vert x(s) \Vert.
\end{align*}
Queremos mostrar que $x(t) = \mathbf{0}$ é uniformemente assintoticamente estável, 
ou seja, deseja-se que exita $\delta>0$ tal que para todo 
$\epsilon >0$, exista $\tau=\tau(\epsilon)>0$ tais que, sendo
$\Vert x(s) \Vert < \delta \Rightarrow \Vert x(t) \Vert <\epsilon \quad \forall 
t\geq s+\tau$. Dessa forma, devemos ter
\begin{align*}
 \delta ke^{-\eta(t-s)} & <\epsilon \quad \Rightarrow \\
-\eta(t-s) & < \ln{(\epsilon/{\delta k})} \quad \Rightarrow \\ 
t & > s + (1/\eta) \ln{({\delta k}/{\epsilon})}.
\end{align*}
Assim, basta tomar $\tau = (1/\eta) \ln{({\delta k}/{\epsilon})}$ para garantir o que desejamos.

Para provar a ida, considere que a solução $x(t) = \mathbf{0}$ é uniformemente 
assintoticamente estável. Então
existe $\delta>0$ tal que para todo $\epsilon>0$, exite $\tau=\tau(\epsilon)>0$, com $\Vert x(s) 
\Vert \leq \delta$, e 
$\Vert x(t) \Vert \leq \epsilon$,
onde $t\geq s+\tau$.
Escolha então $\epsilon = \delta/c$, $c>1$, e considere a sequência $t_n = s+n\tau$, 
com $n \in \mathbb{N}$.
Perceba então que
$$\Vert x(t) \Vert = \Vert T(t,s)x(s)\Vert \leq\frac{\delta}{c}
$$
para qualquer ${\Vert x(s) \Vert}/\delta \leq 1$, o que implica
$$
\Vert T(t,s) \Vert = \sup_{\Vert u\Vert \leq 1}{\Vert 
T(t,s)u\Vert} \leq \frac{1}{c},$$
onde $u = \Vert x(s)\Vert / {\delta}$.
Pela propriedade 1 do operador de evolução, ou seja, $T(t,u)=T(t,s)T(s,u)$, 
podemos fazer
\begin{align*}
\Vert T(t_2,s)\Vert & = \Vert T(s+2\tau,s+\tau)T(s+\tau,s)\Vert \\
& \leq \Vert T(s+2\tau,s+\tau)\Vert \Vert T(s+\tau,s)\Vert \\
& \leq \frac{1}{c^2}.
\end{align*}
Assim, por indução, temos
$$
\Vert T(t_n,s)\Vert \leq \frac{1}{c^n},
$$
onde, podemos impor que $1/c^n = e^{-\eta(t_n-s)}$. Então $-\eta(t_n-s) = \ln{c^{-n}}$, mas  
$t_n-s=n\tau$,
o que implica que podemos tomar $\eta =\ln{(c)}/\tau$, de forma à conseguir que
$$
\Vert {T}(t_n,s)\Vert \leq e^{-\eta(t_n-s)}.
$$

Fica ainda um último ponto a resolver, que é exatamente que neste caso, temos uma 
contração uniforme para o operador de evolução, porém
escrita em tempo discreto. Para torná-lo em tempo contínuo, utilizamos
o Teorema do Valor Intermediário e consideramos o caso geral em que $t=s+v+n\tau$, 
onde $0 \leq v <\tau$. Assim
\begin{align*}
\Vert T(t,s)\Vert & \leq e^{-\eta(t_n-s)} \\
& = e^{-\eta(t-s-v)} \\
& = ke^{-\eta(t-s)}
\end{align*}
com $k=e^{\eta v}$.
$\blacksquare$

\section{Persistências da Solução Trivial}

Dizemos que o sistema $\dot{x} = f(t,x)$, $f:\R_+\times\R^m\rightarrow\R^m$, 
é persistente, quando o mesmo mantém suas propriedades
qualitativas sob pequenas pertubações ou mudanças no campo de vetores $f$,
veja no livro \cite{guckenheimer1983nonlinear} uma discussão mais detalhada.
Mas, o que entendemos por ``pequenas perturbações''? Muitos sistemas
físicos interessantes possuem, tipicamente, uma dependência de parâmetros, os quais aparecem
nas definições das equações como é o caso do sistema de Lorenz \eqref{eq:lorenzsystem}.
Quando esses parâmetros são variados podem acontecer mudanças na estrutura qualitativa das
soluções de forma que os pontos fixos podem perder sua estabilidade a medida que essa variação
ocorre - essas mudanças são chamadas de bifurcações. Então essa leve mudança nos parâmetros
é um dos exemplos do que pode ser uma pequena perturbação no campo de vetores.

No caso particular de nosso interesse temos que $f$ é um campo linear com $f(t,\mathbf{0})=
\mathbf{0}$. Seja $g(t,x)$ uma perturbação de $f$ tal que $g(t,\mathbf{0})=\mathbf{0}$ e
considere o sistema perturbado $\dot{y} = \tilde{f}(t,y)$ onde $\tilde{f} = f+g$.
Note que $y=\mathbf{0}$ é sempre solução do sistema perturbado. A pergunta então
é: sendo $x=\mathbf{0}$ estável, $y=\mathbf{0}$ é também estável? Sob quais condições
a estabilidade é preservada? Esta é uma pergunta interessante pois pode acotecer
que tenhamos uma perturbação muito pequena e mesmo assim, a estabilidade
da solução trivial seja destruída, veja \cite{barreira2} para exemplos.
Diremos portanto que a estabilidade de $x=\mathbf{0}$ em \eqref{edo} é persistente quando
a mesma é preservada sob pertubações pequenas no campo de vetores.

Voltando agora à teoria das contrações uniformes, sabemos que a mesma está muito
relacionada com a propriedade de persistência da estabilidade. Os resultados
desta seção tratam dessa relação.

%
%
%
\begin{teorema}[Persistência]\label{persistencia}
Suponha que $A(t)$ é uma função matricial contínua em um intervalo $J\subset \mathbb{R}_+$,
 referente à equação \eqref{edo}. Suponha também
que o operador de evolução $T(t,s)$ da mesma possui contração uniforme. 
Considere $B(t)$ uma função matricial contínua em $J$  satisfazendo
$$
\sup_{t\in J}\Vert B(t)\Vert = \delta_0 < \frac{\eta}{k}.
$$
Então o operador de evolução $\widehat{T} (t,s)$ da equação perturbada
\begin{equation}\label{pert}
 \dot{y} = [A(t)+B(t)]y
\end{equation}
também satisfaz $\Vert\widehat{T} (t,s)\Vert \leq ke^{-\gamma(t-s)}$ com $\gamma = \eta -\delta_0 k$,
onde $\eta$ e $k$ são tais que $\Vert T(t,s)\Vert\leq ke^{-\eta(t-s)}$.
\end{teorema}

Para a prova do teorema sobre a persistência da estabilidade da solução trivial, utilizaremos
dois importantes resultados da teoria de equação diferenciais, a saber,
o Lema de Grönwall e o Método da Variação dos Parâmetros,
que podem ser encontrados no Apêndice \ref{ape:edoB} com mais detalhes.

\noindent
\textbf{Prova:} Usando o método da variação dos parâmetros em \eqref{pert}, obtemos
$$
y(t) = T(t,s)y(s) + 
\int_s^t{T(t,u)B(u)y(u)du}.
$$
Tomando a norma sobre a equação acima, e usando a desigualdade triangular, ficamos com
\begin{align*}
\Vert y(t)\Vert & \leq\Vert T(t,s)y(s)\Vert + \left\Vert\int_s^t{T(t,u)B(u)y(u)du}\right\Vert \\
& \leq ke^{-\eta(t-s)}\Vert y(s)\Vert +\int_s^t{\Vert T(t,u)\Vert \Vert B(u)\Vert \Vert y(u)\Vert du} \\
& \leq ke^{-\eta(t-s)}\Vert y(s)\Vert +\int_s^t{\delta_0 ke^{-\eta(t-u)}\Vert y(u)\Vert du},
\end{align*}
onde as sucessivas desigualdades foram obtidas através da propriedade sub-multiplicativa da
norma induzida (Lema \ref{lema:matrixpropeties}).

Neste ponto, tomamos a função $\omega(t) = e^{\eta t}\Vert y(t)\Vert$. Multiplicando a última desigualdade
por $e^{\eta t}$, obtemos
$$
\omega (t) \leq k\omega (s) + \int_s^t{\delta_0 k \omega (u)du}
$$
a qual podemos usar o lema de Grönwall, de forma que obtemos
$$\omega (t) \leq k\omega(s)e^{\delta_0 k(t-s)}.
$$
Substituindo nesta última desigualdade a função $\omega(t)$, segue que
$$
\Vert y(t)\Vert \leq ke^{(-\eta +\delta_0 k)(t-s)}\Vert y(s)\Vert.
$$
Logo, podemos tomar $\Vert \widehat{T}(t,s)\Vert \leq ke^{-\gamma(t-s)}$, com 
$\gamma = \eta -\delta_0 k$.
$\blacksquare$
\\

O próximo resultado é importante pois o mesmo dá as hipóteses sobre a qual o resto de Taylor,
advindo da linearização pela expansão em série de Taylor, não destrói a estabilidade
da solução trivial.

\begin{proposicao}[Princípio da Linearização]\label{pro:principio2}
Considere que o operador de evolução da equação \eqref{edo} tem contração uniforme. Considere
a equação perturbada
\begin{equation}\label{eq:atmaisr}
\dot{y} = A(t)y + R(t,y)
\end{equation}
onde $R(t,y)$ é tal que $R(t,\mathbf{0})=\mathbf{0}$, $\forall \mbox{ } t\geq0$, e a seguinte 
propriedade é satisfeita
$$
\forall\mbox{ } \epsilon>0, \exists \mbox{ } \delta>0 : \mbox{se } \Vert y\Vert \leq \delta 
\Rightarrow  \Vert R(t,y)\Vert \leq \epsilon\Vert y\Vert, \forall \mbox{ } t\geq0.
$$
Então o operador de evolução
de \eqref{eq:atmaisr} tem contração uniforme.
\end{proposicao}

\noindent
\textbf{Prova:}
Utilizando a variação dos parâmetros na equação \eqref{eq:atmaisr} tem-se
\begin{equation}\label{eq:vplema}
y(t) = T(t,s)y(s) + \int_s^t T(t,u)R(u,y(u))du
\end{equation}
onde $T(t,s)$ é o operador de evolução de \eqref{edo}, que por hipótese
satisfaz $\Vert T(t,s)\Vert \leq k e^{-\eta(t-s)}$. Calculando a norma em \eqref{eq:vplema}
e utilizando as hipóteses dadas tem-se
\begin{align}
\Vert y(t)\Vert & \leq \Vert T(t,s)\Vert \Vert y(s)\Vert +\int_s^t \Vert T(t,u)\Vert \Vert R(u,y(u)\Vert du\\
& \leq ke^{-\eta(t-s)}\Vert y(s)\Vert +\int_s^t \epsilon k e^{-\eta(t-u)} \Vert y(u)\Vert du.\label{eq:vplema2}
\end{align}


Considere a função real $\omega(t) = e^{\eta t}\Vert y(t)\Vert$. Multiplicando $e^{\eta t}$
em \eqref{eq:vplema2} teremos
\begin{equation}\label{eq:vplema3}
\omega(t) \leq k\omega(s)+k\epsilon\int_s^t \omega(u) du.
\end{equation}
A partir daqui, utilizamos o Lema de Grönwall em \eqref{eq:vplema3} e seguimos os mesmos
passos utlizados na prova do Teorema \ref{persistencia} mostrando que
$$
\Vert y(t)\Vert \leq ke^{(-\eta +\epsilon k)(t-s)}\Vert y(s)\Vert.
$$
Escolhendo $\epsilon<\eta /k$, encerra-se a demonstração. $\blacksquare$

Algo interessante que a Proposição \ref{pro:principio2} nos diz é que a estabilidade da solução
trivial de \eqref{edo} não muda mesmo considerando a pertubação não-linear $R(t,y)$
com a propriedade descrita, isso porque podemos tomar $\epsilon$ tão pequeno quando desejarmos.
Note ainda que a condição dada sobre $R(t,y)$ é a mesma condição dada sobre o resto de Taylor
\cite{lima1987curso}, isto é, $\lim_{y\rightarrow \mathbf{0}} R(t,y)/{\Vert y\Vert} = 0$.

%

\section{Critérios para Contração Uniforme}

O resultado a seguir fornece um importante critério a respeito de como garantir que o operador 
$T(t,s)$ de \eqref{edo} tenha contração uniforme. O Teorema da Diagonal Dominante garante que os 
coeficientes da diagonal da matriz $A(t)$ de \eqref{edo} controlam a estabilidade do sistema.

\begin{teorema}[\textbf{da Diagonal Dominante}]\label{teo:criterio1}
Seja $A(t) = [A_{ij}(t)]_{i,j=1}^{m}$ uma função matricial contínua e limitada com
$t\in \mathbb{R}_+$. Suponha que existe uma constante $\eta >0$ tal que
\begin{equation}\label{eq:dd}
 A_{ii}(t)+ \sum_{j=1, j\neq i}^{m}{\vert A_{ij}(t)\vert} \leq -\eta <0
\end{equation}
para todo $t\in \mathbb{R}_+$ e $i=1,\dots,m$. Então o operador de evolução de \eqref{edo}
tem contração uniforme.
\end{teorema}

\noindent
\textbf{Prova:}
Seja $x=(x_1,\cdots,x_m)\in\mathbb{R}^m$. Usaremos a norma $\Vert x\Vert_\infty = 
\max_i{\vert x_i\vert}$ para estimar $x(t)$. Fixe $u>0$, então para
algum $1\leq i\leq m$ tem-se que $\Vert x(u)\Vert_\infty = \vert x_i(u)\vert$ ou ainda
$\Vert x(u)\Vert_\infty^2 = {x_i}^2(u)$. Analisaremos a norma $\Vert x(t)\Vert_\infty^2$
pois tal função é diferenciável. Note porém que poderia acontecer de que exatamente
neste $u$ fixado a norma $\Vert x(u)\Vert_\infty^2$ ser atingida em mais de um $i$. Suponha
por exemplo, sem perder a generalidade, que esse máximo é atingido em $h$ e $k$, $h,k\in\{1,\cdots,m\}$.
Então surge a pergunta: para $t$ numa vizinhança de $u$, qual seria maior,
${x_h}^2(t)$ ou ${x_k}^2(t)$? A resposta é: independente de quem seja maior,
esse máximo continuará sendo o máximo por pelo menos um intervalo $(u,u+\epsilon)=I$
com $\epsilon>0$ suficientemente pequeno. Suponha que o que acabamos de afirmar não seja verdade,
isto é, suponha que para $t\in I$, o máximo alterna-se entre as componentes $h$ e $k$ um 
número arbitrário de vezes. Neste caso,
teríamos que $\lim_{\epsilon\rightarrow 0}{\Vert x(u+\epsilon)\Vert_\infty^2}$ não existe e portanto
$d/dt \Vert x(t)\Vert_\infty^2$ não existiria $t=u$, o que é absurdo pois todas
as componentes ${x_i}^2(t)$ são diferenciáveis visto que $x(t)$, solução
não-trivial de \eqref{edo}, é diferenciável. Portanto podemos impor, sem perder a generalidade,
que se $ {x_k}^2(u) =\Vert x(u)\Vert_\infty^2$ então ${x_k}^2(t) = 
\Vert x(t)\Vert_\infty^2$ para todo $t\in I$.
Então, tomando $t\in I$ tem-se
$$
\frac{1}{2}\frac{d}{dt}\Vert x(t)\Vert_\infty^2 = x_k(t)\frac{d}{dt}x_k(t).
$$
Mas, note que de \eqref{edo} tem-se $$\frac{d}{dt}x_k(t) = \sum_{j=1}^m A_{kj}(t)x_j(t),$$
assim
\begin{align*}
\frac{1}{2}\frac{d}{dt}\Vert x(t)\Vert_\infty^2 & = \left[ A_{kk}(t){x_k}^2(t) + \sum_{j\neq k}{A_{kj}(t)x_j(t)x_k(t)}\right] \\
& \leq \left[ A_{kk}(t){x_k}^2(t) + \sum_{j\neq k}{\vert A_{kj}(t)\vert {x_k}^2(t)}\right]
=  \left[ A_{kk}(t) + \sum_{j\neq k}{\vert A_{kj}(t)\vert}\right] {x_k}^2(t),
\end{align*}
onde por hipótese vale \eqref{eq:dd}. Então
$$
\frac{d}{dt}\Vert {x}(t)\Vert_\infty^2 \leq -2\eta  {x_k}^2(t)  =-2\eta\Vert x(t)\Vert_\infty^2.
$$
Tomando a integral, definida em $[s,t]\subset I$, temos
$$
\Vert x(t)\Vert_{\infty}^2 \leq \Vert x(s)\Vert_\infty^2 
-2\eta\int_s^t{\Vert x(v)\Vert_\infty^2dv},
$$
da qual podemos usar a desigualdade de Grönwall, obtendo
\begin{equation}\label{eq:xinf}
\Vert x(t)\Vert_\infty \leq e^{-\eta(t-s)}\Vert x(s)\Vert_\infty
\end{equation}

Como a derivada de $\Vert x(t)\Vert_\infty^2$ é estritamente negativa,
então $\Vert x(t)\Vert_\infty^2$ não possui máximo local em $t\in I$ e portanto
vale que $\Vert x(t)\Vert_\infty^2 < \Vert x(s)\Vert_\infty^2$ para $t>s$. 
Mais do que isso, neste caso, é possível garantir que $\Vert x(t)\Vert_\infty$ é uma função 
monotonicamente decrescente em $\mathbb{R}_+$ (veja \cite{coppel}, a partir da página 55 para 
mais detalhes), então a desigualdade \eqref{eq:xinf} será válida para todo $t\in\mathbb{R}_+, t\geq s$.
Como em $\Mat(\mathbb{R},m)$ todas as normas são equivalentes, então existem $c>0$ e $k>0$ tais
que $c\Vert T(t,s)\Vert_\infty\leq \Vert T(t,s)\Vert \leq k\Vert T(t,s)\Vert_\infty$, portanto
$\Vert T(t,s)\Vert\leq k e^{-\eta(t-s)}$. $\blacksquare$
\\

Agora, enunciamos um novo critério que está relacionado com o espectro da parte 
simétrica da matriz $A(t)$ de \eqref{edo}.

\begin{teorema}[\textbf{da Parte Simétrica}]\label{teo:criterio2}
Considere o sistema \eqref{edo}. 
Suponha que $A(t)$ é uniformemente limitada em $t\in\mathbb{R}_+$ e que sua parte simétrica é negativa-definida.
Então o operador de evolução de \eqref{edo} satisfaz $\Vert T(t,s)\Vert_2 \leq  e^{\eta(t-s)}$,
onde $\eta<0$ é uma cota superior para os autovalores da parte simétrica de $A(t)$.
\end{teorema}

\noindent
\textbf{Prova} 
Considere $x\in\mathbb{R}^m$ e $\Vert x\Vert_2^2 = \langle x,x\rangle$. Considere
a solução não-trivial $x(t)$ de \eqref{edo} definida em $\mathbb{R}_+$. A função
$\Vert x(t)\Vert_2^2$ é contínua e diferenciável. Então
\begin{align*}
\frac{1}{2}\frac{d}{dt}\left \| x(t) \right \|_2^2 & =  
1/2\left\langle A(t)x,x\right\rangle +1/2\left\langle x,A(t)x\right\rangle\\
& = \left\langle  x,\frac{A(t)+A^\dagger(t)}{2} x\right\rangle =  \left\langle  x,B(t)x\right\rangle
\end{align*}
onde $B(t)$ é a parte simétrica de $A(t)$. Neste momento consideraremos a análise para
cada $t$ fixado.

Como $B$ é uma matriz simétrica então, pelo Teorema Espectral, seus autovetores ortonormalizados
formam uma base de $\mathbb{R}^m$. Podemos representar $x$ nessa base, $x=\sum_{i=1}^m c_iv_i$, de 
forma que faz-se valer 
\begin{align*}
\langle x,Bx\rangle & \leq \left\langle \sum_{i=1}^m c_iv_i , \sum_{i=1}^m c_i\lambda_i v_i\right\rangle\\
& \leq \lambda_m \langle x,x\rangle,
\end{align*}
onde $\lambda_1\leq \cdots \leq \lambda_m <0$ são os autovalores de $B$.
Como $B(t)$ é uniformemente limitada em $t$ e negativa-definida podemos garantir que existe $\eta<0$ tal que
$\lambda_m(t) < \eta$ para todo todo $t\in\mathbb{R}_+$.
Dessa forma, temos a relação
$$
\frac{d}{dt}\left \| x(t) \right \|_2^2 \leq 2\eta \Vert x(t)\Vert_2^2
$$
à qual podemos tomar a integral, definida em $[s,t]\subset\mathbb{R}_+$, $s<t$, e utilizar
a desigualdade de Grönwall, produzindo
$$
\Vert {x}(t)\Vert_2  \leq e^{\eta(t-s)}\Vert {x}(s)\Vert_2.
$$
$\blacksquare$

\chapter{Tópicos em Teoria dos Grafos}\label{cap:grafos}


Abordaremos os fundamentos de teoria dos grafos para estruturar e modelar
redes de elementos interagentes afim de fundamentar matematicamente
o objetivo de estudar a sincronização e suas estabilidades em tais redes. 
Por todo o presente capítulo, focaremos
o estudo na teoria espectral dos grafos, ou seja, o estudos dos autovalores da uma das matrizes 
que caracteriza a estrutura de conexão, a saber, o Laplaciano. 
%
Inicialmente daremos algumas definições
necessárias ao longo do presente capítulo. Consideraremos apenas grafos finitos.

\begin{definicao}
Um \textbf{grafo} $G=(V(G),E(G))$ consiste de um conjunto finito de vértices $V(G)$ e um conjunto,
consequentemente finito, de
arestas $E(G)$, onde cada aresta consiste de um par de vértices.
\end{definicao}
Quando ao referirmos que uma aresta $\{u,v\}\in E(G)$ vamos utilizar
a notação $u\sim v$ e dizer que $u$ é vizinho de $v$.

Um grafo \textit{direcionado} consiste de um conjunto de vértices $V(G)$
e um conjunto de arestas $E(G)$, onde cada aresta consiste de um par de vértices
ordenados.

\begin{definicao}
Dizemos que um grafo é \textbf{não-direcionado} se os pares de vértices $\{u,v\}\in E(G)$
são não-ordenados.
\end{definicao}
Assim,
em grafos não-direcionados, se $u$ é vizinho de $v$ então reciprocamente $v$ é
vizinho de $u$.

\begin{definicao}
Diremos que um grafo é \textbf{simples} quando o mesmo não possui auto-conexões ou
múltiplas arestas entre dois vértices. 
\end{definicao}

\begin{definicao}
O \textbf{grau} de um vértice $i$ é numericamente igual ao número de vizinhos que o mesmo possui.
Indicaremos o grau do vértice $i$ por $g_i$.
\end{definicao}

Um \textit{caminho} em um grafo é uma sequência de vértices distintos e conectados.

\begin{definicao}
Um grafo é \textbf{conexo} se para quaisquer dois vértices $u$ e $v$ existe pelo
menos um caminho de $u$ até $v$.
\end{definicao}
De forma semelhante, dizemos que $G$ é \textit{não-conexo} se em $G$ existe algum vértice que não
pode ser alcançado por qualquer outro.
Para todo grafo não-conexo $G$, sempre podemos decompô-lo
da forma simbólica $G=G_1\cup G_2\cup\cdots\cup C_r$, se $G$ tem $r$ componentes conexas,
de forma que os subgrafos  $G_1,G_2,\cdots,G_r$ são chamados componentes conexas de $G$.

Trabalharemos apenas com grafos que são simples, {não-direcionados} e conexos. Um grafo pode ser 
simples e ainda assim representar matematicamente uma rede complexa.
Para que o leitor não faça confusão desses conceitos, note
que a rede é simples como traz a definição acima, porém pode ter uma estrutura não regular de conexão
entre seus elementos. Um grafo simples e regular é aquele 
o qual podemos identificar precisamente a regra que gera o conjunto de vértices.

Para quaisquer dois pares de vértices $u,v\in V(G)$, a \textit{distância} entre $u$ e $v$ é o menor
caminho dentre todos os possíveis entre $u$ e $v$.

\begin{definicao}
Em um grafo $G$, o \textbf{diâmetro} de $G$ é a máxima distância sobre todos os pares
de vértices em $G$. Utilizaremos a letra ``$d$'' para indicar o diâmetro de $G$.
\end{definicao}



%
%

\section{Matrizes de Adjascência e Laplaciana}

Seja $G$ um grafo simples e não-direcionado
com $n$ vértices. A \textbf{matriz de adjacência} de $G$, denotada por ${A}_G$ ou
simplesmente ${A}=[A_{ij}]_{i,j=1}^n$, quando a indicação do tipo de grafo não for necessária, 
é uma matriz
$n\times n$ definida da seguinte forma
$$
A_{ij} = \left\{
\begin{array}{cc}
1 & \mbox{se } i \mbox{ e } j  \mbox{ são vizinhos;} \\
0 & \mbox{ caso contrário. }
\end{array}
\right.
$$
Observe que por construção, uma matriz de adjacência é uma matriz simétrica.
Em termos dessa matriz, o grau $g_i$, do vértice $i$ de $G$ pode ser escrito como
$$
g_i = \sum_{j=1}^n{A_{ij}}.
$$

A \textbf{matriz laplaciana} de $G$, denotada por ${L}_G$ ou simplesmente 
${L}=[L_{ij}]_{i,j=1}^n$, 
é uma matriz
$n\times n$ definida como segue:
$$
L_{ij} = \left\{
\begin{array}{cc}
g_{i} & \mbox{se } i=j ;\\
-1 & \mbox{se } i \mbox{ e } j  \mbox{ são vizinhos;} \\
0 & \mbox{ caso contrário. }
\end{array}
\right.
$$

%
%
%
%

Seja ${D}_G$ a matriz diagonal formada pelos graus dos vértices de $G$. Então podemos reescrever
a matriz laplaciana na forma
$$
{L}_G = {D}_G -{A}_G.
$$
De forma mais precisa, note que os elementos da matriz laplaciana e da matriz de adjacência estão
relacionados da seguinte maneira
\begin{equation}\label{eq:lijaij}
L_{ij} = \delta_{ij}g_i - A_{ij}
\end{equation}
onde o $\delta_{ij}$ é o delta de Kronecker.

A matriz laplaciana, por construção, também é uma matriz simétrica. Na seção a seguir, discorremos
algumas propriedades relacionadas à essa importante matriz.

\section{Propriedades Espectrais do Laplaciano}

Resumimos as principais propriedades relacionadas ao espectro do laplaciano no

\begin{teorema}\label{espectrolaplaciano}
Seja $G$ um grafo simples e não-direcionado e ${L}$ sua respectiva matriz laplaciana. Então

\begin{enumerate}
\item[(a)] Todos os autovalores de ${L}$ são reais.
\item[(b)] $0$ (zero) é autovalor de ${L}$ associado ao autovetor
$\mathbf{c}=(c,c,\cdots,c), c\in\mathbb{R}\setminus\{0\}$ o qual é único. Em particular,
podemos considerar o autovetor $\mathbf{1} = (1,1,\cdots,1)$.
\item[(c)] ${L}$ é positiva semi-definida e seus autovalores podem ser ordenados,
de acordo com suas multiplicidades algébricas, da forma
$$
0=\lambda_1\leq\lambda_2\leq\cdots\leq\lambda_n.
$$
\item[(d)] A multiplicidade do autovalor $0$ é igual ao número de componentes
conexas de $G$.
\end{enumerate}
\end{teorema}

\noindent
\textbf{Prova:} 

\textit{(a)} Sabemos que a matriz ${L}$ é simétrica, ou seja, ${L}={L}^\dagger$. 
Considere $\lambda$
um autovalor de ${L}$ e ${v}$ o autovetor associado à $\lambda$. Assim,
$
{L}{v} = \lambda {v},
$
donde segue que ${v}^\dagger {L}{v} = 
\lambda {v}^\dagger {v}$. Note que $({v}^\dagger {L}{v})^\dagger 
= {v}^\dagger {L}{v}$. Logo,
$
{v}^\dagger {L}{v} = \lambda ^\dagger v^\dagger {v}
=\lambda {v}^\dagger {v},
$
ou seja, $\lambda ^\dagger = \lambda$.

\textit{(b)} Note que
$$
L\mathbf{c} = \left( \sum_{j=1}^n L_{1j} c, \cdots, \sum_{j=1}^n L_{nj} c\right) = 
(0,\cdots,0) = 0\mathbf{c}.
$$
Por outro lado, se ${L}{v}=\mathbf{0}$ então
$$
\sum_{j=1}^n{L_{1j}v_j} =0, \quad \sum_{j=1}^n{L_{2j}v_j} =0, \quad \cdots, 
\quad \sum_{j=1}^n{L_{nj}v_j} =0
$$
mas
$$
\sum_j{L_{1j}v_j}  =L_{11}v_1 + \sum_{j=2}^n{L_{1j}v_j}
= -\sum_{j=2}^n{L_{1j}v_1} + \sum_{j=2}^n{L_{1j}v_j}
= \sum_{j=1}^n{L_{1j}(v_j-v_1)}=0,
$$
assim, seguindo o mesmo argumento temos que
$$
\sum_j{L_{1j}(v_j-v_1)} =0, \quad \sum_j{L_{2j}(v_j-v_2)} =0, 
\quad \cdots, \quad \sum_j{L_{nj}(v_j-v_n)} =0
$$
o que produz
$$
\sum_j{L_{1j}(v_j-v_1)} + \sum_j{L_{2j}(v_j-v_2)} +
\cdots + \sum_j{L_{nj}(v_j-v_n)} 
 = \sum_i{\sum_j{L_{ij}(v_j-v_i)}} =0
$$
implicando $v_j=v_i \quad \forall i,j =1,\cdots,n$.

\textit{(c)} Mostraremos que ${L}$ é positiva semi-definida. Seja ${x}
=(x_1,\cdots,x_n)\in\mathbb{R}^n$. Então,
queremos mostrar que ${x}^\dagger {L} {x}\geq0$.
Note que
\begin{align*}
{x}^\dagger {L} {x} & = \sum_i{x_i\sum_j{L_{ij}x_j}}\\
& = \sum_i{x_i\sum_j{(D_{ij}-A_{ij})x_j}}
= \sum_i{x_i\sum_j{D_{ij}x_j}} - \sum_i{x_i\sum_j{A_{ij}x_j}}
\end{align*}
mas
$$
\sum_i{x_i\sum_j{D_{ij}x_j}} = \sum_i{g_i {x_i}^2}
 = \sum_i{\sum_j{A_{ij}{x_i}^2}}.
$$
Então
$$
{x}^\dagger L {x} = \sum_i{\sum_j{A_{ij}{x_i}^2}} - \sum_i{\sum_j{A_{ij}{x_i x_j}}},
$$
mas como $A$ é simétrica tem-se que
$$
\sum_i{\sum_j{A_{ij}{x_i}^2}} =\sum_i{\sum_j{A_{ji}{x_j}^2}} 
 = \sum_i{\sum_j{A_{ij}{x_j}^2}},
$$
assim,
\begin{align*}
2{x}^\dagger {L} {x} & =
\sum_i{\sum_j{A_{ij}{x_i}^2}} + \sum_i{\sum_j{A_{ij}{x_j}^2}} - 2\sum_i{\sum_j{A_{ij}x_i x_j}} \\
& = \sum_i{\sum_j{A_{ij}( {x_i}^2 - 2x_i x_j +{x_j}^2)}}
= \sum_i{\sum_j{A_{ij}{({x_i} - {x_j})^2}}}\geq0
\end{align*}
portanto,
$$
{x}^\dagger {L} {x}\geq0.
$$
Como ${L}$ é positiva semi-definida, segue que seus autovalores são não-negativos e portanto podemos
ordená-los da forma
$$
0=\lambda_1\leq\lambda_2\leq\cdots\leq\lambda_n.
$$

Para provar a
afirmação \textit{(d)}, devemos notar que se o grafo $G$ tem $r$ componentes conexas, sendo $G_{1}, \dots ,G_{r}$ 
suas componentes, então é possível representar ${L}$ como uma divisão em 
blocos ${L}_{1}, \dots , {L}_{r}$. Seja $m$ a multiplicidade algébrica do $0$. 
Então cada ${L}_{i}$ tem um autovetor $\mathbf{c}_{i}$ associado ao autovalor $0$. 
Note que $\mathbf{c}_{i}=(c_{i1},c_{i2}\dots ,c_{in})$ pode ser definida como
$$
c_{ij} = \left\{
\begin{array}{cc}
1 & \mbox{se } j \mbox{ pertence a componente } i \\
0 & \mbox{ caso contrário },
\end{array}
\right.
$$
portanto $m\geq r$. Como o autovetores $\mathbf{c}_i$ são únicos então $m=r$.
$\blacksquare$

Analisar o espectro do Laplaciano é de grande importância ao cumprimento do
objetivo de estudar a estabilidade da sincronização em redes. 
De forma mais precisa, sabemos que o segundo autovalor do laplaciano
está intimamente relacionado ao valor do parâmetro crítico de acoplamento entre os
elementos da rede para a sincronização da mesma e sua respectiva estabilidade.

\section{Espectro de Alguns Grafos Regulares}

Vamos agora examinar os autovalores e os autovetores do Laplaciano de 
grafos com estrutura regular de conexão. 
De forma mais precisa, vamos examinar
\begin{enumerate}
\item O grafo \textit{complete} de $n$ vértices, $K_n$, tal que 
$E(K_n)=\{ (u,v) : u\neq v\}$;
\item O grafo \textit{star} de $n$ vértices, $S_n$, tal que 
$E(S_n)=\{ (1,v) : 2\leq v\leq n\}$;
\item O grafo \textit{path} de $n$ vértices, $P_n$, tal que
$E(P_n)=\{ (u,u+1) : 1\leq u< n\}$ e
\item O grafo \textit{ring} de $n$ vértices, $R_n$, que além de ter 
todas as arestas do grafo
\textit{path} possui também a aresta $(1,n)$.
\end{enumerate}

\begin{proposicao}
O Laplaciano de $K_n$ tem autovalor $0$ (zero) com multiplicidade de $1$ e autovalor $n$ com multiplicidade
$n-1$.
\end{proposicao}

\noindent
\textbf{Prova:} Como o grafo $K_n$ é conexo, segue pelo Teorema \ref{espectrolaplaciano} \textit{(d)}
que o autovalor $0$ (zero) tem multiplicidade $1$.
Seja ${v}=(v(1),\cdots,v(n))$ um vetor qualquer, não-nulo, 
ortogonal ao vetor $\mathbf{1}=(1,\cdots,1)$, então
$$
\sum_i v(i) = 0.
$$
Assuma sem perder a generalidade que $v(1)\neq0$. Vamos calcular a primeira coordenada de
${L}_{K_n}{v}$. Ficamos com
\begin{align*}
\left({L}_{K_n}{v}\right)(1) & = g_1 v(1) - \sum_{i=2}^n v(i)\\
& = (n-1)v(1) - \sum_{i=2}^n v(i)
= n v(1),
\end{align*}
onde $g_1=\cdots =g_n = n-1$ é o grau do vértice 1.
Estendendo a análise para todas as coordenadas de ${v}$ temos que
${L}_{K_n}{v} = n{v}$.
$\blacksquare$

\begin{lema}\label{lemasn}
Seja $G$ um grafo conexo de $n$ vértices e sejam 
$i,j \in \{2,\cdots,n\}$ vértices de grau $1$ tais que ambos estão conectados à um outro vértice
$k$. Então o vetor ${v}(u)$ dado por
$$
{v}(u)=
\left \{
\begin{matrix}
1 & para & u=i;\\
-1 & para & u=j;\\
0 & \mbox{caso contrário,} & { }
\end{matrix}
\right.
$$
é um autovetor do Laplaciano de $G$ com autovalor $1$.
\end{lema}



\begin{proposicao}
O grafo $S_n$ tem autovalor $0$ (zero) com multiplicidade $1$, autovalor $1$ com multiplicidade 
$n-2$ e autovalor
$n$ com multiplicidade $1$.
\end{proposicao}

\noindent
\textbf{Prova:} Como $S_n$ é um grafo conexo então pelo Teorema \ref{espectrolaplaciano} $(d)$
segue que o autovalor $0$ tem multiplicidade $1$. Aplicando o Lema \ref{lemasn} aos 
vértices $i$ e $i+1$ encontramos $n-2$ autovetores linearmente independentes. 
Afim de encontrar o último autovalor utilizamos o fato de que o traço de uma matriz é
igual à soma de seus autovalores. 
Notamos que $\mbox{tr}({L}_{S_n}) = 2n-2$. Como identificamos
$n-1$ autovalores cuja soma é $n-2$, então o último autovalor é igual a $n$.
$\blacksquare$

Antes de enunciar a próxima proposição vamos discorrer algumas palavras sobre o Laplaciano
de $R_n$. De forma explícita o mesmo se lê
$$
{L}_{R_n} =
\left[
\begin{matrix}
2 & -1 & 0 & \cdots & 0 & -1\\
-1 & 2 & -1 & 0 & \cdots & 0\\
0 & -1 & 2 & -1 & \cdots & 0\\
\vdots & \mbox{ } & \mbox{ } & \ddots & \mbox{ } & \vdots\\
0 & 0 & \cdots & -1 & 2 & -1\\
-1 & 0 & \cdots & 0 & -1 & 2
\end{matrix}
\right]_{n\times n}.
$$
Dado um vetor de ${v}=(v_1,v_2,\cdots,v_n)\in\mathbb{R}^n$, notamos que
a ação de ${L}_{R_n}$ sobre ${v}$ produz 
$$
\left({L}_{R_n}{v}\right)_k = -v_{k+1} - v_{k-1} + 2v_k,
$$
para $1\leq k\leq n$.
Note ainda que dada uma função real $f(t)$ de classe $C^2$, a sua segunda derivada pode ser escrita da forma
$$
f''(t) = \lim_{h\rightarrow 0}\frac{f(t+h) + f(t-h) - 2f(t)}{h^2}.
$$
Sendo assim,vemos que o resultado da ação de ${L}_{R_n}$ sobre ${v}$ se assemelha a 
segunda derivada discreta, escolhendo $h=1$. As funções reais tais que suas segunda derivada
é um múltiplo de si mesmo são os senos, cossenos e exponenciais, que por sua vez pode
ser expressa em função de senos e cossenos. 
Podemos então enunciar a seguinte

\begin{proposicao}\label{espectrorn}
O grafo $R_n$ tem autovetores
$$
\begin{matrix}
{x}_k(u) = \sin\left(\frac{2ku\pi }{n}\right)\\
\\
{y}_k(u) = \cos\left(\frac{2ku\pi}{n}\right)
\end{matrix}
$$
com $0\leq k\leq n/2$ se $n$ é par e $0\leq k\leq (n-1)/2$ se $n$ é ímpar e $1\leq u \leq n$. 
Ambos os autovetores têm autovalor $2-2\cos(2k\pi /n)$. O vetor
${x}_0(u) = \mathbf{0}$ deve ser desconsiderado e se $n$ é par então da mesma forma o vetor
${x}_{n/2}(u)=\mathbf{0}$ deve ser desconsiderado. Note que 
${y}_0(u) =\mathbf{1}$ é o autovetor
constante.
\end{proposicao}




\noindent
\textbf{Prova:} Para encontrar o autovalor 
vamos considerar apenas
o vértice $1$, ou seja, $u=1$ sem perder a generalidade. Então, 
sabendo que $\sin(2\alpha) = 2\sin(\alpha)\cos(\alpha)$ e que
$\cos(2\alpha)= 2\cos^2(\alpha) -1$, temos
\begin{align*}
{L}_{R_n}{x}_k (1) & = 2{x}_k(1) - {x}_k(0) - 2{x}_k(2)\\
& = 2\sin{(2k\pi/n)} - 2\sin{(2k\pi/n)}\cos{(2k\pi/n)}\\
& = \left[ 2-2\cos\left(\frac{2k\pi}{n}\right) \right]{x}_k (1).
\end{align*}
Da mesma forma
\begin{align*}
{L}_{R_n}{y}_k (1) & = 2{y}_k(1) - {y}_k(0) - 2{y}_k(2)\\
& = 2\cos{(2k\pi/n)} - 1 - \cos{(2\cdot 2k \pi/n)}\\
& = 2\cos{(2k\pi/n)} - 1 - \cos^2{(2k \pi/n)}-\sin^2{(2k \pi/n)}\\
& = \left[ 2-2\cos\left(\frac{2k\pi}{n}\right) \right]\cos{(2k\pi/n)} =
\left[ 2-2\cos\left(\frac{2k\pi}{n}\right) \right]{y}_k (1).
\end{align*}
$\blacksquare$

Antes de enunciar o próximo resultado vamos verificar uma importante relação
entre os grafos $P_n$ e $R_{2n}$.
Vamos olhar para $P_n$ como um quociente de $R_{2n}$ pela seguinte relação
de equivalência: Dados $u$ e $v$ vértices de $R_{2n}$ temos que $u$ é equivalente a $v$ se
e somente se $u+v = 2n+1$. Dessa forma, por exemplo considerando $n=5$, ou seja,
considerando os grafos $P_5$ e $R_{10}$, em $R_{10}$ estamos identificando os vértices
$1$ e $10$, $2$ e $9$ e assim por diante, tornando o quociente de ${R_{2n}}$ pela
relação de equivalência dada igual a $P_n$.

Argumentos análogos aos utilizados na Proposição \ref{espectrorn} nos mostram que

\begin{proposicao}
O Laplaciano de $P_n$ tem os mesmos autovalores de $R_{2n}$ e autovetores
$$
{v}_k(u) = \sin (ku\pi/2n+\pi/2n),
$$
com $0\leq k \leq n$.
\end{proposicao}

\begin{teorema}\label{teo:lambda2cotas}
Seja $G$ um grafo simples de $n$ vértices, $d$ o seu diâmetro, $g_1$ o menor grau
entre todos os vértices e $\lambda_2$ o segundo autovalor do laplaciano. Então
 \begin{enumerate}
  \item $\lambda_2 \geq \dfrac{4}{nd}$
\item $\lambda_2 \leq \dfrac{ng_1}{n-1}$
 \end{enumerate}
\end{teorema}
As provas para estas desigualdades podem ser encontradas em \cite{bojan} e \cite{fiedler} 
respectivamente.
Agrupamos os principais resultados da análise espectral do Laplaciano na Tabela \ref{tab:laplacian}, 
onde $g_1$ e $g_n$ representam respectivamente o menor e o maior grau no grafo.

\begin{table}[h!]
\centering
\begin{tabular}{ c | c | c | c | c }
\hline Grafo & $\lambda_{2}$ & $g_{n}$ & $g_{1}$ & $d$ \\ 
\hline                        
\hline $K_n$ & $n$ & $n-1$ & $n-1$ & $1$ \\
\hline $R_n$ & ${2-2\cos\left(\frac{2\pi}{n}\right)}$ & $2$ & $2$ & 
$
\begin{array}{c}
(n+1)/{2} \mbox{ se $n$ é ímpar}\\
n/2 \mbox{ se $n$ é par}
\end{array}
$\\
\hline $S_n$ & 
$
\begin{array}{c}
{2} \mbox{ se $n=2$}\\
1 \mbox{ se $n>2$}
\end{array}
$ 
& $n-1$ & $1$ & $2$ \\
\hline $P_n$ & ${2-2\cos\left(\frac{\pi}{n}\right)}$ & 1 & 1 & $n-1$\\
\hline
\end{tabular}
\caption{Propriedades espectrais do Laplaciano de algumas redes regulares.}\label{tab:laplacian}
\end{table}

\section{Redes Complexas}\label{sec:redescomplexas}

Nesta seção abordaremos alguns exemplos de redes complexas, como elas são formandas, suas
propriedades e aplicações.
Inicialmente daremos a definição formal desse termo.

\begin{definicao}
Uma rede, modelada a partir de um grafo G, é dita ser \textbf{complexa} se G não possui
uma estrutura regular de conectividade.
\end{definicao}

Os exemplos que queremos abordar de redes complexas estão separados em subseções como segue.

\subsection{Redes Aleatórias}

As redes aleatórias, ou redes que seguem o modelo Erdös-Rényi, tem a seguinte característica
de geração: para um dado $p$ fixo, $0\leq p \leq 1$, cada aresta em potencial é escolhida com
probabilidade $p$, independentemente das outras arestas. Dessa forma, todos os vértices de
um grafo aleatório tem a mesma esperança para o grau. Nesse modelo
de rede, podem existir vértices isolados, ou seja, a rede pode não ser efetivamente conexa.
Quando o número de vértices $n\rightarrow\infty$ então a probabilidade de que a rede seja conexa
tende a $1$ \cite{chung2006complex}.

\subsection{Redes Pequeno Mundo}

O termo fenômeno de pequeno mundo é usado para se referir à duas propriedades distintas,
a saber, a propriedade de pequena distância (dois estranhos estão tipicamente ligados por uma
pequena cadeia de conhecidos mútuos) e a propriedade de efeito de agrupamento (duas pessoas que
compartilham o mesmo vizinho possuem uma maior probabilidade de também serem vizinhos). Esse
tipo de rede foi introduzido por Watts e Strogatz em \cite{wattsandstrogatz}. Alguns exemplos
desse tipo de rede complexa são as redes neurais, redes elétricas,  redes de
co-atuação de atores em filmes americanos, dentre outras. As redes de
pequeno mundo encontram-se em um meio termo entre as redes regulares (não complexas) e
redes aleatórias. Portanto, uma forma de gerar esse tipo de rede é considerando
inicialmente uma rede regular tipo \textit{ring} com $n$ vértices e cada vértice com grau $g$, em
seguida reconecta-se cada aresta de forma aleatória com probabilidade $p$, onde $p=0$ gera a rede
regular descrita e $p=1$ uma rede aleatória. Essas redes apresentam, significante aumento na
velocidade de propagação de um sinal e consequente sincronizabilidade. 

\subsection{Redes \textit{Scale-Free}}

Nem as redes aleatórias ou as de pequeno mundo tem uma propriedade frequentemente observada
em redes do mundo real, a saber, o comportamento de que os vértices possuem probabilidade
de ter grau $g$ seguindo a lei de potência
$$
P(g) \approx g^{-\beta}.
$$

Como sabemos, o grau de um vértice é o número de vizinhos que o mesmo possui. A lei de
potência assegura que o número de vértices com grau $g$ é proporcional a $g^{-\beta}$ com 
$\beta\geq1$. As redes que seguem essa lei de distribuição de potência são conhecidas como
\textit{scale-free networks} \cite{dehmer2010structural}. Este modelo de rede é devido a
Barabasi e Albert \cite{barabasiandalbert1999}. O termo scale está relacionado tanto ao
espaço quanto ao tempo e de fato, escalas de espaço e tempo podem coexistir simultaneamente.
Por exemplo, as redes de chamadas telefônicas tem formas muito similares mesmo em diferentes
regiões geográficas e em diferentes horários \cite{chung2006complex}. 
As redes que seguem o modelo Barabasi-Albert
possuem uma grande quantidade de heterogeneidade, ou seja, enquanto a maioria dos vértices
possuem apenas alguns vizinhos, alguns poucos vértices, chamados de \textit{hubs} possuem
muitos vizinhos.
Um exemplo claro desse tipo de rede é a rede de estradas que conectam as cidades de um país.
As grandes cidades possuem muitas alternativas de estradas que estão ligadas à ela, e as
pequenas cidades poucas alternativas de ligação com outras cidades, porém o número de cidades
pequenas é muito maior do que o de mega-cidades.
Por ter essa característica de heterogeneidade, esse tipo de rede complexa, pode não ter
sincronização global, isto é, neste caso apenas os \textit{hubs} podem apresentar sincronização
porém os demais elementos da rede ficarem fora desse estado \cite{tiago}.



\chapter{Redes de Osciladores com Acoplamento Difusivo}\label{cap:osciladores}

O comportamento oscilatório desempenha um papel importante na natureza. Toda forma
de vida exibe tal comportamento em cada nível de organização biológica com períodos
que podem variar desde milissegundos, como é o caso dos neurônios, a anos, como é
o caso da interação presa-predador na ecologia \cite{pogrosmky1999}.
%
Considere um grafo $G$ que modela uma certa rede. Dizemos então
que os elementos da rede (vértices de $G$) possuem um comportamento oscilatório, ou 
chamamos os elementos de {osciladores}, quando os mesmos possuem alguma dinâmica.
%
Um oscilador pode possuir qualquer tipo de dinâmica, isto é, pode ter um comportamento
periódico, não-periódico ou caótico. No nosso caso, para qualquer que seja
o tipo de comportamento que o oscilador possuir, representamos esse comportamento 
através de equações diferenciais ordinárias.

A presente dissertação trata apenas o caso de comportamento oscilatório quando as
equações diferenciais ordinárias são todas idênticas, em outras palavras, do ponto
de vista isolado, os osciladores são todos idênticos, a menos das condições iniciais.

Neste momento, podemos introduzir o nosso modelo de rede. Considere que um grafo
$G$ modela uma rede qualquer. Por padrão,
consideramos sempre que $n$ é o número de vértices em $G$. 
Então, para que possamos introduzir o comportamento
oscilatório, sobre cada vértice da rede introduzimos uma cópia da Equação \eqref{x=ftx}.
A regra de nomeação dos vértices será a seguinte: 
considerando que os graus dos vértices de $G$ podem ser dispostos da forma
\begin{equation}
g_1\leq g_2\leq\cdots\leq g_n,
\end{equation}
então o vértice $i$ será aquele que possui grau $g_i$.
O modelo de interação entre os elementos de uma rede, que utilizamos,
é o modelo de acoplamento chamado difusivo:
\begin{definicao}
Dizemos que uma rede possui modelo de \textbf{acoplamento difusivo} 
quando a dinâmica do vértices $i$ é influenciada por seu vizinho $j$ de forma
proporcional à uma função da diferença de seus estados.
\end{definicao}
Chamaremos essa função acoplamento de $H:\mathbb{R}^m\rightarrow\mathbb{R}^m$ e sem perder 
a generalidade vamos considerar que a mesma
é uma matriz em $\Mat(\mathbb{R},m)$. A função de acoplamento
poderia ter até mesmo um comportamento não-linear, porém por toda a dissertação,
vamos considerar apenas
o caso em que $H$ é uma matriz positiva-definida. Se supormos que os vértices $i$ e $j$ 
são vizinhos, então a influência que 
o vértice $j$ exerce sobre seu vizinho $i$ é igual
$$H(x_j) - H(x_i),$$ de modo que a análise é recíproca.
Em muitos exemplos vamos considerar que $H=I_m$, nesse caso esse
modelo de acoplamento é chamado de acoplamento totalmente difusivo.

Sendo assim, considerando uma rede de $n$ osciladores idênticos, com acoplamento
difusivo, tem-se que a dinâmica de um vértice $i$ qualquer da rede passa a ser descrita por 
\begin{equation}\label{eq:mdoeloacoplamentoadsj}
\dot{x}_i  = f(x_i) + \alpha \sum_{j=1}^n{A_{ij}[H(x_j) - H(x_i)]}
\end{equation}
onde $\alpha$ é o parâmetro global de acoplamento, ou seja, para todo
vértice $i=1,\cdots,n$ temos que $\alpha$ é o mesmo valor.
Note que a equação \eqref{eq:mdoeloacoplamentoadsj} está em função dos elementos
da matriz de adjacência do grafo e note ainda que o modelo está bem definido 
pois o mesmo considera que a dinâmica do vértice $i$
é influenciada apenas por seus vizinhos e de fato, $A_{ij}=0$ se $j$ não é vizinho de $i$.

Nos podemos também representar o nosso modelo através da matriz laplaciana do grafo, pois se
considerarmos apenas o termo de acoplamento e lembrando que $L_{ij} = \delta_{ij}g_i-A_{ij}$
(Equação \eqref{eq:lijaij}) tem-se
\begin{align*}
\sum_{j}{A_{ij}[H(x_j) - H(x_i)]} & =  \sum_j A_{ij}H(x_j) - g_i H(x_i)\\
& = -\sum_j L_{ij} H(x_j).
\end{align*}
Sendo assim, reescrevemos a dinâmica do vértice
$i$ da forma
\begin{equation}\label{eq:laplacianmodel}
\dot{x}_i  = f(x_i) - \alpha \sum_{j=1}^n{L_{ij}H(x_j)}.
\end{equation}

O modelo de acoplamento difusivo desempenha um papel importante para a sincronização devido
exatamente à sua natureza difusiva, ou seja, a tentativa de igualar os estados dos vértices da rede, 
influenciando assim positivamente a rede à um estado síncrono.
Observe que se todos os elementos da rede possuem a mesma condição inicial, então os termos
de acoplamento irão desaparecer identicamente. Além disso, se a rede em algum momento $t_0$
está em um estado de sincronização, então é possível mostrar que a variedade de sincronzação, 
ou seja
\begin{equation}\label{eq:estadosincrono}
x_1(t) = x_2(t)  = \cdots = x_n(t) = s(t)
\end{equation}
para todo $t\geq t_0$ é uma variedade invariante no tempo para qualquer que seja o parâmetro
de acoplamento $\alpha > 0$ e para qualquer que seja a escolha da função de acoplamento $H$.
Queremos caracterizar as sincronizações que são persistentes e portanto introduzimos formalmente
o termo que é essencial nesta dissertação.
\begin{definicao}[\textbf{Sincronização}]
Dizemos que uma rede está sincronizada se existe $\delta>0$ tal que para todo $\epsilon>0$ existe
$\tau = \tau(\epsilon) >0$ satisfazendo $\Vert x_i(u) - x_j(u)\Vert \leq \delta$, para algum
$u\geq0$, e $\Vert x_i(t) - x_j(t)\Vert \leq \epsilon$ para todo $t\geq u +\tau$, e para quaisquer
$i$ e $j$.
\end{definicao}

Esta noção de sincronização produz uma variedade de sincronização que é uniformemente assintoticamente
estável, o que nos garante a persistência da sincronização sob perturbações.

O termo sincronização global, no nosso contexto, está relacionado 
à dois aspectos, a saber, global no sentido de que a sincronização não é local (o que pode acontecer, 
por exemplo em rede complexas tipo \textit{scale-free}), e global no sentido de que a coincidência
dos estados (Equação \eqref{eq:estadosincrono}) é preservada com o tempo $t$. Formaliza-se essa
propriedade no Teorema \ref{teo:variedadeinvariante}.
De modo mais preciso, vamos mostrar a invariância do movimento
globalmente sincronizado \eqref{eq:estadosincrono}, para isso precisamos reescrever o 
nosso modelo \eqref{eq:laplacianmodel}
em uma forma compacta onde possamos agrupar as equações do movimento de todos os vértices,
de forma que passemos a visualizar apenas uma equação de movimento em $\mathbb{R}^{nm}$ e não $n$
equações em $\mathbb{R}^m$.
Considere então $X=(x_1,x_2,\cdots,x_n)\in\mathbb{R}^{nm}$
onde estamos considerando o empilhamento dos vetores $x_1, x_2, \cdots,x_n\in\mathbb{R}^m$. Da mesma
forma considere
$F(X) = (f(x_1),f(x_2),\cdots,f(x_n)).$
Então o modelo \eqref{eq:laplacianmodel} pode ser reescrito da forma

\begin{equation}\label{eq:modelocompacto}
\dot{X} = F(X) - \alpha (L\otimes H)X
\end{equation}
onde $\otimes$ representa o produto de Kronecker. 

Sejam  $\Phi_t(\cdot)$ é o fluxo da equação \eqref{eq:modelocompacto} e 
\begin{equation}\label{eq:variedadedesincronizacao}
{N} = \{ \mathbf{1}\otimes s(t) \in \mathbb{R}^{n\times m}  : 
s(t)=x_i(t) \in \mathbb{R}^m, \mbox{ } \forall \mbox{ } 1\leq i\leq n  \}
\end{equation}
a variedade em $\mathbb{R}^{n\times m}$ que representa o estado síncrono global. Enunciamos
então o

\begin{teorema}\label{teo:variedadeinvariante}
${N}$ é uma variedade invariante pelo fluxo.
\end{teorema}

\noindent
\textbf{Prova:} 
Vamos considerar apenas o caso em que $N$ é positivamente invariante pelo fluxo.
Dada a condição inicial $X_0\in N$, isto é, $X_0 = \um\otimes s(0)$, afirmamos que $X(t) = 
\um\otimes s(t)$, $t\geq0$, é solução de \eqref{eq:modelocompacto}. De fato pois
\begin{align*}
\dot{X}(t) & = F(X(t)) - \alpha (L\otimes H)X(t)\\
& = F(\mathbf{1}\otimes s(t)) - \alpha (L\otimes H)\left(\mathbf{1}\otimes s(t)\right)\\
& = \mathbf{1}\otimes f(s(t)) - \alpha (L \mathbf{1})\otimes H s(t)
= \mathbf{1}\otimes f(s(t))
\end{align*}
visto que $L \mathbf{1} = \mathbf{0}$, $F(\mathbf{1}\otimes s(t))=\mathbf{1}\otimes f(s(t))$ e
$\dot{X}(t) = \mathbf{1}\otimes \dot{s}(t)$. Além disso, se $x_i(t) = s(t)$ para todo $i$, de 
\eqref{eq:mdoeloacoplamentoadsj} segue que $\dot{s} = f(s)$.
Dessa forma,
o fluxo pode ser escrito como $\Phi_t(\mathbf{1}\otimes s(0) ) = X(t) = \mathbf{1}\otimes s(t)$,
logo $\Phi_t(X_0 ) \in {N}$. $\blacksquare$
%
\\

Esse teorema nos fornece uma forte garantia sobre o modelo que estamos utilizando, ou seja,
ao considerar que, isoladamente, os osciladores que possuem a mesma dinâmica. Se
os osciladores tem as mesmas condições iniciais, então os mesmos estarão
automaticamente sincronizados e assim permanecerão por todo o tempo futuro. Se por
outro lado, os osciladores não possuem as mesmas condições iniciais, mas para algum $t>0$,
o parâmetro global de acoplamento garante que os osciladores irão sincronizar, então
da mesma forma, sincronizados eles estarão por todo o tempo futuro.
Para que esse cenário seja coerente, estamos considerando que não existem influências
de pertubações sobre a variedade de sincronização. Porém, na prática,
essa coincidência de estados só estará garantida se tal regime é estável, ou seja, se
esse subespaço permanece inalterado após uma pequena pertubação \cite{rosenblum03}.

\section{Parâmetro Global de Acoplamento}

Analisaremos agora o quão importante é o parâmetro de acoplamento $\alpha$ para a
sincronização e sua estabilidade.

Considere que o campo de vetores $f$ da dinâmica isolada dos osciladores está associado
à uma dinâmica caótica. Considere também que $\alpha = 0$. Se impormos que as condições
iniciais dos osciladores são próximas, porém distintas, teremos que as trajetórias associadas
aos vértices irão divergir exponencialmente rápido até o ponto em que as diferenças
$\Vert x_j(t) - x_i(t)\Vert$ para todo $i,j=1,\cdots,n$ com $i\neq j$ serão tão grandes quanto
o próprio diâmetro do atrator. Em outras palavras, se considerarmos a rede desacoplada
($\alpha = 0$) e condições iniciais distintas, então não haverá
sincronização entre os osciladores. Além disso, se considerarmos agora que a rede
está desacoplada mas todos os osciladores possuem as mesmas condições inciais, teremos que
a variedade de sincronização, sob qualquer pequena perturbação, também irá começar a divergir
exponencialmente rápido e dessa forma levar os vértices a terem comportamentos
distintos em suas dinâmicas, ou seja, neste caso, a variedade de sincronização perde
sua estabilidade sob qualquer perturbação. Não consideramos $\alpha<0$, isso por que
neste caso, as soluções de \eqref{eq:mdoeloacoplamentoadsj} podem não ser limitadas, e portanto
não faria sentido falar de sincronização. Considere por exemplo $n=2$, $f=\mathbf{0}$ e $H=I_m$. 
Neste contexto, a diferença $x_1(t)-x_2(t)=z(t)$ cumpre
$$
z(t) = e^{-2\alpha t}z(0),
$$
de forma que a mesma só irá decrescer exponencialmente em $t$ se $\alpha>0$. Portanto,
consideramos sempre $\alpha>0$.

\section{Principais Resultados}

A questão da estabilidade da variedade de sincronização ${N}$ é um dos principais
objetivos do presente trabalho. Portanto, os nossos principais resultados estão relacionados
à mesma. Tais resultados são três, e chamá-los-emos de Teorema da Existência Global das Soluções,
o qual está relacionado ao modelo \ref{eq:laplacianmodel}, Teorema da Estabilidade
da Variedade de Sincronização, que garante a sincronização e sua consequente estabilidade e finalmente,
Teorema da Persistência da Sincronização, que segue praticamente como um corolário do Teorema
do Parâmetro Crítico de Acoplamento. Enunciamos os mesmos no presente capítulo, porém as suas 
respectivas provas serão dadas no Capítulo \ref{cap:estabilidadedassolucoessincronas}.

A construção da função de Lyapunov, estabelecida no capítulo \ref{cap:dinamicanaolinear},
é usada para garantir a existência global das soluções, tanto no contexto das equações desacopladas
quanto e principalmente para os oscilados em interação - estamos nos referindo ao modelo 
\eqref{eq:laplacianmodel}. Isso porque para tratar de sincronização dos osciladores, primeiramente
precisamos garantir a existência global das soluções e que as mesmas sejam limitadas. Portanto, 
independentemente da rede considerada, a hipótese de que a dinâmica isolada tem uma função de
Lyapunov garante que a dinâmica coletiva também cumpre a mesma propriedade. Enunciamos então o

\begin{teorema}[\textbf{da Existência Global das Soluções}]\label{teo:egscaposc}
Considere o modelo de acoplamento difusivo
\begin{equation}\label{eq:modeloredecaposc}
\dot{x}_i = f(x_i) - \alpha\sum_{j=1}^n{L_{ij}H(x_j)}.
\end{equation}
Assuma que a dinâmica isolada
tem função de Lyapunov satisfazendo a Suposição \ref{sup:atrator}.
Então, independentemente da rede, as soluções do modelo \eqref{eq:modeloredecaposc} certamente
entram, em tempo finito, em um domínio absorvente $\Omega$, e portanto as mesmas são 
limitas e existem por todo o tempo futuro.
\end{teorema}

A ideia da prova do Teorema da existência global das soluções, que encontra-se
na Seção \ref{sec:egs} do Capítulo \ref{cap:estabilidadedassolucoessincronas}, está baseada na 
hipótese de que a dinâmica isolada, de todos os osciladores, tem função de Lyapunov e consequentemente um domínio absorvente
onde as soluções ``moram'' por todo o tempo. Então, tomamos como principal hipótese
a Suposição \ref{sup:atrator}. A partir então dessa hipótese construiremos uma
função de Lyapunov para a dinâmica coletiva. Essa função para a coletividade é dada pela soma
das funções de Lyapunov para as dinâmicas individuais, mais precisamente, consideramos
na soma apenas as funções de Lyapunov das dinâmicas cujas trajetórias ainda estão fora 
do domínio absorvente. Note que se todas as trajetórias já estão num domínio absorvente então
não há nada a se fazer, pois tal domínio é um conjunto compacto, positivamente invariante, e
portanto a existência global das soluções já estaria automaticamente garantida. Porém,
através dessa função de Lyapunov para a dinâmica coletiva, conseguimos garantir se se todas as $n$
as trajetórias estão fora do domínio absorvente, ou $n-1$ ou $\cdots$ ou $1$, está(ão) fora do domínio
absorvente, a tempo finito todas as trajetórias, coletivamente, isto é, considerando
o acoplamento difusivo, entram num domínio absorvente e lá permanecem por todo o tempo futuro.

Graças a construção da função de Lyapunov para a dinâmica isolada e consequentemente para a
dinâmica coletiva, sabemos onde moram as soluções. Considerando o modelo de acoplamento difusivo
\eqref{eq:laplacianmodel} e sua forma compacta \eqref{eq:modelocompacto}, condicionamos
que a variedade de sincronização é localmente atratora se as trajetórias possuem
condições iniciais próximas, ou seja, para quaisquer $i,j=1,\cdots,n$ com
$\Vert x_i(0) - x_j(0)\Vert\leq\delta$ com $\delta >0$ pequeno, então
$x_i(t)\rightarrow x_j(t) \mbox{ } \forall
i,j$. Essa convergência é determinada pela força de acoplamento $\alpha$ entre os osciladores. Além
disso, o parâmetro crítico que produz essa propriedade advém das contribuições da
dinâmica isolada dos osciladores, de propriedades espectrais da função de acoplamento $H$
e do segundo autovalor do laplaciano da rede.

\begin{teorema}[\textbf{da Estabilidade da Variedade de Sincronização}]\label{teo:principal1}
Considere o modelo de acoplamento difusivo \eqref{eq:modeloredecaposc}. Considere também
a existência global das soluções de acordo com o Teorema \ref{teo:egscaposc}.
Assuma que as condições iniciais estão em uma vizinhança de $\Omega$ e considere
\begin{equation}\label{eq:alphacriticogeralcaposc}
\alpha_c = \frac{\beta}{\lambda_2\mu_1}
\end{equation}
onde $\beta  = \sup_{x\in\Omega} \Vert P^\dagger Df(x)P\Vert_{\infty}$,
de forma que a função de acoplamento $H$ do modelo \eqref{eq:modeloredecaposc}
é uma matriz positiva-definida com a representação espectral $H=PDP^\dagger$, $\lambda_2$ é o segundo
autovalor do laplaciano da rede e $\mu_1$ é o primeiro autovalor de $H$.

Se $\alpha > \alpha_c$, então o estado síncrono é garantido e o mesmo é uniformemente
assintoticamente estável.
\end{teorema}

A prova é apresentada no capítulo \ref{cap:estabilidadedassolucoessincronas} Seção
\ref{sec:sgr}.
Abordamos aqui ideia da mesma. Um dos principais ingredientes da prova é considerar a
análise de uma só equação diferencial em $\mathbb{R}^{nm}$ e não de $n$ diferentes equações
em $\mathbb{R}^m$. Sob as hipóteses de que as matrizes laplaciana e da função
de acoplamento são simétricas, é possível realizar uma mudança linear de coordenadas
para uma nova base ortonormalizada formada pelos autovetores das matrizes anteriormente mencionadas.
Considerando essa mudança linear de coordenadas, um vetor solução
arbitrário é dado por $X = \mathbf{1}\otimes s + U$, onde $s=s(t)$ é uma dada solução síncrona fixada e
o incremento $U$ é tal $X$ está numa vizinhança da variedade de sincronização ${N}$.
Iremos então conseguir condições para que as contribuições dos modos transversais em $X$
convirjam à zero com o tempo, isto é,
que os operadores de evolução dos modos transversais tenham todos contração uniforme. 
As nomenclaturas ``modos normais'' e ``modos transversais'' estão 
relacionadas com os autovetores da matriz laplaciana da rede. O modo normal é exatamente
o sub-espaço gerado pelo autovetor $\mathbf{1}$.
Como consideramos que a dinâmica individual de cada elemento da rede é essencialmente uma
dinâmica não-linear, utilizamos a expansão em série de Taylor e analisamos a parte linear
do campo de vetores separadamente do resto de Taylor. Após conseguir condições para
que os operadores de evolução dos modos transversais tenham todos contração uniforme, consideramos
a influência do resto de Taylor sobre a análise. Utilizamos diretamente a Proposição
\ref{pro:principio2} que garante que mesmo considerando o resto de Taylor, uma perturbação
não-linear, o estado síncrono não sofre alteração. Para conseguir as condições mencionadas acerca dos
modos transversais, projetamos a equação em blocos nos modos normais
e transversais afim de analisar apenas as equações diferenciais da mudança linear de coordenadas,
que por sua vez são equações diferenciais lineares não-autônomas,
e utilizamos o Teorema da diagonal dominante (\ref{teo:criterio1}) para finalmente obter uma cota para a força de acoplamento
$\alpha$ que por sua vez nos dá o critério enunciado.

O próximo resultado segue praticamente como um corolário do Teorema \ref{teo:principal1}, isso
devido ao tipo de estabilidade que estamos considerando, isto é, tem-se a propriedade de
``robustez'' sobre a variedade invariante de sincronização global. Assim,
devido a forte garantia do Teorema \ref{teo:principal1} podemos introduzir pertubações
no modelo considerado e esperar que a estabilidade da variedade de sincronização não seja 
destruída, isso sob certas hipóteses na magnitude da pertubação.

Vamos trabalhar com um modelo de pertubação que age na função de acoplamento do modelo
de rede \eqref{eq:mdoeloacoplamentoadsj}. Consideramos um modelo de pertubação que age nas funções
de acoplamento pois com algumas manipulações algébricas é possível utilizar precisamente
o Teorema \ref{persistencia}, e então estabelecer a magnitude que as perturbações podem ter 
para que o estado síncrono seja persistente.
Mais precisamente vamos considerar, para cada
oscilador $1\leq i\leq n$, $n$ funções de perturbação $V_{ij}: \mathbb{R}^m\rightarrow\mathbb{R}^m$,
totalizando então $n^2$ possíveis funções de pertubação.
Consideraremos apenas o caso em que cada $V_{ij}$ é um operador linear em $\mathbb{R}^m$, 
os quais podem ou não ter uma dependência temporal, porém a análise é a mesma, e sem perder 
a generalidade, representamos tais perturbações na forma matricial. Considere ainda 
o modelo de rede \eqref{eq:mdoeloacoplamentoadsj} o qual será escrito na forma do modelo de 
perturbação considerado como
\begin{align}
\dot{x}_i & = f(x_i) + \alpha\sum_{j=1}^n A_{ij}  H(x_j-x_i) + \sum_{j=1}^n A_{ij} V_{ij} (x_j-x_i)\\
& = f(x_i) + \sum_{j=1}^n A_{ij}  (\alpha H+V_{ij})(x_j-x_i) \label{eq:mpa}
\end{align}
ou ainda, sabendo que $L_{ij}=\delta_{ij}g_i - A_{ij}$ tem-se
\begin{align}
\dot{x}_i & = f(x_i) + \sum_{j=1}^n (\delta_{ij}g_i-L_{ij}) ( \alpha H+V_{ij})(x_j-x_i)\\
& = f(x_i) - \sum_{j=1}^n L_{ij} (\alpha H+V_{ij})(x_j-x_i).\label{eq:mpl}
\end{align}
E na forma de blocos, a equação \eqref{eq:mpl} é escrita como
\begin{equation}\label{eq:mpb}
\dot{X} = F(X) -   \alpha(L\otimes H)X + 
\left[ \sum_{i,j=1}^n (B_{i}LB_j)( D_{ji}-I_n) \otimes V_{ij} \right]X,
\end{equation}
onde $D_{ij}\in\Mat(\mathbb{R},n)$ é uma matriz cuja $ij$-ésima entrada é igual
a $1$ e todas as demais são iguais $0$ (zero), $A_i\in\Mat(\mathbb{R},n)$ é a matriz diagonal 
tal que a $i$-ésima
entrada da digonal é igual a $1$ e todas as demais são nulas, $L\in\Mat(\mathbb{R},n)$ é o 
laplaciano da rede e $I_n\in\Mat(\mathbb{R},n)$ é a matriz identidade.
Maiores detalhes de como a equação \eqref{eq:mpb} é obtida estão disponíveis no Apêndice 
\ref{apes:blocos}.

Temos então o
\begin{teorema}[\textbf{da Persistência da Sincronização}]\label{teo:persistenciaglobal}
Considere o modelo de acoplamento difusivo \eqref{eq:mdoeloacoplamentoadsj} o qual é reescrito na
na forma de perturbação \eqref{eq:mpa}
e escrito na forma de blocos como descrito em \eqref{eq:mpb}.
Considere que $\alpha> \alpha_c$ \eqref{eq:alphacriticogeralcaposc}. Então a variedade de sincronização
é persistente se
\begin{equation}\label{eq:criterioperturbacao1}
\sup_t \sum_{i,j=1:i\sim j}^n \Vert V_{ij}\Vert_{\infty} < \frac{\eta}{ 2  \Vert L\Vert_{\infty}}
\end{equation}
onde $\eta = \alpha\lambda_2\mu_1-\beta$, e $m$ é a dimensão do sistema isolado \eqref{eq:mpa}.
\end{teorema}

A notação $i\sim j$, significa dizer que o vértice $i$ da rede é vizinho do vértice $j$, ou então
que o oscilador $i$ está acoplado com o oscilador $j$. Antes de
enunciar o teorema falamos que é possível introduzir ao total $n^2$ funções lineares de
perturbação, porém na prática, algumas delas, ou muitas delas, não produzirão qualquer influência
na estabilidade da variedade de sincronização. Precisamente, como consideramos que o modelo de perturbação age na função
de acoplamento entre os osciladores, então só deverá existir perturbação para aqueles osciladores
que estão acoplados.

O principal ingrediente da prova, que estará disponível na Seção \ref{sec:pfa} do Capítulo 
\ref{cap:estabilidadedassolucoessincronas}, é a utilização do Teorema \ref{persistencia}.
Além disso, a prova de referido teorema segue praticamente os mesmos passos da prova do
Teorema \ref{teo:principal1}. Ou seja, utilizamos novamente a teoria de equações diferenciais
lineares não-autônomas abordada no Capítulo \ref{cap:edolinear}. O critério sobre a magnitude
das pertubações é dado então principalmente pelo Teorema \ref{persistencia} e secundariamente pelo
Teorema \ref{teo:criterio1}. O Teorema da persistência da sincronização poderia até mesmo
ser enunciado de uma forma mais generalizada, isto é, generalizada no que se refere
à utilização de uma norma induzida arbitrária substituindo \eqref{eq:criterioperturbacao1} 
por uma desigualdade geral, porém se assim o 
fosse, o seu enunciado envolveria sucessivos e dependentes parâmetros os quais aparecem na 
demonstração do Teorema \ref{teo:principal1}, portanto preferimos enunciar o Teorema 
\ref{teo:persistenciaglobal} utilizando apenas a norma infinito.

Uma das coisas interessantes que a cota \eqref{eq:criterioperturbacao1} nos diz é que
podemos comparar a ``robustez'' de diferentes redes, isso porque
tal cota é dependente e inversamente proporcional a norma infinito da matriz laplaciana da rede.
Então, por exemplo, é de se esperar que redes complexas aleatórias sejam mais ``robustas'' do
que redes complexas \textit{scale-free}, pois notoriamente,
tem-se que $\Vert L\Vert_\infty = 2\max_i g_i$ onde $g_i$ é o grau do vértice $i$ da rede. E como,
as redes tipo \textit{scale-free} tem uma grande quantidade de heterogeneidade, ao contrário das
redes aleatórias, é de se esperar que essa última seja mais resistente à perturbações
lineares. Portanto, pelo Teorema \ref{teo:persistenciaglobal}, podemos comparar a robustez
de diferente tipos de redes olhando apenas para as suas topologias.
De uma forma mais geral, note que $\eta = \eta (\lambda_2)$ e que, considerando redes 
suficientemente grandes, $\lambda_2 \approx g_1$ (Teorema \ref{teo:lambda2cotas}, 
onde $g_1$ é o menor grau da rede, tornando
o lado direito da desigualdade \eqref{eq:criterioperturbacao1} proprorcional à 
$g_1/g_n$ e portanto, para redes
scale-free muito grandes qualquer pequena perturbação poderia destruir a estabilidade
da variedade de sincronização.

Queremos enfatizar que a condição \eqref{eq:criterioperturbacao1} dada pelo Teorema 
\ref{teo:persistenciaglobal} é apenas uma condição necessária, isto é, a mesma não garante 
que a estabilidade da variedade de sincronização 
não seja resistente a perturbações com magnitudes fora da cota estabelecida. Uma pergunta
interessante, a qual não iremos tratar nesta dissertação, é: 
Qual conjunto de perturbações o qual o resultado apresentado é justo? Ou seja,
quais seriam as perturbações tais que, tendo valores de magnitude  fora do estabelicido em
\eqref{eq:criterioperturbacao1}, dostroem a estabilidade da variedade de sincronização?


\chapter{Ilustrações}\label{cap:doisosciladoresacoplados}

Este capítulo aborda o caso de rede $n=2$ que é o exemplo mais trivial de rede.
Apesar de ser o exemplo mais simplório, será possível
extrair a essência dos Teoremas \ref{teo:principal1} e \ref{teo:persistenciaglobal}.

O objetivo então é transformar o problema de estudar a sincronização
entre osciladores difusivamente acoplados num problema de estudar a estabilidade da solução
trivial de uma equação diferencial linear não-autônoma. 
%
Sobre cada oscilador, introduzimos uma cópia da equação \eqref{x=ftx}.
Dessa maneira, todas as variáveis do sistema podem ser escritas da  seguinte forma:
\begin{equation}\label{eq:modelo2osciladores}
\begin{array}{c}
\dot{x}_1=f(x_1)+\alpha[H(x_2)-H(x_1)]\\
\dot{x}_2=f(x_2)+\alpha[H(x_1)-H(x_2)]
\end{array} 
\end{equation}
onde $\alpha$ é o parâmetro global de acoplamento.

\section{Equação Variacional}

Diferentemente da abordagem feita na prova do Teorema \ref{teo:principal1},
utilizaremos o vetor diferença $z(t) = x_1(t) - x_2(t)$ para estudar a estabilidade
da variedade de sincronização, porém,
o resultado é o mesmo descrito pelo Teorema em questão para $n=2$.
Utilizando o teorema da existência global das soluções (\ref{teo:principal1}) garantimos que
as soluções $x_1(t)$ e $x_2(t)$ de \eqref{eq:modelo2osciladores} existem sempre e são 
limitadas valendo-se das hipóteses do referido teorema.
Sabendo que as soluções existem globalmente, o interesse foca-se na determinação do valor do 
parâmetro de acoplamento $\alpha$, tal que
a diferença de estados, ou seja, ${x}_1(t)-{x}_2(t)$, converge ao vetor nulo.
Iremos omitir a dependência temporal das variável por uma questão de conforto na
notação.
Derivando, obtemos
\begin{align*}
\dot{{z}} & = \dot{{x}}_1-\dot{{x}}_2\\
& = {f}({x}_1)+\alpha[H({x}_2)-H({x}_1)] - {f}({x}_2)-\alpha[H({x}_1)-H({x}_2)]\\
& = f(x_1) - f(x_2) - 2\alpha H(z).
\end{align*}
Como por hipótese, devemos ter que $\Vert {z}(0)\Vert \ll 1$, podemos realizar a linearização
de  ${f}({x}_1) ={f}({x}_2+{z})$, baseando-se na expansão em série de Taylor 
ficando com
\begin{align*}
\dot{{z}} & =  f(x_2+z) -f(x_2) - 2\alpha H(z)\\
& = Df(x_2)z   - 2\alpha H(z)+ r(z)
\end{align*}
onde $D{f}({x}_2)$ é a matriz jacobiana de ${f}$ em ${x}_2$.
Assim, 
como por hipótese $f$ é de classe $C^r$, $r\geq2$, utlizamos o Princípio da Linearização 
(Proposição \ref{pro:principio2}),
pois é possível garantir, pelo Teorema da Taylor com Resto de Lagrange \cite{apostol1962calculus}
que $r(z)= \mathcal{O}(\Vert z\Vert^2)$, de forma que o resto de Taylor
não destrói a estabilidade da solução trivial da equação variacional

\begin{equation}\label{eq:equacaovariacional}
 \dot{{z}} = [D{f}({x}_2) - 2\alpha H] {z}.
\end{equation}

\section{Parâmetros Críticos de Acoplamento}\label{section:alphacritico}

Vamos utilizar a Equação \eqref{eq:equacaovariacional}
 para estudar a estabilidade da 
solução trivial ${z}(t)=\mathbf{0}$, à qual representa a variedade de sincronização.
Consideraremos a partir deste momento, e por todo o restante deste capítulo, que $H=I_m$.

Queremos estabelecer condições sobre $\alpha$ para que o operador de evolução
da Equação \eqref{eq:equacaovariacional} tenha contração uniforme e consequentemente
a sua solução trivial seja uniformemente assintoticamente estável.
Podemos assim, utilizar um dos critérios para as contrações uniformes abordados. De maneira mais
específica, vamos utilizar o Teorema da diagonal dominante (\ref{teo:criterio1}).
Portanto, devemos ter que
$$
{D{f}({x}_2)}_{ii} - 2\alpha + 
\sum_{j=1,j\neq i}^m{\left\vert{D{f}({x}_2)}_{ij}\right\vert} <0 
$$
ou
$$
\alpha>\frac{1}{2}\left[ {D{f}({x}_2)}_{ii} +
\sum_{j=1,j\neq i}^m{\left\vert{D{f}({x}_2)}_{ij}\right\vert}\right]<\frac{1}{2}
\left[\sum_{j=1}^m{\left\vert{D{f}({x}_2)}_{ij}\right\vert}\right]
$$
ou ainda, de forma superestimada, e considerando-se
que as condições iniciais estão em uma vizinhança de um domínio absorvente $\Omega$ para
a dinâmica $f$, tomamos
$$
\alpha > \frac{1}{2}\sup_{ \tiny{\begin{matrix}
{x}\in\Omega \\ 
1\leq i\leq m
\end{matrix}}}\left[{\sum_{j=1}^m{\left\vert {D{f}({x})}_{ij}\right\vert}}\right]
$$
e portanto, escrevemos 
\begin{equation}\label{eq:alphacritico2}
\alpha_c = \beta/2
\end{equation}
 onde
\begin{equation}\label{eq:beta}
\beta = \sup_{{x}\in\Omega}{\Vert D{f}({x})\Vert_{\infty}},
\end{equation}
de forma que o conjunto $\Omega$ é advindo da existência da função de Lyapunov para
a dinâmica isolada dos osciladores. Como $\Omega$ é um conjunto compacto, pelo teorema de
Weierstrass (Teorema \ref{weierstrass}) segue que $\beta$ sempre existe.
Note que devido as superestimações sobre o valor de $\alpha$, devemos esperar que
a cota para o $\alpha_c$ que conseguimos seja bem maior do que a realmente necessária
para garantir a sincronização entre os osciladores.

Dessa forma, podemos enunciar o resultado, o qual é um caso particular
do Teorema \ref{teo:principal1}:

\begin{teorema}\label{teo:alphac2}
Considere o  modelo de rede \eqref{eq:modelo2osciladores} com $H=I_m$.
Assuma que a dinâmica isolada tem uma função de Lyapunov com um domínio
absorvente $\Omega$. Além disso, assuma que 
as condições iniciais estão em uma vizinhança de $\Omega$ e considere
$$
\beta = \sup_{{x}\in\Omega}{\Vert D{f}({x})\Vert_{\infty}}.
$$
Seja $\alpha_c = \beta / 2$. Então se $\alpha>\alpha_c$ a variedade de
sincronização é uniformemente assintoticamente estável.
\end{teorema}





Foram feitas simulações computacionais utilizando a dinâmica de Lorenz
como a dinâmica isolada dos osciladores. Sabemos que a dinâmica de Lorenz
é uma dinâmica caótica, então considerando o sistema desacoplado
e com condições iniciais distintas, o mesmo não irá sincronizar. Apenas com base
na intensidade do parâmetro global de acoplamento $\alpha$, poderemos esperar 
que a variedade de sincronização seja estável.
Além disso, como tratado na Subseção \ref{subsection:flsl}, o sistema de Lorenz admite
uma função de Lyapunov. Dessa forma, tal dinâmica cumpre as hipóteses do Teorema \ref{teo:alphac2},
considerando as condições iniciais sempre numa vizinhança de $\Omega$,
neste caso como definido em \eqref{eq:lorenzatractor}.

Considere o Teorema \ref{teo:alphac2}. 
De posse das principais ferramentas, vamos encontrar o $\alpha_c$ para os dois
osciladores acoplados.
Sendo $\mathbf{x}=(x,y,z)$, a matriz jacobiana do sistema de Lorenz é
$$
Df(\mathbf{x}) =
\begin{pmatrix}
-\sigma & \sigma & 0 \\ 
r-z & -1 & -x \\ 
y & x & -b
\end{pmatrix}
$$
e sua norma infinito é dada por 
$\Vert Df(\mathbf{x})\Vert_\infty = \max \{ 2\sigma, \vert r-z\vert +1 +\vert x\vert, \vert y \vert +\vert x\vert +b\}$,
de forma que queremos calcular $\beta = \sup_{\mathbf{x}\in\Omega}\Vert Df(\mathbf{x})\Vert_\infty$.

As trajetórias encontram-se dentro do domínio absorvente $\Omega$ \eqref{eq:lorenzatractor} 
que por sua vez é uma região limitada por um elipsoide,
então o valor máximo atingido por $x$ acontece quando $y=0$ e $z=2r$, da mesma forma $y$ é
máximo quando $x=0$ e $z=2r$ e $z$ é máximo quando $x=y=0$. Assim, podemos tomar as seguintes cotas:
$$
\vert x\vert\leq\frac{b\sqrt{r}}{\sqrt{b-1}},\quad\vert y\vert\leq\frac{rb}{\sqrt{\sigma(b-1)}}\quad 
\mbox{e}\quad\vert z-2r\vert\leq \frac{rb}{\sqrt{\sigma(b-1)}}.
$$
Mas com relação a $z$ queremos estimar $\vert r-z\vert=\vert z-r\vert$. Usamos então a desigualdade
triangular inversa, assim
$$
\frac{rb}{\sqrt{\sigma(b-1)}}\geq \vert z-2r\vert = \vert (z-r)-r\vert\geq 
\vert \vert z-r\vert -r\vert
$$
então
$$
\vert \vert z-r\vert -r\vert \leq \frac{rb}{\sqrt{\sigma(b-1)}}
$$
e portanto
$$
r- \frac{rb}{\sqrt{\sigma(b-1)}} \leq \vert z-r \vert \leq
\frac{rb}{\sqrt{\sigma(b-1)}} +r.
$$
Assim,
$$\beta= \max\left\{ 2\sigma,\frac{rb}{\sqrt{\sigma(b-1)}} +r+1+ \frac{b\sqrt{r}}{\sqrt{b-1}},
\frac{b\sqrt{r}}{\sqrt{b-1}} +\frac{rb}{\sqrt{\sigma(b-1)}} +b \right\}.$$

\begin{figure}[!ht]
\begin{center}
	\includegraphics[height=6cm]{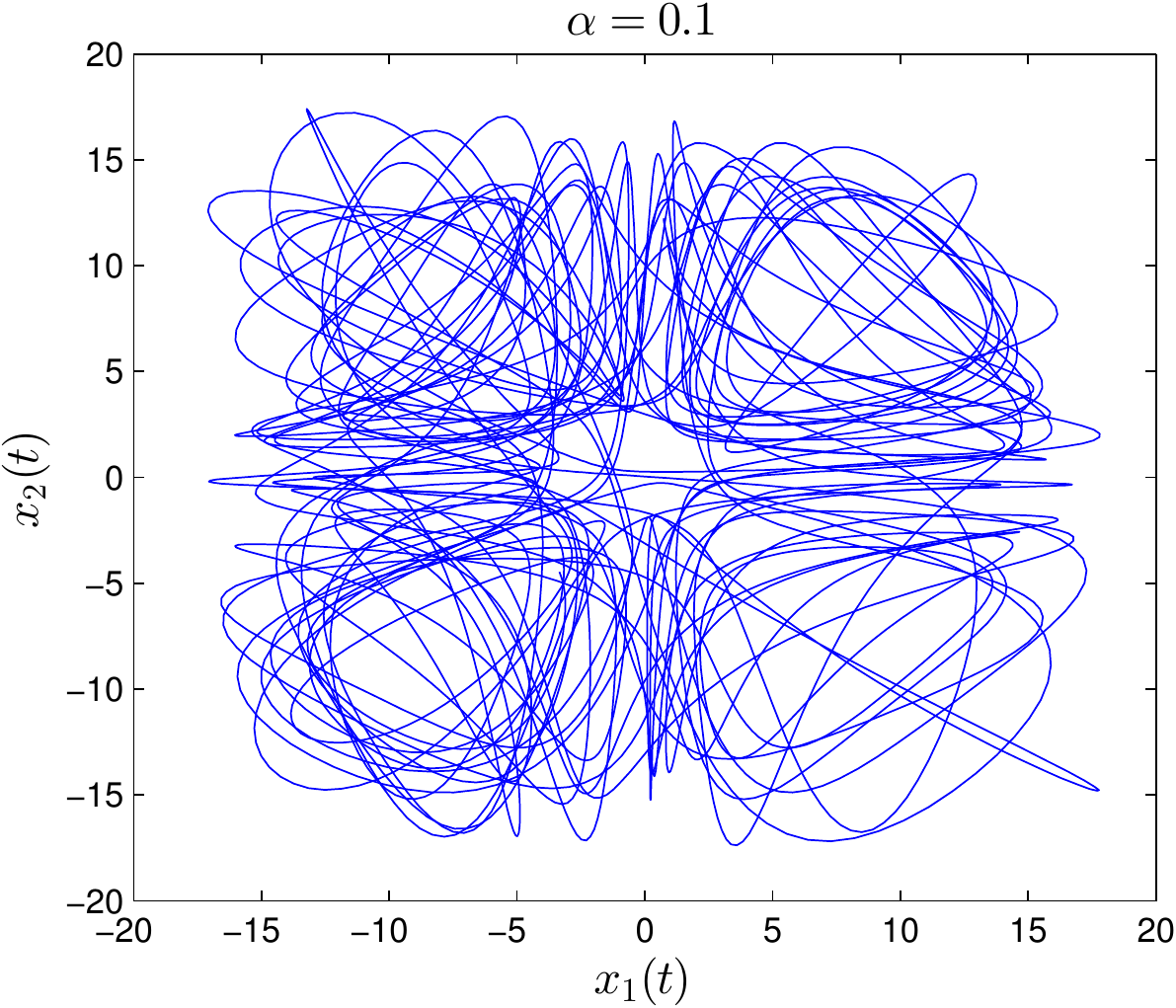}
	\includegraphics[height=6cm]{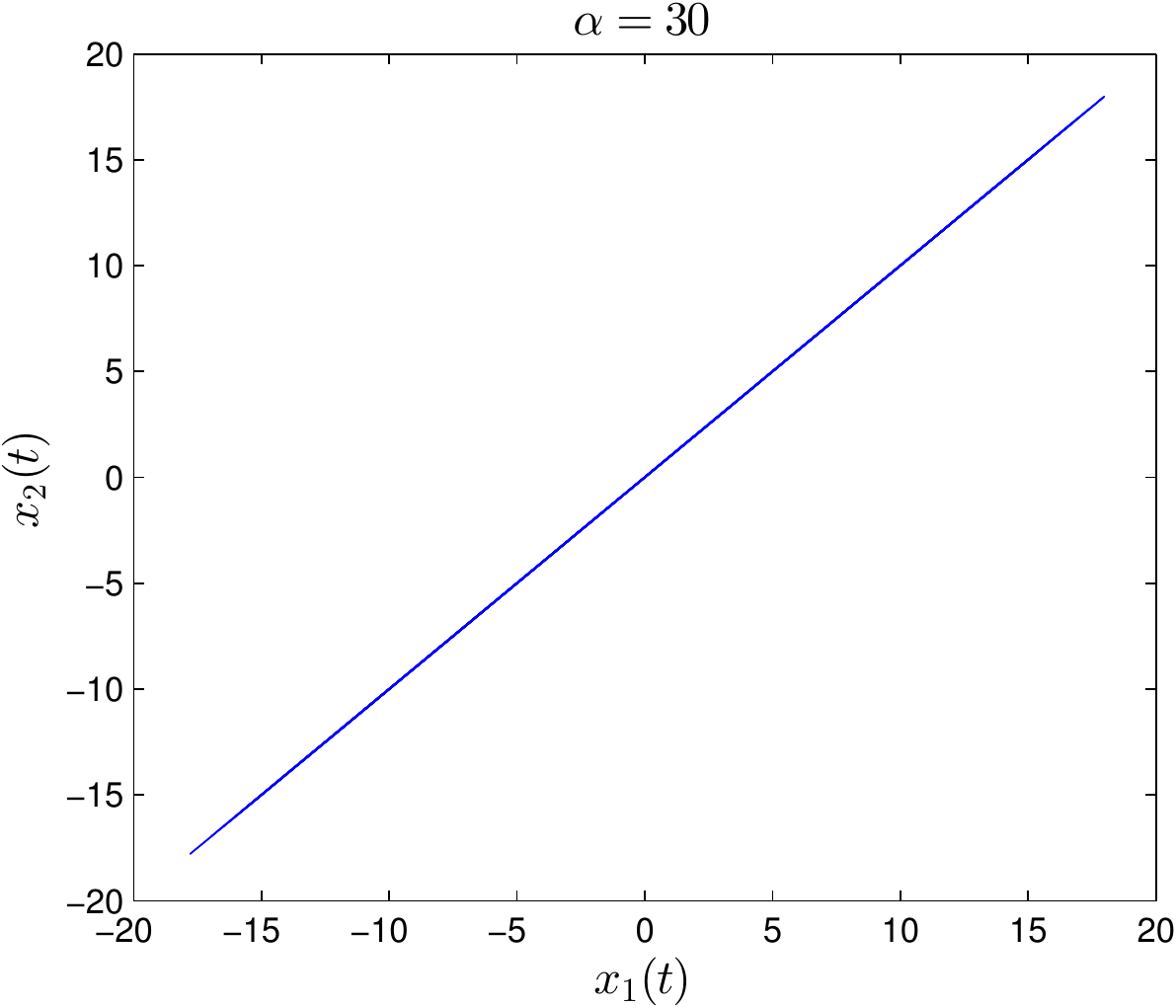}
\caption{Gráficos de $x_1(t)\times x_2(t)$ com $0\leq t\leq 100$ onde $x_1(t)$ e $x_2(t)$ são 
respectivamente
as primeiras componentes das trajetórias de $\mathbf{x}_1(t) = (x_1(t),y_1(t),z_1(t))$ e
$\mathbf{x}_2(t) = (x_2(t),y_2(t),z_2(t))$ 
no modelo de dois osciladores acoplados com dinâmica de Lorenz. Do
lado esquerdo com parâmetro de acoplamento igual a $0.1$ e do lado direito com $30$.} 
\label{fig:x11x21a29a01}
\end{center}
\end{figure}

Então, substituindo os valores clássicos dos parâmetros do modelo de Lorenz,
ou seja, $\sigma=10$, $r=28$ e $b=8/3$ obtemos
$$\beta=\frac{rb}{\sqrt{\sigma(b-1)}} +r+1+ \frac{b\sqrt{r}}{\sqrt{b-1}} \approx 58.22.$$
Portanto, o parâmetro crítico de acoplamento é
$$
\alpha_c = \frac{\beta}{2} \approx 29.11.
$$
A Figura \ref{fig:x11x21a29a01} mostra a exitência da variedade de sincronização 
considerando $\alpha>\alpha_c$.

Novamente, devido as superestimações para se obter o $\alpha_c$,
é de se esperar que o real parâmetro crítico para se observar sincronização entre os
osciladores seja bem menor do que $29.11$ no caso do sistema de Lorenz
considerando os parâmetros usuais. A Figura \ref{fig:alphanorm2} traz um resultado da simulação
de $\Vert \mathbf{x}_1(t) - \mathbf{x}_2(t)\Vert_2$ em função de $\alpha$. Valores de $\alpha>0.5$ produzem
$\lim_{t\rightarrow\infty}\Vert \mathbf{x}_1(t) - \mathbf{x}_2(t)\Vert_2 =0$. 

\subsection{Utilizando o Critério da Parte Simétrica}
Considere o sistema \eqref{eq:equacaovariacional}. Seja $M(\mathbf{x}) = Df(\mathbf{x})-2\alpha I_3$. 
Fazendo uso do Teorema da Parte Simétrica (\ref{teo:criterio2}), queremos estabelecer condições 
sobre $\alpha$ tal que a parte simétrica
de $M(\mathbf{x})$ seja negativa-definida com $\mathbf{x}=(x,y,z)\in\Omega$, para tanto, 
consideraremos também o Corolário \ref{cor:matriznegativadefinida}.
A parte simétrica de $M(\mathbf{x})$, aqui denotada por $P_s(M(\mathbf{x}))$, se lê
\begin{equation}
P_s(M(\mathbf{x})) = \frac{Df(\mathbf{x})+{Df(\mathbf{x})}^\dagger}{2}-2\alpha I_3.
\end{equation}
De forma explícita, temos
$$
P_s(M(\mathbf{x})) =
\begin{pmatrix}
-\sigma -2\alpha & ({r+\sigma - z})/{2} & {y}/{2}\\
({r+\sigma - z})/{2} & -1-2\alpha & 0\\
{y}/{2} & 0 & -b-2\alpha
\end{pmatrix}.
$$

Vamos requer que os autovalores de $P_s(M(\mathbf{x}))$ sejam todos negativos.
Note que $P_s(M(\mathbf{x}))$ depende ainda de $y$ e $z$, o que significa
que o $\alpha_c$ é o maior $\alpha$ tal que os autovalores de $P_s(M(\mathbf{x}))$ sejam
todos negativos com $\mathbf{x}\in\Omega$. Mas, de
posse das cotas $\vert y\vert\leq {rb}/\sqrt{\sigma(b-1)}$ e $\vert z-2r\vert\leq 
{br}/\sqrt{\sigma(b-1)}$, para cada $-{rb}/\sqrt{\sigma(b-1)}\leq y\leq 
{rb}/\sqrt{\sigma(b-1)}$ fixado, variamos $2r-{br}/\sqrt{\sigma(b-1)}\leq z\leq 
2r+{br}/\sqrt{\sigma(b-1)}$ e a medida que fazemos essa variação, verificamos quais valores 
de $\alpha$ produzem $\lambda_1 \leq \lambda_2\leq\lambda_3 <0$.

Foram feitas simulações computacionais para o procedimento indicado acima resultando, com precisão
de uma casa decimal, tem-se seguinte parâmetro crítico de acoplamento:
\begin{equation}\label{eq:criterio2}
\alpha_c = 13.03.
\end{equation}

Agora, utilizando o item $3$ do Corolário \ref{cor:matriznegativadefinida} vamos impor
que os determinantes dos menores principais de $P_s(M(\mathbf{x}))$ alternam os sinais começando
com $\mbox{det}(P_s(M(\mathbf{x}))_1)<0$, onde indiciamos por $P_h(M(\mathbf{x}))_1$, 
$P_s(M(\mathbf{x}))_2$ e $P_s(M(\mathbf{x}))_3$
os respectivos menores principais de ordem $1$, $2$ e $3$. A ideia é que para cada determinante
devemos tomar o maior $\alpha$ do polinômio resultante tal que esse polinômio produz valores positivos
ou valores negativos dependendo da ordem do menor principal. E em seguida, tomamos o maior
valor de $\alpha$ dentre os obtidos na etapa anterior, então esse será o nosso $\alpha_c$
para este caso.

Fazendo os cálculos temos
$$
\mbox{det}(P_s(M(\mathbf{x}))_1) = -\sigma-2\alpha <0
$$
o que implica $\alpha>-\sigma/2$, o que é irrelevante para a análise pois estamos considerando
apenas valores positivos. Para o segundo determinante temos
$$
\mbox{det}(P_s(M(\mathbf{x}))_2) = (-\sigma-2\alpha)(1-2\alpha)-\frac{(r+\sigma-z)^2}{4}>0.
$$
Para facilitar a notação faça $p_2(\alpha) = \mbox{det}(P_s(M(\mathbf{x}))_2)$. Então
\begin{align*}
p_2(\alpha) & = (-\sigma-2\alpha)(1-2\alpha)-\frac{(r+\sigma-z)^2}{4}\\
& = 4\alpha^2+2(\sigma+1)\alpha+\sigma-1/4(r+s-z)^2
\end{align*}
Note que $p_2(\alpha)$
depende ainda de $z$. Então o procedimento é o seguinte: 
verificamos quais as raízes de $p_2(\alpha)$ com $\mathbf{x}\in\Omega$ (\eqref{eq:lorenzatractor})
e tomamos a maior dentre todas as raízes positivas. 
Isso porque o coeficiente de 
$\alpha^2$ em $p_2(\alpha)$ é positivo indicando que $p_2(\alpha)>0$ quando $\alpha$ assume valores 
entre $(-\infty,\alpha_i)\cup (\alpha_j,+ \infty)$, considerando que $\alpha_i$ e $\alpha_j$ são as
raízes de $p_2(\alpha)$ e $\alpha_i<\alpha_j$.
Note ainda que $p_2(\alpha)$ não possui raízes
complexas.
%
Simulações computacionais sugerem que $\alpha>6.5972$ produz $p_2(\alpha)>0$ com
$\mathbf{x}\in\Omega$. De forma mais específica, notamos que a maior raiz de $p_2(\alpha)$ é dada por
$$
\alpha_j = \frac{-2(\sigma+1)+2\sqrt{(\sigma+1)^2-4\sigma+(r+\sigma-z)^2}}{8},
$$
onde $\alpha_j$ atinge o seu valor máximo quando $(r+\sigma-z)^2$ é máximo, implicado que 
esse valor máximo é atingido em $z=2r+rb/\sqrt{\sigma(b-1)}$. Substituindo os valores
de $\sigma$, $b$ e $r$ temos
$$
\alpha_j \approx 6.5972.
$$

Por fim, calculando os valores de $\alpha$ que produzem $p_3(\alpha) = 
\mbox{det}(P_h(M(\mathbf{x}))_3)<0$ temos
\begin{equation}
p_3(\alpha) = -8\alpha^3+[-4b-4(\sigma+1)]\alpha^2+\left[-2\sigma-2b(\sigma+1)+\frac{(r+\sigma-z)^2}{2}
+\frac{y^2}{2}\right]\alpha - \sigma b+\frac{b}{4}(r+\sigma-z)^2+\frac{y^2}{4}.
\end{equation}
É possível verficar que para 
todo $y$ e $z$ as raízes de $p_3(\alpha)$ não assumem valores complexos. Descartamos
as raízes negativas e dentre as positivas verificamos quais produzem $p_3(\alpha)<0$ para
cada $y$ e $z$, e dentre essas tomamos a maior. As simulações computacionais produzem
o valor $\alpha>7.5546$, e além disso, que tal valor é atingido quando $z=2r+rb/\sqrt{\sigma(b-1)}$ e 
$y=-{rb}/{\sqrt{\sigma(b-1)}}$.
Assim, podemos tomar então
\begin{equation}
\alpha_c= 7.5546
\end{equation}
que é o maior $\alpha$ que produz todos as desigualdades relacionadas aos menores principais.

Portanto, temos dois novos parâmetros críticos de acoplamento que são um pouco melhores do que
aquele apresentado na seção anterior. Porém, a vantagem de utilizar o Teorema da diagonal
dominante para encontrar o parâmetro crítico de acoplamento é que a mesma análise pode ser estendida 
para redes em geral como pode-se ver na prova do Teorema \ref{teo:principal1} exibida
no Capítulo \ref{cap:estabilidadedassolucoessincronas}.

\begin{figure}[!ht]
\begin{center}
	\includegraphics[height=7cm]{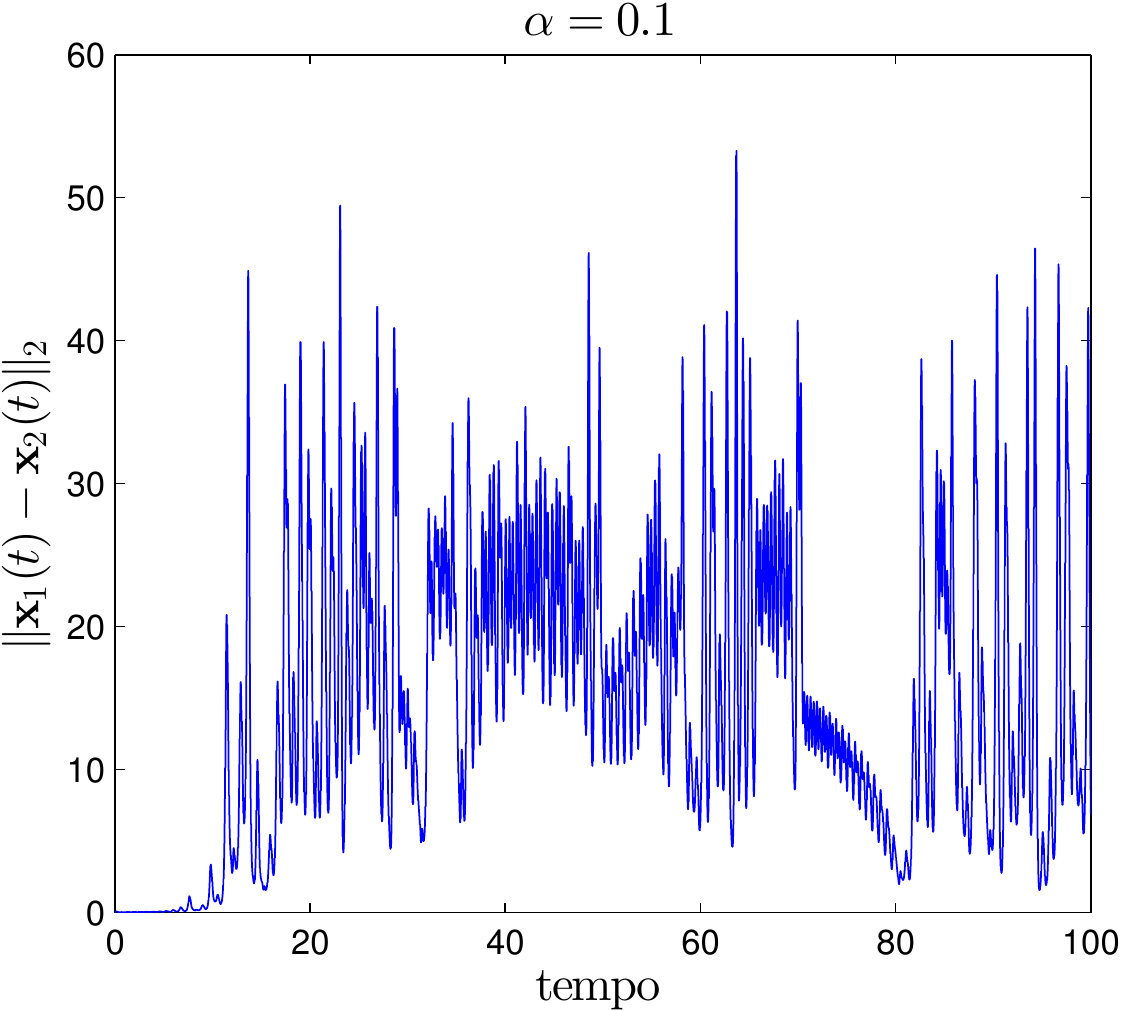}
	\includegraphics[height=7cm]{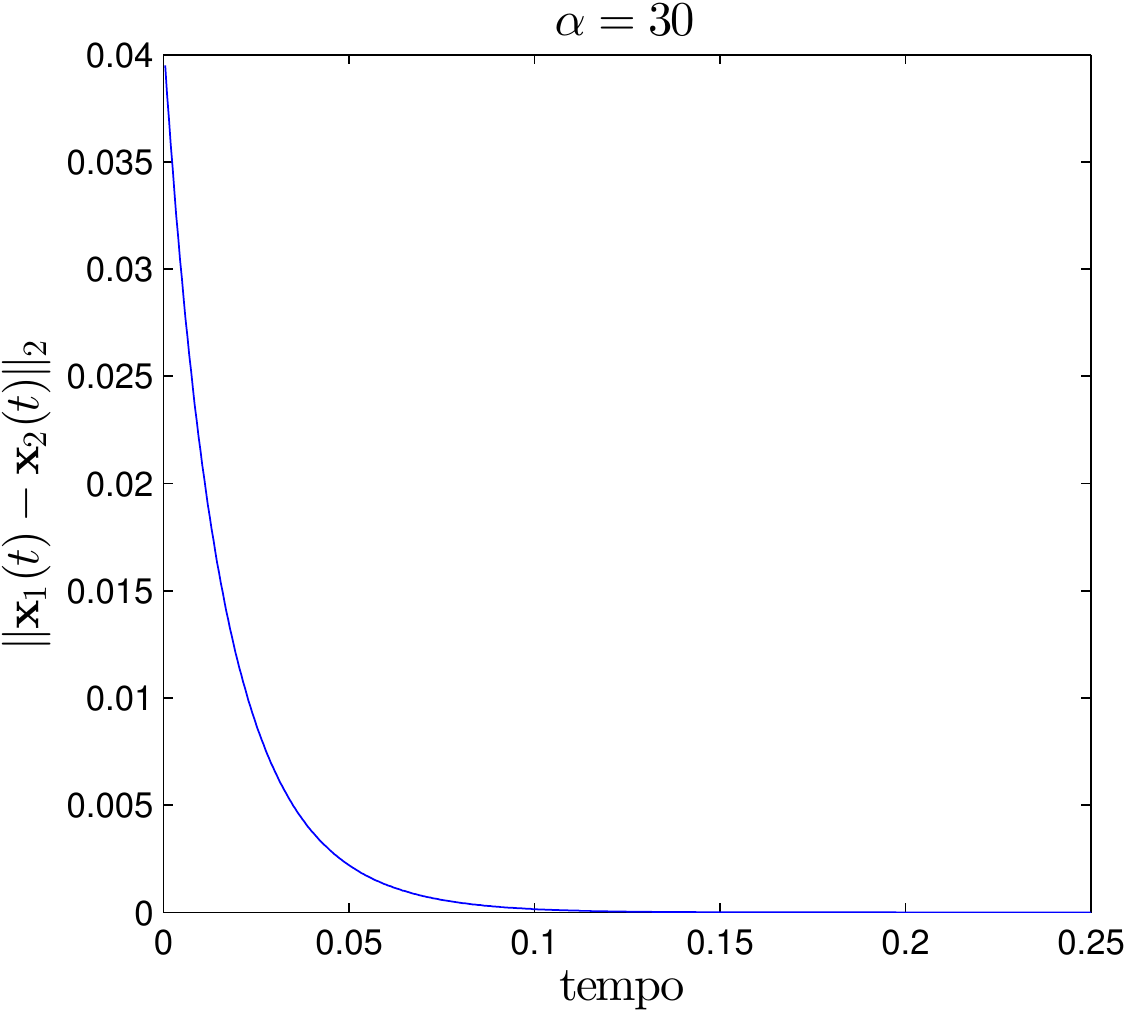}
\caption{Do lado esquerdo temos a norma euclidiana da diferença de estados das trajetórias dos osciladores 
em função do tempo para $\alpha = 0.1$ e do lado direito para $\alpha = 30$.} \label{fig:alpha01ealpha29}
\end{center}
\end{figure}

\begin{figure}[!ht]
\begin{center}
        \includegraphics[height=7cm]{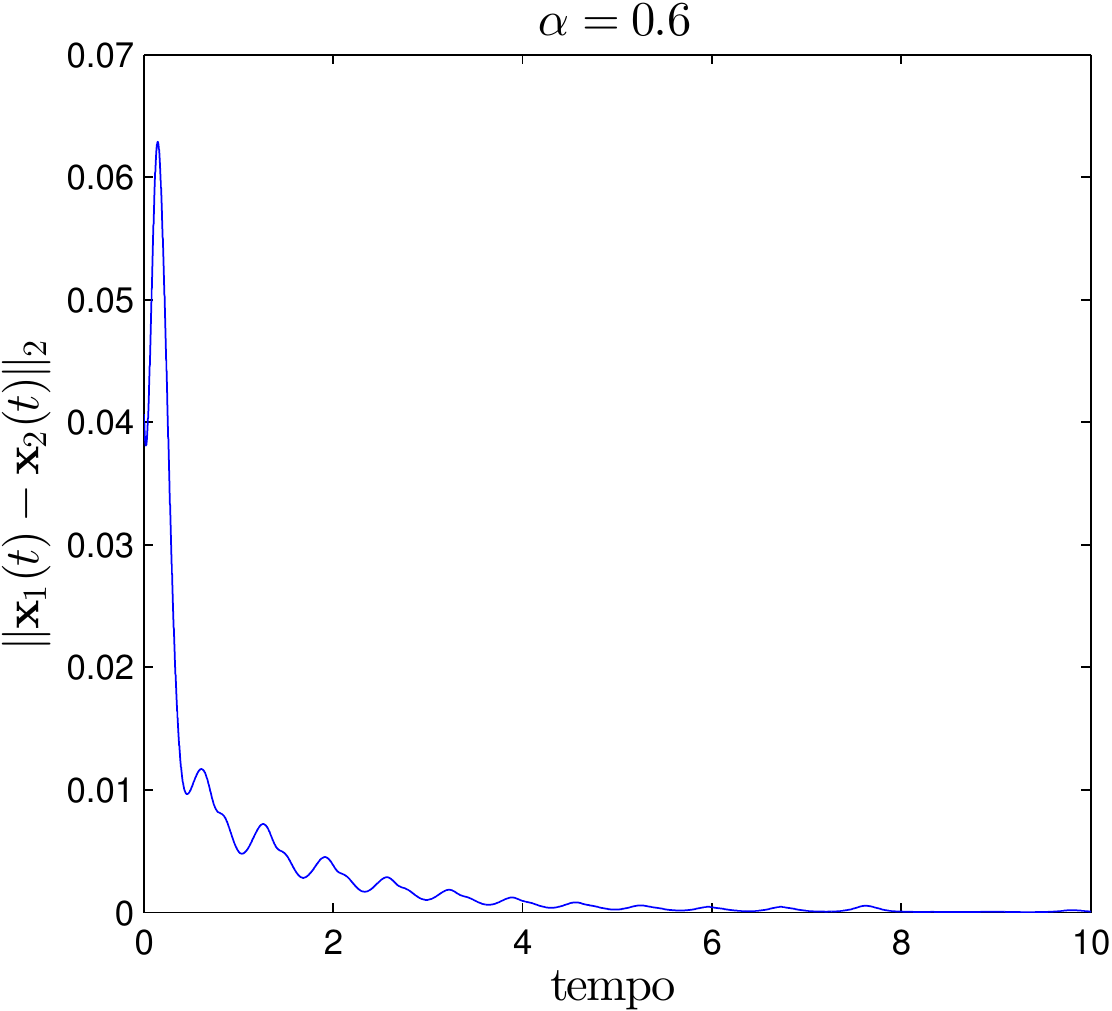}
	\includegraphics[height=6.7cm]{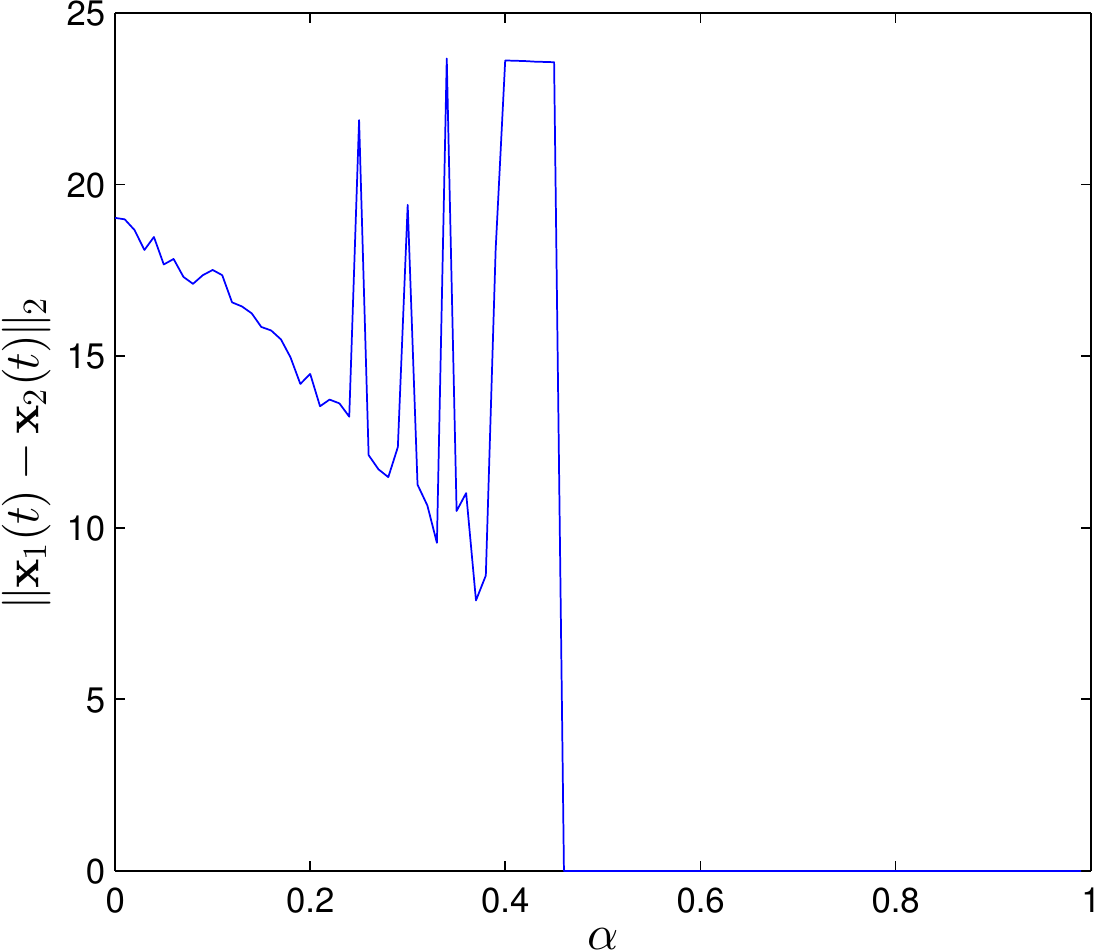}
\caption{Do lado esquerdo temos agora $\alpha=0.6$. Verificamos que
o decaimento é mais lento porém o sistema tende à sincronização.
Do lado direito temos a relação entre a norma euclidiana da diferença
dos estados contra o parâmetro global de acoplamento. Neste caso,
estamos tomando a média temporal da norma euclidiana com $1000\leq t\leq 2000$ para cada $\alpha$. 
Notoriamente, para valores $\alpha >0.5$ o sistema irá sincronizar.} \label{fig:alphanorm2}
\end{center}
\end{figure}

Podemos observar alguns resultados das simulações, os quais são apresentados nas figuras
a seguir.
Para $\alpha = 0.1$ verificamos que não há sincronização entre os osciladores como
mostra a figura \ref{fig:alpha01ealpha29}. 
Em contraste, para $\alpha = 30$ verificamos o decaimento exponencial em $t$.
De forma um pouco mais precisa temos a relação entre o parâmetro de acoplamento $\alpha$
e a norma euclidiana da diferença dos estados $\mathbf{x}_1(t)$ e $\mathbf{x}_2(t)$ como mostra a figura
\ref{fig:alphanorm2}

\newpage
\section{Perturbações em Dois Osciladores Acoplados}

Da mesma forma que sabemos que o parâmetro crítico de acoplamento pode ser muitas vezes
maior do que o realmente necessário para se obter sincronização persistente,
a magnitude das perturbações poderia exceder as cotas estabelecidas 
pelo Teorema \ref{teo:persistenciaglobal} e ainda assim a equação perturbada continuar 
na vizinhança de atração da variedade de sincronização. 

Nesta seção vamos trabalhar com alguns exemplos para o Teorema \ref{teo:persistenciaglobal}
considerando $n=2$, $H=I_m$ e o sistema isolado com a  dinâmica de Lorenz. Assim, como
$\Vert L\Vert_\infty =\Vert L\Vert_2 = 2$ e $m=3$, então o valor estabelecido em
\eqref{eq:criterioperturbacao1} se resume a
$$
\sup_{t}\sum_{i,j=1;i\sim j}^2 \Vert V_{ij}\Vert_\infty < \frac{\eta}{4}
$$

\begin{exemplo}\label{exemplo:V}
Considere o Teorema \ref{teo:persistenciaglobal} para $n=2$, $H=I_3$, a dinâmica isolada sendo
o sistema de Lorenz e as matrizes de perturbação
$V_{11} = V_{22}= V_{21} = \mathbf{0}\in\Mat(\mathbb{R},3)$ e 
$V_{12}$ a matriz constante
$$
V_{12}=
\begin{pmatrix}
\xi & 0 & -\xi\\
0 & -\xi & 0\\
-\xi & 0 & \xi
\end{pmatrix} = \xi
\begin{pmatrix}
1 & 0 & -1\\
0 & -1 & 0\\
-1 & 0 & 1
\end{pmatrix}.
$$
Vamos tomar $\alpha = 30$.
Já sabemos que, para o caso que estamos considerando, temos $\beta = 58.22$, de forma que tomamos 
$\alpha=30>\alpha_c = \beta/2$, produzindo assim $\eta = 30 - 29.11 = 0.89$. 
Então tem-se
$$
\sup_{t}\sum_{i,j=1;i\sim j}^2 \Vert V_{ij}\Vert_\infty  =\Vert V_{12}\Vert_\infty = 2\vert \xi \vert,
$$
que deve satisfazer
$2\vert \xi \vert < {\eta}/{4}  = 0.2225,$
implicando que $-0.11125 < \xi<0.11125$. Porém, simulações computacionais mostram que
para esta faixa de intervalo, a variedade de sincronização existe para $\alpha>0.6$  como
pode ser observado na Figura \ref{fig:a0_10xi-10_20}, portanto e certamente
para $\alpha=30$.
%
%
Tal figura tem como objetivo verificar regiões de sincronização de forma que 
podemos fazer essa análise através de um
mapa de cores. Por padrão, a região em preto, no mapa, é a região
que verificamos a propriedades citada. As cores no mapa representam
a magnitude da média temporal de $\Vert \mathbf{x}_1(t)-\mathbf{x}_2(t)\Vert_2$
com $t$ entre $1000$ e $2000$ para cada $\alpha$ e $\xi$, onde $\mathbf{x}_1 =(x_1,y_1,z_1)$ e
$\mathbf{x}_2 = (x_2,y_2,z_2)$ são as duas trajetórias para o sistema de Lorenz, de forma 
que essa variação acontece continuamente a partir de $0$ (preto)
até o valor máximo atingido pela norma euclidiana da diferença
de estados (amarelo).

\begin{figure}[!ht]
\begin{center}
	\includegraphics[scale=0.6,angle=-90]{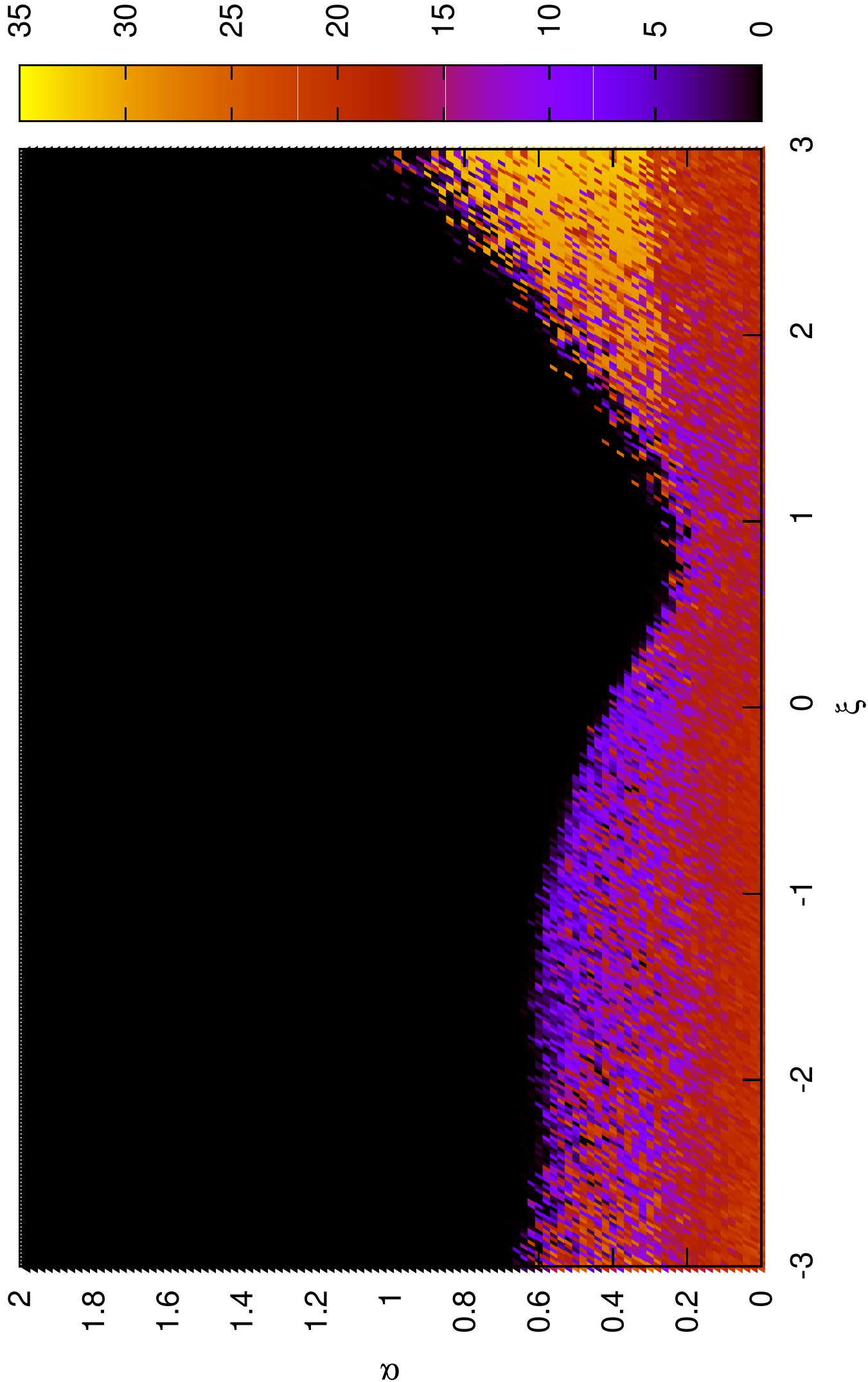}
\caption{Mapa de cores para o Exemplo \ref{exemplo:V}. A escala de cores representa
$\Vert \mathbf{x}_1(t)-\mathbf{x}_2(t)\Vert_2$.} \label{fig:a0_10xi-10_20}
\end{center}
\end{figure}

\end{exemplo}



\begin{lema}\label{lema:pertint}
Sejam $A(t)$ e $B(t)$ funções matriciais contínuas em um intervalo $J\subset\R_+$ tais que
$\Vert A(t)\Vert \leq M$, $\Vert B(t) \Vert \leq M$, e suponha que o operador de evolução
da equação $\dot{x} = A(t)x$ cumpre $\Vert T(t,s) \Vert \leq k e^{\eta(t-s)}$, $\eta<0$ e $k\geq1$. 
Se
\begin{equation*}
\left\Vert \int_{t_1}^{t_2} B(t) dt \right\Vert \leq \delta \quad \mbox{para} \quad 
\vert t_1 - t_2\vert \leq h,
\end{equation*}
então o operador de evolução da equação perturbada $\dot{y} = [A(t)+B(t)]y$ cumpre
\begin{equation*}
\Vert \tilde{T}(t,s)\Vert \leq (1+\delta)ke^{\beta(t-s)}
\end{equation*}
onde $\beta = \eta +3Mk\delta +h^{-1}\ln{(1+\delta)}k$.
\end{lema}

Este resultado, cuja prova pode ser encontrada em \cite{coppel} página 7, mostra que
se os coeficientes das matrizes $A(t)$ e $B(t)$ são limitados, então a estabilidade
uniformemente assintótica é preservada sobre perturbações ``integralmente pequenas''.
Uma aplicação então que pode ser observada é que a estabilidade uniformemente assintótica
não é destruída sobre perturbações que oscilam suficientemente rápido.

\begin{exemplo}\label{exemplo:vt}
Considere agora as mesmas hipóteses feitas no Exemplo \ref{exemplo:V} 
com uma única diferença, ou seja, considere agora a matriz de perturbação não-constante 
$$
V_{12} =
\begin{pmatrix}
\xi \cos{(\omega t)} & \xi \cos{(\omega t)} & \xi \cos{(\omega t)}\\
\xi \cos{(\omega t)} & 0 & 0\\
\xi \cos{(\omega t)} & 0 & 0
\end{pmatrix} =\xi \cos{(\omega t)}
\begin{pmatrix}
1 & 1 & 1\\
1 & 0 & 0\\
1 & 0 & 0
\end{pmatrix},
$$
onde a princípio tomamos $\omega$ sendo a frequência de oscilação do sistema de Lorenz para
os parâmetros usuais $\sigma=10$, $r=29$ e $b=8/3$. 
Verifica-se, através dos dados obtidos pelas simulações 
que a frequência dominante é aproximadamente $\omega = 4.1888$.

Pelo Teorema 
\ref{teo:persistenciaglobal}, para que a magnitude da perturbação não altere a estabilidade
da variedade de sincronização, faz-se necessário que
$$
\sup_{t}\sum_{i,j=1;i\sim j}^2 \Vert V_{ij}\Vert_\infty  = \sup_t \Vert V_{12}\Vert_\infty < \frac{\eta}{4}.
$$
Então, como $\sup_t \Vert V_{12}\Vert_\infty =  \sup_t 3\vert \xi\vert \vert \cos (\omega t) \vert
= 3\vert \xi \vert$, teremos que $-7.42\times 10^{-2} < \xi< 7.42\times 10^{-2}$.
Plotamos, na Figura \ref{fig:colormapw41888}, valores um pouco mais realísticos
para $\xi$ e $\alpha$, produzindo novamente regiões de sincronização (em preto).

É perceptível a simetria que encontramos. Tal
simetria está relacionada com a dependência temporal em função do cosseno.
Além disso, a inclinação da ``língua'' que observamos está intimamente
relacionada com o parâmetro $\omega$.
\begin{figure}[!ht]
\begin{center}
	\includegraphics[scale=0.6,angle=-90]{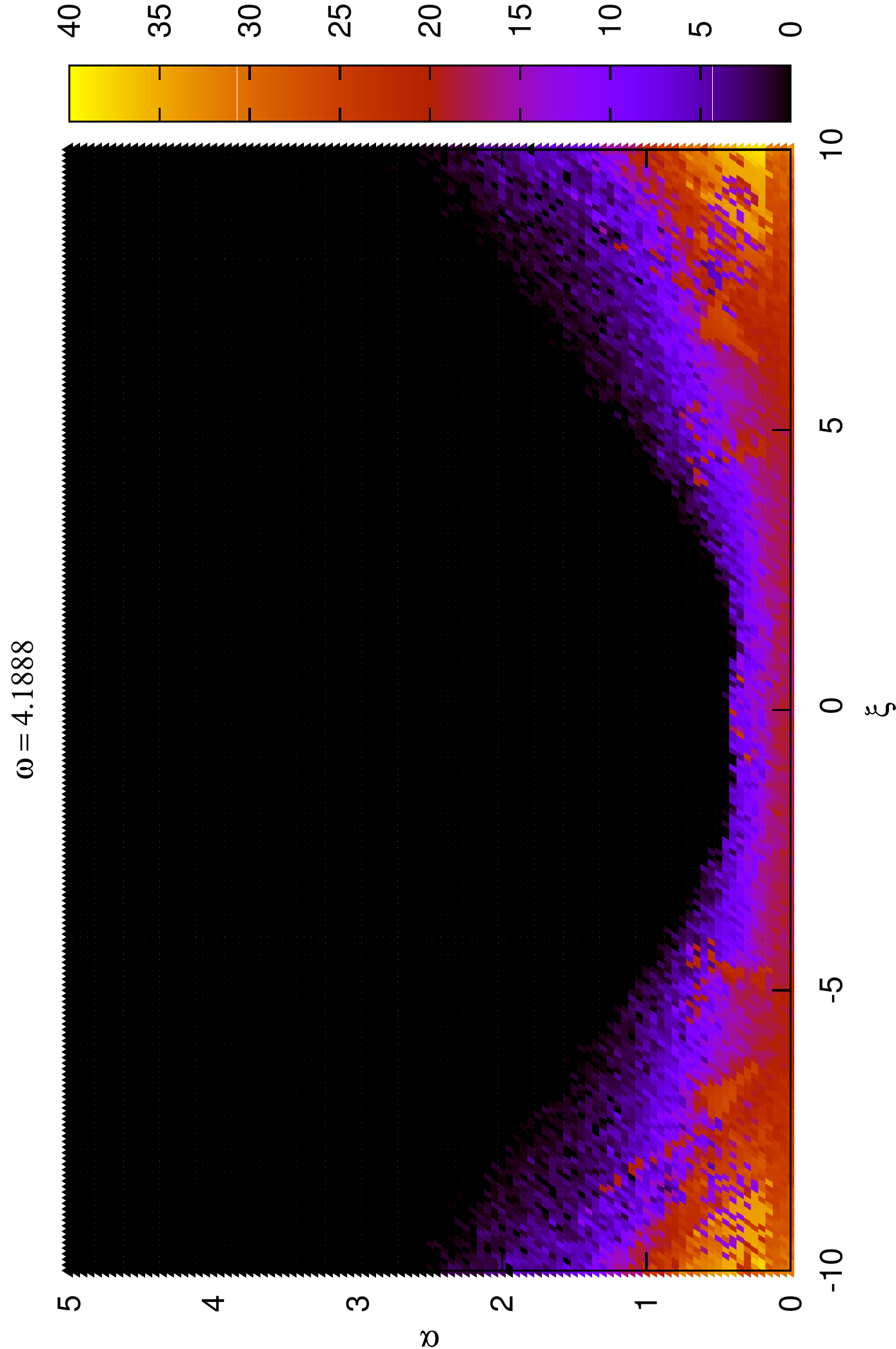}
\caption{Mapa de cores para o Exemplo \ref{exemplo:vt} com $\omega=4.1888$.}\label{fig:colormapw41888}
\end{center}
\end{figure}

Simulações computacionais mostram que os mapas de cores para baixas frequências são bastante 
oscilatórios
e não simétricos, de forma que não é possível verificar inclinação neste caso.
Portanto, se $\omega$ é pequeno, isto é, se $\omega\approx 0$,
então isto equivale a ter uma perturbação $V_{12}$ com a matriz constante 
$$V=
\begin{pmatrix}
1 & 1 & 1\\
1 & 0 & 0\\
1 & 0 & 0
\end{pmatrix},
$$ 
de fato pois $\lim_{\omega\rightarrow 0} \cos(\omega t) =1$.
Iremos considerar então apenas valores ``grandes'' para $\omega$ e
apenas valores não-negativos visto que a análise seria recíproca
pois o cosseno é uma função par.
%

A pergunta interessante é: o que acontece então quando $\omega\rightarrow\infty$ ?
E a resposta é: isso equivale a não ter pertubação, ou seja, tornando $V_{12}=\mathbf{0}$. 
Esta conclusão pode ser verificada pelo Lema \ref{lema:pertint}. Note que
$$
\lim_{\omega\rightarrow\infty}\left\Vert \int V_{12}(t)dt\right\Vert_\infty =
\lim_{\omega\rightarrow\infty} 3 \vert\xi\vert \int \vert \cos(\omega t)\vert dt =
\lim_{\omega\rightarrow\infty} \frac{\sin(\omega t) \operatorname{sgn}(\omega t)}{w} =0
$$
onde ``$\operatorname{sgn}$'' é a função sinal.
Intuitivamente, quando $\omega$ cresce, a perturbação tende a ter uma oscilação muito rápida,
fazendo com que a estabilidade da variedade de sincronização, definida por $x_1(t) = x_2(t)$,
não seja afetada - de acordo com o Lema \ref{lema:pertint}.

Além disso, sabemos que para o modelo de dois osciladores acoplados com dinâmica de Lorenz e sem 
perturbação, o real parâmetro de acoplamento que se observa sincronização é aproximadamente $0.5$.
Então,
a medida que aumenta-se o valor de $\omega$ no modelo de perturbação considerado neste exemplo,  
menor será a inclinação
da ``língua'' no mapa de cores, e portanto quando $\omega \rightarrow \infty$ restará apenas
regiões de não-sincronização abaixo de valores menores que $0.5$ para $\alpha$. Este fato
pode ser observado na Figura \ref{fig:colormaps}.

\begin{figure}[!ht]
\begin{center}
        \includegraphics[scale=0.35,angle=-90]{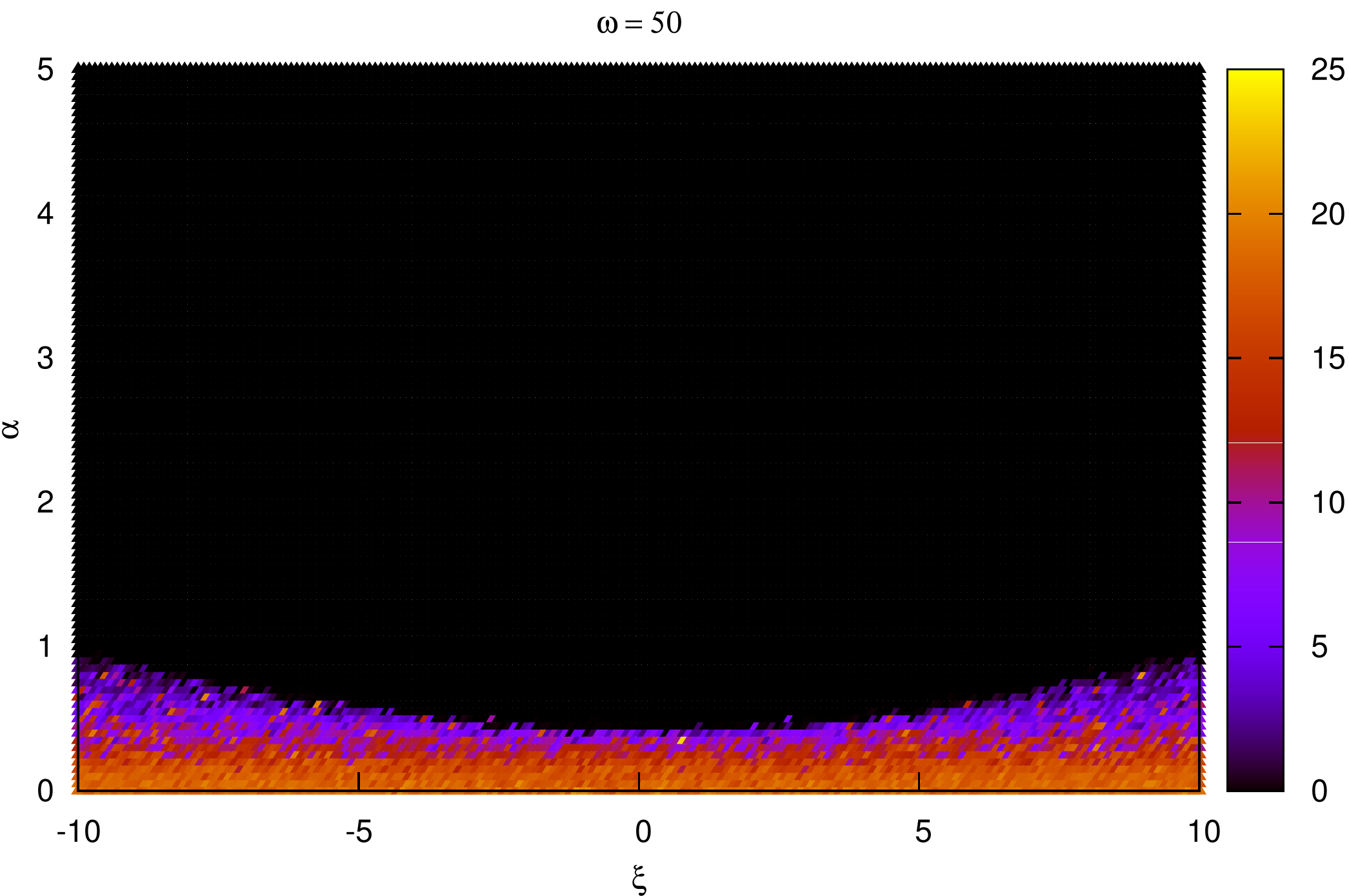}\quad
        \includegraphics[scale=0.35,angle=-90]{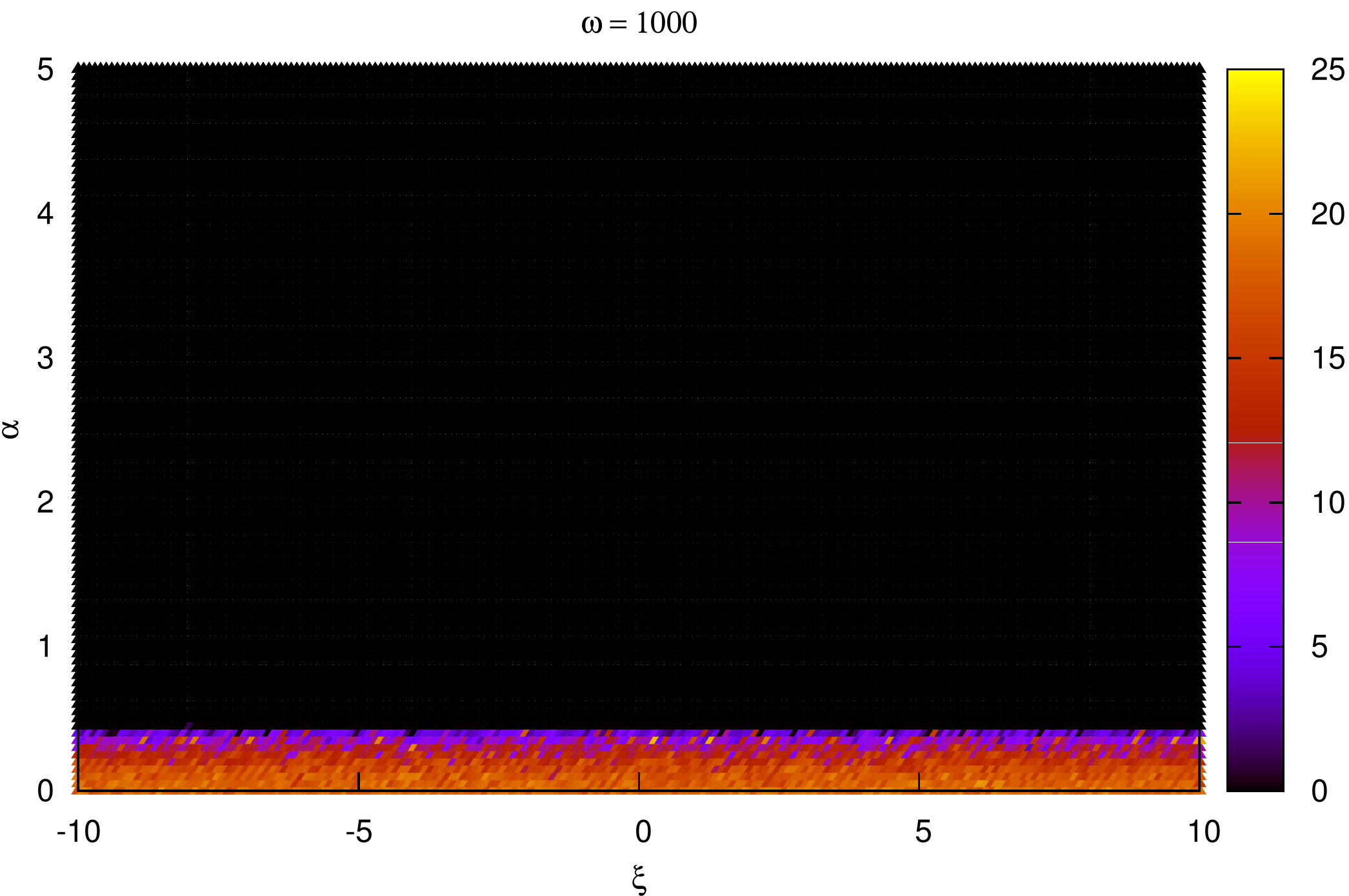}
\caption{Mapas de cores para diferentes valores de $\omega$, relativo ao Exemplo \ref{exemplo:vt}.}
\label{fig:colormaps}
\end{center}
\end{figure}

\end{exemplo}
Dessa forma, possuindo dependência em $t$ ou não, nota-se que a variedade de sincronização 
é ``bastante''
persistente, isso devido a força (magnitude) do parâmetro global de acoplamento $\alpha$, visto
que o $\alpha_c$ é bem maior do que o real $\alpha_c$
para se observar ao menos a sincronização. Então seguindo os mesmos passos
é de se esperar que a região de estabilidade uniformemente assintótica pode ser maior do que
aquela que podemos garantir. 

\chapter{Prova dos Teoremas Fundamentais}
\label{cap:estabilidadedassolucoessincronas}

\section{Existência Global das Soluções}\label{sec:egs}

Iremos construir uma função de Lyapunov para a dinâmica coletiva de redes de 
$n$ osciladores difusivamente acoplados, dada a hipótese de que a dinâmica individual possui
uma função de Lyapunov. O 
resultado é o Teorema \ref{teo:egscaposc}.

\noindent
\textbf{Prova:} 
Omitiremos a dependência temporal de $x_i(t)$ para facilitar a notação.

Considere que a função $V:\mathbb{R}^m\rightarrow\mathbb{R}$ dada por
\begin{equation}
V({x_i}) =
1/2\langle x_i-a,Q(x_i-a) \rangle
\end{equation}
onde $a\in\mathbb{R}^m$ é fixado e $Q$ é uma matriz positiva-definida,
é uma função de Lyapunov para os sistemas isolados, $i=1,\cdots,n$ com um conjunto absorvente
$\Omega$. Considere então o sistema acoplado \eqref{eq:laplacianmodel} que pode ser
escrito na forma compacta
\begin{equation}
\dot{X} = F(X) - \alpha (L\otimes H)X = G(X)
\end{equation}
onde $X=(x_1,\cdots,x_2)$ e $F(X)=(f(x_1),\cdots,f(x_n))$.
Considere também a função $W:\mathbb{R}^{nm}\rightarrow\mathbb{R}$ definida da forma
\begin{equation}
W(X) \colon{=} \left \{ \sum_{i}{V(x_i)} \quad \forall \mbox{ } i \mbox{ tal que } 
x_i\in\mathbb{R}^m\setminus\Omega \right.
\end{equation}
Afirmação: $W$ é função de Lyapunov para $G(X)$.

Por construção, $W$ é positiva-definida em relação a $\Omega^{nm}$.
Portanto, resta-nos estudar o sinal da derivada de $W$. 
Temos 3 casos a analisar, porém antes de mais nada, note que considerando a dinâmica isolada
\begin{align*}
V'({x}_i) & = 1/2\langle dx_i/dt,Q(x_i-a)\rangle + 1/2\langle x_i-a,Qdx_i/dt\rangle\\
& = 1/2\langle f(x_i),Q(x_i-a)\rangle + 1/2\langle x_i-a,Qf(x_i)\rangle =
\langle f(x_i),Q(x_i-a)\rangle.
\end{align*}
E considerando a dinâmica coletiva tem-se que $V'(x_i) = \langle g_i(x_i), Q(x_i-a)\rangle$
onde $g_i(x_i) = f(x_i) - \alpha \sum_{j=1}^n L_{ij}H(x_j)$.

Caso 1: $x_i\in\mathbb{R}^m\setminus\Omega$
para todo $i=1,\cdots,n$. Então, temos que
\begin{align*}
W'(X) & = \sum_{i=1}^n{V'(x_i)} = \sum_i \langle f(x_i)-\alpha\sum_j{L_{ij}H(x_j)},Q(x_i -a) \rangle\\
& = \sum_i \langle f(x_i),Q(x_i-a)\rangle - \alpha\sum_i\sum_j L_{ij} \langle H(x_j),Qx_i\rangle +
\alpha\sum_i\sum_j L_{ij} \langle H(x_j),Qa\rangle\\
& = \sum_i \langle f(x_i),Q(x_i-a)\rangle - \alpha\sum_i\sum_j L_{ij} x_i^\dagger QHx_j +
\alpha\sum_i\sum_j L_{ij} a^\dagger QHx_j,
\end{align*}
porém, note que podemos escrever estas duas últimas somas da forma
$$
\sum_i\sum_j L_{ij} x_i^\dagger QHx_j = X^\dagger(L\otimes QH)X \quad \mbox{e}\quad
\sum_i\sum_j L_{ij} a^\dagger QHx_j = A^\dagger(L\otimes QH)X
$$
onde $A = \mathbf{1}\otimes a$. Assim,
\begin{align*}
A^\dagger(L\otimes QH)X & = (\mathbf{1}\otimes a)^\dagger (L\otimes QH)X\\
& = (\mathbf{1}^\dagger\otimes a^\dagger)(L\otimes QH)X\\
& = (\mathbf{1}^\dagger L) \otimes (a^\dagger QH)X,
\end{align*}
porém, $\mathbf{1}^\dagger L = L\mathbf{1} = \mathbf{0}$, logo $A^\dagger(L\otimes QH)X=0$. Ficando com
$$
W'(X) = \sum_i \langle f(x_i),(x_i-a)\rangle -\alpha X^\dagger(L\otimes QH)X.
$$
Como $L$ é positiva semi-definida e $Q$ e $H$ são positivas-definida então tem-se que
$L\otimes QH$ é positiva semi-definida (Teorema \ref{teo:positivedifinitekronecker}). 
E por hipótese $\sum_i \langle f(x_i),Q(x_i-a)\rangle<0$
visto que $x_i\in\mathbb{R}^m\setminus\Omega$ para todo $i$, portanto $W'(X)<0$.

Caso 2: $x_i\in\Omega$ para todo $i$. 
Como $\Omega$ é um conjunto
positivamente invariante, então não temos nada a fazer pois as trajetórias
permanecem em $\Omega$ por todo o tempo futuro.

Caso 3:
Vamos impor, sem perder a generalidade e
para facilitar as manipulações algébricas,
que $a=\mathbf{0}$.  Considere inicialmente que
$x_1\in\mathbb{R}^m\setminus\Omega$ e $x_2,\cdots,x_n\in\Omega$. 
Assim,
\begin{align*}
W'(X) & = V'(x_1) = \langle f(x_1) - \alpha \sum_j L_{1j}H(x_j), Qx_1\rangle\\
& = \langle f(x_1),Qx_1\rangle - \alpha \left\langle \sum_j L_{1j}H(x_j), Qx_1 \right\rangle.
\end{align*}
Queremos mostrar nesse ponto que $\left\langle \sum_j{L_{1j}H(x_j)}, Qx_1 \right\rangle\geq0$. Perceba que
$$
\left\langle \sum_j{L_{1j}H(x_j)}, Qx_1 \right\rangle = x_1^\dagger\sum_j{L_{1j}QHx_j} =
X^\dagger(L\otimes QH)X - \sum_{i=2}^n\sum_j L_{ij} x_i^\dagger QHx_j.
$$
Em outras palavras, queremos mostrar que cada parcela da soma $X^\dagger(L\otimes QH)X$
é um número não-negativo começando com $i=1$, de forma que o argumento utilizado pode ser
generalizado. Então, sabendo que $L_{ij} = \delta_{ij}k_i - A_{ij}$ onde $k_i$ é o grau
do vértice $i$, temos
\begin{align*}
\left\langle \sum_j{L_{1j}H(x_j)}, Qx_1 \right\rangle & = \left\langle \sum_j{(\delta_{1j}k_1-A_{1j})H(x_j)},Qx_1 \right\rangle\\
& = \left\langle \sum_j{(\delta_{1j}k_1H(x_j)},Qx_1\right\rangle - \left\langle \sum_j A_{1j}H(x_j),Qx_1\right \rangle\\
& = k_1 \langle H(x_1),Qx_1 \rangle -\left\langle \sum_j A_{1j}H(x_j),Qx_1 \right\rangle\\
& = \sum_{j\sim1} \left( \left\langle H(x_1),Qx_1 \right\rangle - \left\langle H(x_j),Qx_1\right\rangle \right)\geq0
\end{align*}
visto que para $j>1$, $x_j\in\Omega$, então, utilizando a desigualdade de Schwarz 
\cite{lima1987curso} é possível
mostrar que $\left\langle H(x_1),Qx_1 \right\rangle \geq \left\langle H(x_j),Qx_1\right\rangle$ para todo $j=1,\cdots,n$.
Assim, cada parcela da soma
$$
X^\dagger(L\otimes QH)X = \sum_i\sum_j L_{ij} x_i^\dagger QHx_j
$$
é um número positivo, pois basta considerar que o vetor $x_i\in\mathbb{R}^m\setminus\Omega$ e
$x_j\in\Omega$ para $j=1,\cdots,n$, $j\neq i$. Então o argumento será o mesmo feito para $x_1$.
Sendo assim, podemos afirmar que para qualquer $1\leq k\leq n$, sendo $x_1,\cdots,x_k\in
\mathbb{R}^m\setminus\Omega$ e $x_{k+1},\cdots,x_n\in\Omega$ teremos que $W'(X)<0$ pois
\begin{align*}
W'(X) & = \sum_{i=1}^k V'(x_i)\\
& = \sum_{i=1}^k \left\langle f(x_i),Qx_i \right\rangle - \alpha \sum_{i=1}^k\sum_j L_{ij} x_i^\dagger QHx_j
\end{align*}
onde $\sum_{i=1}^k\sum_j L_{ij} x_i^\dagger QHx_j\geq0$ e por hipótese,
$\sum_{i=1}^k \left\langle f(x_i),Qx_i \right\rangle <0$.
$\blacksquare$

\section{Sincronização Global em Redes}\label{sec:sgr}


O resultado a ser mostrado é o Teorema \ref{teo:principal1}. Antes porém
considere o seguinte lema cuja prova está disponível no Apêndice \ref{ape:provalemarus}.
\begin{lema}\label{lema:rus}
Considere o sistema \eqref{x=ftx} onde $f:\R^m\rightarrow\R^m$ é de classe $C^d$, $d\geq2$. Suponha
que as trajetórias de \eqref{x=ftx} moram num domínio compacto $\Omega\subset\R^m$. Então,
o resto $R(s,u)$ de uma linearização de $f$, em torno de uma curva solução $s=s(t)$ de 
\eqref{x=ftx}, cumpre
$$
\Vert R(s,u)\Vert \leq M \Vert u\Vert^2
$$
onde $u\in\R^m$ é tal que $x=s+u$ está suficientemente próximo de $s$ e $M$ é uma 
constante uniforme em $s$.
\end{lema}


Ainda antes de iniciar a prova, vamos considerar uma pequena observação sobre o espaço euclidiano
$\R^{nm}$.
Os espaços $\mathbb{R}^n\otimes \mathbb{R}^m$ e $\mathbb{R}^{nm}$ são isomorfos e 
$\mathbb{R}^{nm}$ pode ser decomposto em uma soma direta da forma
$$
\mathbb{R}^{nm} = {N}\oplus {T}
$$
com as projeções
$$
\pi_{{N}} : \mathbb{R}^{nm}\rightarrow {N} \quad \mbox{e} \quad  \pi_{{T}} :
\mathbb{R}^{nm} \rightarrow {T}
$$
onde $\pi_{{N}} +\pi_{{T}} = Id_{\mathbb{R}^{nm}}$. Os subespaços
${N}$ e ${T} \subset\R^{nm}$ são determinados pelos mergulhos de 
$\mathbb{R}^{m}$ e $\mathbb{R}^{(n-1)m}$
induzido pelo laplaciano $L$ da rede. Isto é, tais subespaço são dados por
$$
{N} = \{ \um \}\otimes \R^{m} \quad \mbox{e} \quad {T} = \{ v_2,\cdots,v_n\}
\otimes \R^m
$$
onde $\{ \um \}$ representada o subespaço de $\R^n$ gerado pelo vetor $\um$ e da mesma forma
para $\{ v_2,\cdots,v_n\}$.
\\
\\
\noindent
\textbf{Prova (do Teorema \ref{teo:principal1}):}
Queremos mostrar que a variedade de sincronização ${N}$ \eqref{eq:variedadedesincronizacao} é
localmente atratora. 
Vamos olhar para as equações do movimento $(x_i(t)\in\mathbb{R}^m)$ não separadamente, mas como um
único vetor obtido do empilhamento $X(t)=(x_1(t),\cdots,x_n(t))\in\mathbb{R}^{nm}$, o
que produz o modelo escrito na forma de blocos
\begin{equation}\label{eq:modelocompacto2}
\dot{X} = F(X) - \alpha (L\otimes H)X 
\end{equation}
com $F(X)=(f(x_1),\cdots,f(x_n))$, $L\in\operatorname{Mat}(\mathbb{R},n)$ é o laplaciano
da rede. Para facilitar a notação e a análise
omitiremos a dependência em $t$ de $X(t)$.

Como $L$ é simétrica, e por hipótese $H$ também é simétrica, então essas matrizes 
são ortogonalmente diagonalizáveis (Corolário \ref{cor:teoremaespectral}). Considere então a
representação $L=OMO^\dagger$ e $H = PDP^\dagger$. Da mesma forma, $L\otimes H$ tem a decomposição
$L\otimes H = (O\otimes P) (M\otimes D) (O\otimes P)^\dagger$.
Então, pelo Teorema Espectral,
 como os autovetores de $L\otimes H$ formam uma base de $\mathbb{R}^{nm}$ podemos representar o vetor
$X$ nessa base efetuando a seguinte mudança linear de coordenadas
$$
X = (O\otimes P)Y,
$$
onde $Y= (y_1,\cdots,y_n), y_i\in\mathbb{R}^m$.
Mas, de forma explícita pode-se verificar que
$$
(O\otimes P)Y = \sum_{i=1}^n v_i \otimes Py_i,
$$
onde os $v_i$'s são os autovetores de $L$. 
Estamos considerando
que os autovetores de $L$, $v_1=\mathbf{1}, v_2,\cdots, v_n$, já estão ortonormalizados. 
Para facilitar a notação colocamos $s=Py_1$ e escrevemos então o vetor $X$ na forma
\begin{equation}\label{eq:mudancalineardecoordenadas}
X = \mathbf{1}\otimes s + \sum_{i=2}^n v_i\otimes Py_i = \mathbf{1}\otimes s + U
\end{equation}
com $U=\sum_{i=2}^n v_i\otimes Py_i$.

Introduzimos neste momento a seguinte nomenclatura: chamamos o autovetor $\mathbf{1}$ de $L$
de modo normal e os demais $(v_2,\cdots,v_n)$ de modos transversais. 
Queremos então obter condições para que as contribuições dos modos transversais convirjam ao vetor 
nulo de $\mathbb{R}^{nm}$:
\begin{equation}\label{eq:uazero}
\lim_{t\rightarrow\infty} U = \mathbf{0}.
\end{equation}
Isto é, a rede irá sincronizar globalmente pois  neste caso tem-se que $\lim_{t\rightarrow\infty}X 
= \mathbf{1}\otimes s$, o que significa
$$
x_1(t) = x_2(t) = \cdots = x_n(t) = s(t).
$$

Reescrevendo a equação \eqref{eq:modelocompacto2} na nova representação 
\eqref{eq:mudancalineardecoordenadas} e sabendo que 
$(L\otimes H)(\mathbf{1}\otimes s) = (L\mathbf{1}\otimes H s)=\mathbf{0}$ temos
\begin{align*}
\dot{X} & = F(X) - \alpha (L\otimes H)X\\
\mathbf{1}\otimes \dot{s} + \dot{U} & = F(\mathbf{1}\otimes s + U)-\alpha (L\otimes H)
(\mathbf{1}\otimes s + U)\\
& = F(\mathbf{1}\otimes s + U)-\alpha (L\otimes H)U.
\end{align*}
Neste ponto, nós assumimos que $U$ é pequeno e realizamos a expansão em série de Taylor
de $F(\mathbf{1}\otimes s + U)$ em torno de uma vizinhança da variedade de sincronização 
${N}$, onde $\mathbf{1}\otimes s + U$ pertence a essa vizinhança, produzindo
$$
F(\mathbf{1}\otimes s + U) = F(\mathbf{1}\otimes s) + DF(\mathbf{1}\otimes s)U + R(s,U).
$$
Mas, não é difícil ver que $F(\mathbf{1}\otimes s) = \mathbf{1}\otimes f(s)$.
Note ainda que $DF(X)$ é uma matriz diagonal em blocos onde cada bloco $1\leq i\leq n$ é a matriz
jacobiana $Df(x_i)$ o implica que podemos escrever
$$
DF(\mathbf{1}\otimes s) = I_n \otimes Df(s).
$$
Então, ficamos com
\begin{equation}\label{eq:linearizacaoparcial}
\mathbf{1}\otimes \dot{s}  + \dot{U} = 
\mathbf{1}\otimes f(s) + [I_n \otimes Df(s)]U -\alpha(L\otimes H)U + R(s,U).
\end{equation}

Considere o operador projeção nos modos normais $\pi_{{N}}=\mathbf{1}\cdot
\mathbf{1}^\dagger\otimes I_m$,
e o seu complemento, o operador projeção nos modos transversais $\pi_{{T}} = 
Id_{\R^{nm}}-\pi_{{N}}$.

Assim, analisamos a equação \eqref{eq:linearizacaoparcial} projetando-a
nos modos normais e transversais.
%
Aplicando o operador $\pi_{{N}}$ tem-se, de um lado
$\pi_{{N}}(\mathbf{1}\otimes \dot{s})+\pi_{{N}}(\dot{U})  = 
\mathbf{1}\otimes \dot{s}$. E do outro
$$
\pi_{{N}}\left(\mathbf{1}\otimes f(s)\right) + \pi_{{N}}[I_n \otimes Df(s)]U 
 -\alpha\pi_{{N}}(L\otimes H)U + \pi_{{N}}R(s,U) = \mathbf{1}\otimes f(s)+ \pi_{{N}}R(s,U),
$$
e portanto
\begin{equation}\label{eq:umcomresto}
\mathbf{1}\otimes \dot{s} = \mathbf{1}\otimes f(s)+\pi_{{N}}R(s,U).
\end{equation}

Aplicando-se o operador $\pi_{{T}}$ em \eqref{eq:linearizacaoparcial} tem-se, de um lado
$$
\pi_{{T}}(\mathbf{1}\otimes \dot{s}) + \pi_{{T}}(\dot{U}) = \dot{U}.
$$
Do outro lado tem-se
$$
\pi_{{T}}(\mathbf{1}\otimes f(s)) = \mathbf{0}, \quad 
\pi_{{T}}[I_n\otimes Df(s)]U = [I_n\otimes Df(s)]U\quad
\mbox{e}\quad -\alpha \pi_{{T}}(L\otimes H)U = -\alpha (L\otimes H)U.
$$
Assim ficamos com
\begin{equation}\label{eq:ucomresto}
\dot{U} = [I_n\otimes Df(s) -\alpha (L\otimes H)]U + \pi_{{T}}R(s,U).
\end{equation}
A equação \eqref{eq:ucomresto} é uma equação diferencial linear não-autônoma com perturbação.
Então, como queremos conseguir condições para garantir que $U=\mathbf{0}$ seja uniformemente
assintoticamente estável em \eqref{eq:ucomresto}, faremos o seguinte procedimento: Vamos
considerar a respectiva equação homogênea, ou seja, 
\begin{equation}\label{eq:ulinear}
\dot{U} = [I_n\otimes Df(s) -\alpha (L\otimes H)]U,
\end{equation}
obter uma cota para que o operador de evolução desta equação linear tenha contração uniforme
e a partir de então, considerar a equação perturbada \eqref{eq:ucomresto} e mostrar
que a cota obtida na equação linear continua sendo válida mesmo sem desprezar o resto.

Reescrevendo $U$ na base de $T = \{v_2,\cdots,v_n\}\otimes \{w_1,\cdots,w_n\}$, isto é, 
$U=\sum_{i=2}^n v_i\otimes Py_i$, onde os $w_i$'s são os autovetores ortonormalizados de $H$,
ficamos com
\begin{align*}
\sum_{i=2}^n v_i\otimes P\dot{y}_i & = [I_n\otimes Df(s) -
\alpha (L\otimes H)]\sum_{i=2}^n v_i\otimes Py_i\\
& = \sum_{i=2}^n v_i \otimes [Df(s) - \alpha \lambda_i H]Py_i
\end{align*}
ou ainda
$$
\sum_{i=2}^n v_i \otimes \left\{P\dot{y}_i - [Df(s) - \alpha \lambda_i H]Py_i\right\}=\mathbf{0}.
$$
Como todos os autovetores $v_i$ são linearmente independentes, segue que 
$P\dot{y}_i - [Df(s) - \alpha \lambda_i H]Py_i = \mathbf{0}$ para todo $i$,
%
ou simplesmente
\begin{equation}\label{eq:coordenadasmudancalinear}
\dot{y}_j = [A(t)-\alpha\lambda_jD]y_j
\end{equation}
com $j=2,\cdots,n$ e $A(t) = P^\dagger Df(s(t))P$. Portanto, reduzimos a análise à estudar as equações 
diferenciais dos coeficientes da mudança linear de coordenadas.

Todas as equações, em \eqref{eq:coordenadasmudancalinear}, tem a mesma estrutura diferenciado-se 
apenas pelo $j$-ésimo autovalor de $L$, $\lambda_j$.
Procedemos então neste momento a mesma análise feita na Seção \ref{section:alphacritico},
utilizando o Teorema da diagonal dominante (\ref{teo:criterio1}).
Para que a solução trivial $y_j=\mathbf{0}$ de cada equação seja uniformemente assintoticamente 
estável devemos impor que
$$
A_{ii}(t) - \alpha\lambda_j D_{ii} + \sum_{k=1, \mbox{ }k\neq i}^m \vert A_{ik}(t)\vert <0
\quad\forall\mbox{ } 1\leq i\leq m \quad \mbox{e}\quad 2\leq j\leq n
$$
ou ainda
$$
\alpha >  \frac{A_{ii}(t)+\sum_{k=1, \mbox{ }k\neq i}^m \vert A_{ik}(t)\vert}{\lambda_j D_{ii}} < 
\dfrac{\sum_{k=1}^m \vert A_{ik}(t)\vert}{\lambda_j D_{ii}}.
$$
Como os autovalores de $L$ podem ser
ordenados da forma $0=\lambda_1<\lambda_2\leq\cdots\leq\lambda_n$, 
e, por hipótese, $H$ é positiva-definida, então os elementos da diagonal de $D$ podem ser ordenados 
da forma $0<\mu_1\leq\mu_2\leq\cdots\leq\mu_m$, e as condições iniciais de todas as trajetórias
estão numa vizinhança do domínio absorvente $\Omega$, somos motivamos
a tomar o parâmetro crítico de acoplamento
\begin{equation}
\alpha_c  = \dfrac{\sup_{x\in\Omega}\Vert P^\dagger Df(x) P\Vert_{\infty}}{\lambda_2\mu_1},
\end{equation}
ou simplesmente
$\alpha_c  = {\beta}/{\lambda_2\mu_1}$,
onde $\beta = \sup_{x\in\Omega}\Vert P^\dagger Df(x) P \Vert_{\infty}$.

Portanto, tomando $\alpha>\alpha_c$, 
garantimos que os operadores de evolução $T_j(t,s)$ associados aos $y_j$'s de 
\eqref{eq:coordenadasmudancalinear}
possuem contração uniforme e portanto,
o Teorema \ref{teo:criterio1} assegura que
\begin{align*}
\Vert y_j(t)\Vert & \leq \Vert T_j(t,s)y(s)\Vert\\
& \leq \Vert T_j(t,s)\Vert \Vert y_j(s)\Vert \leq k_j e^{-\eta (t-s)}\Vert y_j(s)\Vert,
\end{align*}
onde $\eta = \alpha\lambda_2\mu_1 - \beta$. Queremos então, conseguir condições para o 
operador de evolução $T(t,s)$ de $U(t)$ em \eqref{eq:ulinear} tenha contração uniforme,
para tanto, vamos utilizar norma euclidiana, sem perder a generalidade, 
para estimar $U(t)$. Note que $\Vert U(t)\Vert_2^2 = \langle U(t),U(t)\rangle$ e que
\begin{align*}
\langle U(t),U(t)\rangle & = \left\langle \sum_{i=2}^n  v_i \otimes Py_i(t), 
\sum_{j=2}^n  v_j \otimes Py_j(t) \right\rangle\\
& = \left\langle \sum_{i=2}^n  v_i \otimes P T_i(t,s)y_i(s), \sum_{j=2}^n  v_j \otimes P T_j(t,s)y_j(s) \right\rangle\\
& = \sum_i\sum_j v_j^\dagger v_i \otimes [y_j(s)^\dagger T_j^\dagger (t,s) P^\dagger ]
P T_i(t,s) y_i(s)\\
& = \sum_i v_i^\dagger v_i \otimes y_i(s)^\dagger T_i^\dagger(t,s) T_i(t,s) y_i(s) =
\sum_{i=2}^n \Vert y_i(t)\Vert^2_2.
\end{align*}
De forma análoga, mostra-se que $\Vert U(s)\Vert_2^2 = \sum_{i=2}^n \Vert y_i(s)\Vert^2_2$.

Como $\Vert U(t)\Vert_2^2 = \sum_{i=2}^n \Vert y_i(t)\Vert^2_2$ e 
$\Vert y_i(t)\Vert^2_2 \leq k_i^2 e^{-2\eta(t-s)}\Vert y_i(s)\Vert^2_2$, onde os $k_i$'s são
constantes que dependem apenas da dimensão do sistema isolado, então
\begin{align*}
\Vert U(t)\Vert_2^2 & \leq \sum_{i=2}^n k_i^2 e^{-2\eta(t-s)}\Vert y_i(s)\Vert^2_2\\
& \leq \left(\max_i k_i^2 \right)e^{-2\eta(t-s)}  \sum_{i=2}^n \Vert y_i(s)\Vert^2_2\\
& = \kappa^2 e^{-2\eta(t-s)} \Vert U(s)\Vert^2_2
\end{align*}
onde $\kappa = \max_i k_i$.
\label{pag:k}
Note que $T(t,s)U(s) = U(t)$. Então
$\Vert T(t,s)U(s)\Vert_2^2 \leq \kappa^2 e^{-2\eta(t-s)} \Vert U(s)\Vert^2_2$, ou ainda
$$
\left\Vert T(t,s)\frac{U(s)}{\Vert U(s)\Vert_2^2}\right\Vert_2^2 \leq \kappa^2 e^{-2\eta(t-s)}.
$$
Seguindo a definição de norma induzida temos que
$$
\Vert T(t,s)\Vert^2_2 = \sup_{\Vert w\Vert_2^2\leq1}\left\Vert T(t,s)
\frac{U(s)}{\Vert U(s)\Vert_2^2}\right\Vert_2^2 \leq \kappa^2 e^{-2\eta(t-s)}
$$
onde $w= U(s)/\Vert U(s)\Vert_2^2$. Portanto, $\Vert T(t,s)\Vert_2\leq k e^{-\eta(t-s)}$.

Considere a equação \eqref{eq:ucomresto} e sua respectiva equação homogênea \eqref{eq:ulinear},
cujo operador de evolução $T(t,s)$ satisfaz $\Vert T(t,s)\Vert \leq \kappa e^{-\eta(t-s)}$.
Como o resto $R(s,U)$ cumpre $\Vert R(s,U)\Vert \leq M\Vert U\Vert^2$ em virtude do Lema
\ref{lema:rus}, então as hipóteses da Proposição
\ref{pro:principio2} são satisfeitas e portanto a estabilidade não muda, isto é,
$U(t)=\mathbf{0}$ é uniformemente assintoticamente estável.
Logo, garantimos que a variedade de sincronização 
é uniformemente assintoticamente estável. $\blacksquare$

\section{Perturbações na Função de Acoplamento}\label{sec:pfa} 

Provamos o Teorema \ref{teo:persistenciaglobal}.


\noindent
\textbf{Porva:} Considere o modelo de perturbação \eqref{eq:mpa}. 
Como os operadores de perturbação
$V_{ij}:\mathbb{R}^m\rightarrow\mathbb{R}^m$, são operadores lineares,
consideramos, sem perder a generalidade, que para todo $i$ e $j$, $V_{ij}$ é uma matriz de 
dimensão $m$. Então, na forma de blocos o modelo \eqref{eq:mpa} 
é escrito como 
\begin{equation}\label{eq:modelopertubacaogeral}
\dot{X} = F(X) - \alpha (L\otimes  H)X +
\left[\sum_{i,j=1}^n (B_iLB_j)(D_{ji}-I_n)\otimes V_{ij}\right] X,
\end{equation}
onde $X=(x_1,\cdots,x_n)$, $F(X) = (f(x_1),\cdots,f(x_n))$, $L$ é o laplaciano da rede,
$B_i\in\operatorname{Mat}(\mathbb{R},n)$ é uma matriz diagonal cuja $i$-ésima entrada da diagonal é 
igual a $1$ e todas as demais entradas são $0$ (zero) e $D_{ij}\in\operatorname{Mat}(\mathbb{R},n)$ é uma
matriz cuja $ij$-ésima entrada é igual a $1$ e todas as demais são nulas.

Realizamos uma análise similar à feita na prova do Teorema \ref{teo:principal1}, isto é, utilizamos
a representação $X = \mathbf{1}\otimes s + U$ %
em \eqref{eq:modelopertubacaogeral} 
e analisamos a parte linear do campo de vetores $F$ com a influência do resto de Taylor,
de forma que procedemos a linearização de
$F(\mathbf{1}\otimes s + U)$ em torno de $\mathbf{1}\otimes s$ considerando que $\mathbf{1}\otimes s + U$
pertence a uma vizinhança de $\mathbf{1}\otimes s$. Assim
\begin{eqnarray*}
 \lefteqn{\mathbf{1}\otimes \dot{s} + \dot{U} = F(\mathbf{1}\otimes s) + DF(\mathbf{1}\otimes s)U - \alpha(L\otimes H)(\mathbf{1}\otimes s + U) +}\\
   &&   +\left[\sum_{i,j=1}^n (B_iLB_j)(D_{ji}-I_n)\otimes V_{ij}\right] 
(\mathbf{1}\otimes s + U)+R(s,U)
\end{eqnarray*}
ou
\begin{equation}
\mathbf{1}\otimes \dot{s} + \dot{U} =
\mathbf{1}\otimes f(s) +  [I_n\otimes Df(s)]U -\alpha (L\otimes H)U +
  \left[\sum_{i,j=1}^n (B_iLB_j)(D_{ji}-I_n)\otimes V_{ij}\right]U +R(s,U) .\label{eq:geralcomperturbacao}
\end{equation}
visto que $F(\mathbf{1}\otimes s)=\mathbf{1}\otimes f(s), 
$ $DF(\mathbf{1}\otimes s)=I_n\otimes Df(s)$, 
$\alpha(L\otimes  H )(\mathbf{1}\otimes s) =\alpha(L\mathbf{1}\otimes  H(s) ) =\mathbf{0}$ e
$$
\left[\sum_{i,j=1}^n (B_iLB_j)(D_{ji}-I_n)\otimes V_{ij}\right] (\mathbf{1}\otimes s) =\mathbf{0}
$$
pois $(B_iLB_j)(D_{ji}-I_n)\mathbf{1} = (B_iLB_j)D_{ji}\mathbf{1} - (B_iLB_j)\mathbf{1}$
onde $(B_iLB_j)D_{ji}$ é uma matriz diagonal onde a $i$-ésima entrada da diagonal é o número
$L_{ij}$ e todas as demais entradas são nulas e portanto o vetor $(B_iLB_j)D_{ji}\mathbf{1}$
tem na sua $i$-ésima entrada o número $L_{ij}$ e todas as demais entradas nulas. Por outro lado,
o vetor $(B_iLB_j)\mathbf{1}$ é exatamente o vetor que tem o número $L_{ij}$ na sua $i$-ésima entrada
e todas as demais nulas. Portanto $(B_iLB_j)D_{ji}\mathbf{1} - (B_iLB_j)\mathbf{1} = \mathbf{0}$.

Novamente, fazemos as projeções de \eqref{eq:geralcomperturbacao} nos modos normais e tranversais
aplicando os operadores $\pi_{{N}}$ e $\pi_{{T}}$ como descritos na Seção 
\ref{sec:sgr}. Aplicando $\pi_{{N}}$ ficamos com $\mathbf{1}\otimes \dot{s} = 
\mathbf{1}\otimes f(s) + \pi_{{N}}R(s,U)$. Aplicando $\pi_{{T}}$ ficamos com
\begin{equation}\label{eq:ucomresto2}
\dot{U} = [I_n\otimes Df(s) -\alpha (L\otimes H)]U +
\left[\sum_{i,j=1}^n (B_iLB_j)(D_{ji}-I_n)\otimes V_{ij}\right]U+\pi_{{N}}R(s,U).
\end{equation}
Porém, perceba que se olharmos apenas para a equação não-perturbada
\begin{equation}\label{eq:ulinear2}
\dot{U} = [I_n\otimes Df(s) -\alpha (L\otimes H)]U
\end{equation}
já sabemos, pelo Teorema \ref{teo:principal1},
que $U(t)$ cumpre $\lim_{t\rightarrow\infty}U(t)=\mathbf{0}$ quando
$\alpha> \beta/(\lambda_2\mu_1)$.
Então, consideramos a equação \eqref{eq:ulinear2} com a perturbação linear 
$ \left[\sum_{i,j=1}^n (B_iLB_j)(D_{ji}-I_n)\otimes V_{ij}\right]U$, 
utilizamos diretamente o Teorema \ref{persistencia} para determinar a magnitude que a referida
perturbação pode ter.
Assim $ \sum_{i,j=1}^n (B_iLB_j)(D_{ji}-I_n)\otimes V_{ij} $ deve cumprir,
$$
\sup_t \left\Vert \sum_{i,j=1}^n (B_iLB_j)(D_{ji}-I_n)\otimes V_{ij}\right\Vert 
 < \frac{\eta}{\kappa}.
$$
onde $\eta = \lambda_2\mu_1 - \beta$ e $\kappa = \max_{2\leq i\leq n} k_i$ 
(ver página \pageref{pag:k} para mais detalhes).
Mas
\begin{align*}
\sup_t  \left\Vert \sum_{i,j=1}^n (B_iLB_j)(D_{ji}-I_n)\otimes V_{ij} \right\Vert 
& \leq \sup_t  \sum_{i,j=1}^n \Vert (B_iLB_j)(D_{ji}-I_n)\otimes V_{ij} \Vert\\
& = \sup_t  \sum_{i,j=1}^n \Vert (B_iLB_j)(D_{ji}-I_n)\Vert \Vert V_{ij} \Vert\\
& \leq \sup_t  \Vert L\Vert (\Vert D_{ji}\Vert +\Vert I_n\Vert) \sum_{i,j=1}^n \Vert V_{ij}\Vert,
\end{align*}
onde $\Vert B_i LB_j\Vert \leq \Vert B_i\Vert \Vert L\Vert \Vert B_j\Vert$, $\Vert B_i\Vert = 1$,
$\Vert D_{ji}\Vert =1$, 
e $\Vert I_n\Vert =1$, isto é, podemos impor que 
\begin{equation}\label{eq:cotaperturbacao}
\sup_t  \sum_{i,j=1}^n\Vert V_{ij}\Vert <  \frac{\eta}{ 2\kappa \Vert L\Vert}, \quad
\mbox{ou melhor,}\quad
\sup_t  \sum_{i,j=1 : i\sim j}^n\Vert V_{ij}\Vert <  \frac{\eta}{ 2\kappa \Vert L\Vert}
\end{equation}
pois, olhando para o modelo de perturbação \eqref{eq:modelopertubacaogeral} verifica-se que 
quando $i$ não é vizinho de $j$ então $B_iLB_j = \mathbf{0}\in\operatorname{Mat}(\R,n)$ de forma que
a equação \eqref{eq:modelopertubacaogeral} pode ser melhor escrita como
$$
\dot{X} = F(X) - \alpha (L\otimes H)X + 
\left[\sum_{i,j=1 : i\sim j}^n (B_iLB_j)(D_{ji}-I_n)\otimes V_{ij}\right] X
$$
e ainda, as contas realizadas não sofrem alterações.
Pelo Teorema \ref{teo:criterio1}, tem-se que $\kappa=1$
se usarmos a norma $\Vert \cdot\Vert_\infty$.
Considerando agora a influência do resto na equação \eqref{eq:ucomresto2}, utilizamos
novamente o Lema \ref{lema:rus} e a Proposição \ref{pro:principio2}, de forma que garantimos
que a estabilidade da variedade de sincronização não é destruída.
$\blacksquare$



\chapter{Conclusões}
\label{cap:conclusoes}

\section{Considerações Finais} 
Ao se abordar o problema da sincronização em redes complexas utilizamos a teoria de
contrações uniformes em equações diferenciais lineares não-autônomas, caracterizando assim uma
nova abordagem sobre um tema tão importante e que está enraizado nos principais ramos
das ciências naturais, exatas, sociais, e enfim na vida humana.

Observa-se que o estado síncrono global dos osciladores difusivamente acoplados
é uma variedade invariante pelas equações do
movimento ao considerar o modelo utilizado \eqref{eq:laplacianmodel}. 
Tal estado
é atingido e garantido ser uniformemente assintoticamente estável devido à força do
parâmetro global de acoplamento $\alpha$, de forma que o seu valor
crítico é determinado unicamente pelas
contribuições da dinâmica intrínseca dos elementos, de propriedades espectrais
da matriz de acoplamento e do segundo autovalor do laplaciano da rede.

O parâmetro crítico de acoplamento, com relação à sincronização, é observacionalmente, e como
já era de se esperar, devido as várias superestimações sobre o mesmo,
muito maior do que o realmente necessário para que tal estado ocorra. Porém, não se garante
que durante o intervalo entre o real $\alpha_c$ e $\beta/(\lambda_2\mu_1)$ o estado síncrono
global seja uniformemente assintoticamente estável - a variedade invariante $N$
\eqref{eq:variedadedesincronizacao} é localmente atratora se $\alpha>\beta/(\lambda_2\mu_1)$
(Teorema \ref{teo:principal1}). 
É também possível verificar esse resultado através das
simulações realizadas em dois osciladores acoplados.

Da mesma forma que verifica-se
a superestimação de $\alpha_c$, a vizinhança de atração de $N$ pode ser bem maior do que
aquela à qual podemos garantir. Além disso, a cota \eqref{eq:criterioperturbacao1} 
para a magnitude da 
perturbação pode apresentar algumas falhas,
uma delas está relacionada à utilização da desigualdade triangular para fazer
superestimações. Além disso, a condição \eqref{eq:criterioperturbacao1} 
é apenas uma condição necessária, isto é, a mesma não garante 
que a estabilidade da variedade de sincronização 
não seja resistente a perturbações com magnitudes fora da cota estabelecida.


\section{Sugestões para Pesquisas Futuras} 

A principal sugestão seria a abordagem da
persistência da sincroniza\c{c}\~ao em redes de sistemas n\~ao-idênticos.
Este caso já foi estudado por Pereira e colaboradores em \cite{tiagocoherence}.
A variedade invariante $N$ 
n\~ao existe neste caso, mas pode-se utilizar resultados 
da teoria de equa\c{c}\~oes diferenciais n\~ao-autônomas para estabelecer a persistência 
do comportamento 
coletivo sob perturba\c{c}\~oes intrínsecas na dinâmica individual dos elementos.

Esta situa\c{c}\~ao \'e mais desafiadora por várias raz\~oes. Por exemplo, como mencionado,
o estado completamente  
sincronizado n\~ao \'e um estado invariante do sistema, e portanto os elementos podem n\~ao 
sincronizar mesmo para acoplamentos fortes.
Pode-se considerar que a não-identidade na dinâmica dos osciladores é na verdade
uma quase identidade, ou seja, seria o mesmo que considerar que
o campo do sistema isolado está sujeito a perturba\c{c}\~oes e então mapeamos
$$
f({ x}_i) \mapsto f({ x}_i) + {w}_i({ x}_i),
$$
onde $\max_i\Vert {w}_i({ x_i)}\Vert \leq \gamma$ seria a máxima magnitude das perturbações.
Mesmo se todos osciladores tem um estado 
inicial igual a perturba\c{c}\~ao ${w}_i$, eles podem ser levados  a estados completamente distintos.
Uma questão em aberto na literatura \'e se as perturba\c{c}\~oes do campo de vetores implicam na 
n\~ao exist\^encia  de solu\c{c}\~oes síncronas.

\renewcommand{\chaptermark}[1]{\markboth{\MakeUppercase{\appendixname\ \thechapter}} {\MakeUppercase{#1}} }
\fancyhead[RE,LO]{}
\appendix

\chapter{Álgebra Linear}
\label{ape:algebralinearA}

\section{Menores Principais}

Nós podemos aferir a definitude\footnote{Consiste da classe dos conceitos das matrizes
positiva-definida, positiva semi-definida, negativa-definida e negativa semi-definida.} 
de uma matriz simétrica calculando os autovalores da mesma. Porém,
outra forma consiste em usar os menores principais.

\begin{definicao}
Seja $A=[A_{ij}]\in\Mat(\mathbb{R},n)$. Para cada inteiro $1\leq k\leq n$, o $k$-ésimo \textbf{menor principal}
de $A$, denotado por $D_k$, é o determinante da matriz
$$
A_k =\left[
\begin{matrix}
A_{11} & A_{12} & \cdots & A_{1k}\\
A_{21} & A_{22}  & \cdots & A_{2k}\\
\vdots & \vdots & \quad & \vdots\\
A_{k1} & A_{k2} & \cdots & A_{kk}
\end{matrix}\right]
$$
a qual se obtém deletando-se as linhas e colunas de ordem $k+1$ até $n$.
\end{definicao}

\begin{teorema}
Seja $A\in\Mat(\mathbb{R},n)$ uma matriz simétrica e $1\leq k\leq n$. Então
\begin{enumerate}
 \item $A$ é positiva-definida se e somente se $D_k >0$ para todo $k$.
\item $A$ é positiva semi-definida se e somente se $D_k \geq0$ para todo $k$.
\item $A$ é negativa-definida se e somente se $(-1)^kD_k <0$ para todo $k$.
\item $A$ é negativa semi-definida se e somente se $(-1)^k D_k \leq 0$ para todo $k$.
\end{enumerate}
\end{teorema}

\noindent
\textbf{Prova:} Ver \cite{zhang1999matrix}, p. 200.

\begin{corolario}\label{cor:matriznegativadefinida}
Seja $A\in \Mat(\R,n)$ uma matriz simétrica. As seguintes afirmações são equivalentes:
\begin{enumerate}
\item $A$ é negativa-definida.
\item Todos os autovalores de $A$ são negativos.
\item Os menores principais de A alternam o sinal começando com 
$D_1=\mbox{det}(A_1)<0$.
\end{enumerate}
\end{corolario}

\newpage
\section{Obtenção do Modelo de Perturbação em Blocos}\label{apes:blocos}

Esta secção trata sobre a obtenção do modelo de perturbação na função de acoplamento o qual,
quando escrito na forma de blocos é resumido na equação \eqref{eq:mpb}.

Considere a equação \eqref{eq:mpl} que pode ser escrita como
\begin{align}
\dot{x}_i & = f(x_i) - \alpha \sum_j L_{ij} H(x_j) - \alpha\sum_j L_{ij}H(x_i) -
 \sum_j L_{ij} V_{ij} x_j + \sum_j L_{ij} V_{ij} x_i\\
& = f(x_i)- \alpha \sum_j L_{ij} H(x_j)- \sum_j L_{ij} V_{ij} x_j + 
 \sum_j L_{ij} V_{ij} x_i.\label{eq:mplextendido}
\end{align}
O objetivo então é emplihar os $n$ vetores de $\mathbb{R}^m$ de cada parcela de 
\eqref{eq:mplextendido}.
Quando do empilhamento da percela $\sum_j L_{ij} H(x_j)$, isto é,
o vetor $\left( \sum_j L_{1j} H(x_j),\cdots,\sum_j L_{nj} H(x_j)\right)\in\mathbb{R}^{nm}$, já sabemos
que o mesmo pode ser escrito como $(L\otimes H)X$ onde $L$ é o laplaciano da rede e
$X= (x_1,\cdots,x_n)$, $x_i\in\R^m$, $i=1,\cdots,n$. Do empilhamento dos vetores $\sum_j L_{ij} V_{ij} x_j$,
pode-se observar que o mesmo é obtido através da multiplicação
$$
\begin{pmatrix}
L_{11}V_{11} & L_{12}V_{12} & \cdots & L_{1n}V_{1n}\\
L_{21}V_{21} & L_{22}V_{22} & \cdots & L_{2n}V_{2n}\\
\vdots & \vdots & \ddots & \vdots\\
L_{n1}V_{n2} & L_{n2}V_{n2} & \cdots & L_{nn}V_{nn}
\end{pmatrix}
\begin{pmatrix}
x_1\\
x_2\\
\vdots\\
x_n
\end{pmatrix}.
$$
Mas a matriz acima descrita pode ser escrita como $\sum_{i,j=1}^n D_{ij}\otimes L_{ij}V_{ij}$ ou ainda
como $\sum_{i,j=1}^n B_{i}LB_j\otimes V_{ij}$ onde
$D_{ij}\in\Mat(\mathbb{R},n)$ é a matriz cuja $ij$-ésima entrada é igual a $1$ e todas as 
demais são nulas e $B_i\in\Mat(\mathbb{R},n)$ é uma
matriz diagonal tal que a $i$-ésima entrada da diagonal é igual a $1$ e todas as demais são nulas.
Note que o produto $B_{i}LB_j$ resulta numa matriz onde o $ij$-ésimo elemento é igual a $L_{ij}$ 
e todos os demais elementos são nulos.

De forma semelhante, o empilhamento das parcelas $\sum_j L_{ij} V_{ij} x_i$ pode ser escrito
pela multiplicação
$$
\begin{pmatrix}
\sum_j L_{1j}V_{1j} & \mathbf{0} & \cdots & \mathbf{0}\\
\mathbf{0} & \sum_j L_{2j}V_{2j} & \cdots & \mathbf{0}\\
\vdots & \vdots & \ddots & \vdots\\
\mathbf{0} & \mathbf{0} & \cdots & \sum_j L_{nj}V_{nj}
\end{pmatrix}
\begin{pmatrix}
x_1\\
x_2\\
\vdots\\
x_n
\end{pmatrix}
$$
onde $\mathbf{0}\in\Mat(\mathbb{R},m)$ é a matriz nula. Esta última matriz pode ser escrita
como $\sum_{i=1}^n D_{ii}\otimes \sum_{j=1}^n L_{ij}V_{ij} = \sum_{i,j=1}^n D_{ii}L_{ij}\otimes
V_{ij}$, ou ainda como $\sum_{i,j=1}^n (B_iLB_j)D_{ji}\otimes V_{ij}$. Note que o papel de $D_{ji}$
é de trazer o único elemento não nulo, provavelmente, de $B_{i}LB_j$ para a diagonal
de sua respectiva linha. Assim, de posse das parcelas do bloco temos o modelo
\begin{align*}
\dot{X} & = F(X) - \alpha (L\otimes H)X - \left[\sum_{i,j=1}^n B_{i}LB_j\otimes V_{ij}\right]X +
 \left[\sum_{i,j=1}^n (B_iLB_j)D_{ji}\otimes V_{ij}\right]X\\
& = F(X) - \alpha (L\otimes H)X +\left[\sum_{i,j=1}^n (B_{i}LB_j)(D_{ji} - I_n)\otimes V_{ij}\right]X.\\
\end{align*}


\chapter{Equações Diferenciais Ordinárias}\label{ape:edoB}

\section{Existência, Unicidade e Extensão}

Considere o sistema de equações diferenciais autônomas
\begin{equation}\label{eq:edoB}
\dot{x} = f(x)
\end{equation}
onde $f:U\subset\mathbb{R}^m\rightarrow\mathbb{R}^m$ é de classe $C^1$, ou seja, suas derivadas
parciais primeiras existem e são contínuas. Abordamos apenas o caso autônomo, porém
os resultados aqui apresentados são extensíveis ao caso não-autônomo, isto é
$f = f(t,x)$.

Uma solução de \eqref{eq:edoB} é o caminho $x:J\rightarrow\mathbb{R}^m$
definido para algum intervalo $J\subset\mathbb{R}$ tal que $\dot{x}(t) = f(x(t))
\mbox{ }\forall\mbox{ } t\in J$. Um valor
inicial para a solução $x:J\rightarrow\mathbb{R}^m$ é uma especificação da forma
$x(t_0) = x_0$ onde $t_0\in J$ e $x_0\in\mathbb{R}^m$. Podemos considerar $t_0=0$.

Uma equação diferencial não-linear pode ter várias soluções que satisfazem
o mesmo valor inicial. Além da questão da unicidade, também temos a questão
da existência. Então, afim de estabelecer condições para que o sistema \eqref{eq:edoB},
sob uma dada uma condição inicial, possua uma única solução, consideramos o 

\begin{teorema}[da Existência e Unicidade]\label{teo:picardlindelof}
Considere o problema de valor inicial
$$
\dot{x} = f(x), \quad x(0) = x_0
$$
onde $x_0\in\mathbb{R}^m$. Suponha que $f:U\subset\mathbb{R}^m\rightarrow\mathbb{R}^m$ é de
classe $C^1$.
Então, exite uma única solução para o referido problema de valor inicial. De
forma mais precisa, existe $c>0$ e uma única solução
$$
x: (-c,c) \rightarrow\mathbb{R}^m
$$
satisfazendo a condição $x(0)=x_0$.
\end{teorema}

\noindent
\textbf{Prova:} Ver \cite{smale}, p. 385.

O resultado apresentado acima garante a existência e unicidade da solução, porém apenas
localmente. O resultado a seguir dá as condições necessárias para que a solução possa
ser estendida, ou seja, para que a solução esteja definida não apenas para uma vizinhança
do tempo inicial, mas por todo o tempo futuro.

\begin{teorema}[Extensão]\label{teo:extensao}
Seja $\mathcal{C}$ um subconjunto compacto de $U\subset\mathbb{R}^m$ e seja
$f:U\subset\mathbb{R}^m\rightarrow\mathbb{R}^m$, $f\in C^1$. Seja $x_0 \in \mathcal{C}$ e suponha
que toda solução da forma $x:[0,b]\rightarrow U$ com $x(0)=x_0$ encontra-se 
inteiramente em $\mathcal{C}$. Então
existe uma solução $x:[0,\infty)\rightarrow U$ satisfazendo $x(0)=x_0$ e $x(t)\in \mathcal{C}$
para todo $t\geq0$, e portanto a solução está definida por todo o tempo
futuro.
\end{teorema}

\noindent
\textbf{Prova:} Ver \cite{smale}, p. 397.

\section{Desigualdade de Grönwall e Variação dos Parâmetros}

O lema que segue é de fundamental importância para obter várias estimativas. A sua prova
pode ser encontrada em \cite{antsaklis}, página 29.

\begin{lema}[de Grönwall]
Sejam $\omega(t), k(t) : J\subset \R \rightarrow \R_+$ funções contínuas com $\omega(t), k(t)\geq0$
para todo $t\in J$ e $C$ uma constante não-negativa. Se
$$
\omega(t) \leq C + \int_s^t{k(u)\omega(u)du} \quad \forall s,t\in J,
$$
ent\~ao
\begin{equation}
\omega(t)\leq C\exp{\left(\int_0^t{k(u)du}\right)}
\end{equation}
para todo $t\in J$.
\end{lema}

A seguinte proposição também é utilizada em algumas demonstrações, justificando portanto
a sua abordagem.

\begin{proposicao}[Variação dos Parâmetros]
Sejam $A:\mathbb{R} \rightarrow \Mat(\mathbb{R},m)$ e 
$g:\mathbb{R} \rightarrow \mathbb{R}^m$
funções contínuas. Considere a equação não-homogênea
$$
\dot{x} = A(t)x +g(t).
$$
A solução de tal equação com a condição inicial $a(s)=x_0$ é dada por
$$
x(t) = T(t,s)x(s) + \int_s^t{T(t,u)g(u)du},
$$
onde $T(t,s)$ é o operador de evolução da equação homogênea correspondente.
\end{proposicao}

\noindent
\textbf{Prova:} Ver \cite{teschl}, p. 81.

\section{Prova do Lema \ref{lema:rus}}\label{ape:provalemarus}

\noindent
\textbf{Prova:}
Tomando $x=s+u$ suficientemente próximo de $s$, seguimos a expansão em série de Taylor
de $f(s+u)$ até ordem $1$. Vamos utilizar a representação $f = (f_1,\cdots,f_m)$ onde
$f_k:\R^m \rightarrow\R, k=1,\cdots,m$. Assim, na expansão, cada componente de $f$ é da forma
$f_k(s+u) = f_k(s) + Df_k(s) u + r_k(s,u)$, isto é,
$$
f(s+u) = (f_1(s),\cdots,f_m(s)) + (Df_1(s)u,\cdots,Df_m(s)u) + (r_1(s,u),\cdots,r_m(s,u)).
$$
Pelo Teorema de Taylor com resto de Lagrange \cite{apostol1962calculus} tem-se que
$r_k(s,u) = (1/2)u^\dagger D^2f_k(s+\theta u) u$, para algum $0<\theta<1$, $\forall k$.
Assim, $\Vert r_k(s,u)\Vert \leq (1/2) \Vert D^2f_k(s+\theta u)\Vert \Vert u\Vert^2$.
Como $f_k\in C^d, d\geq2$, então $D^2f_k(s+\theta u)$ é contínua e $s+\theta u\in\Omega$ (compacto).
Então, pelo Teorema de Weierstrass existe $M_s$ tal que $\Vert D^2f_k(s+\theta u)\Vert \leq M_s$
para cada $s$ fixado. Assim, $M_s:\R^m\rightarrow\R_+$ é uma função contínua que toma valores num
compacto, logo, existe $M$, uniforme em $s$, tal que $\Vert r_k(s,u)\Vert \leq M \Vert u\Vert^2$, 
$\forall k$.
Podemos utilizar a norma $\Vert \cdot\Vert_\infty$, sem perder a generalidade, para estimar
$R(s,u) = (r_1(s,u),\cdots,r_m(s,u))$. Então
$$
\Vert R(s,u)\Vert_\infty = \max_k \Vert r_k(s,u)\Vert \leq M\Vert u\Vert^2_\infty.
$$
$\blacksquare$

\renewcommand{\arraystretch}{0.85}
\captionsetup{margin=1.0cm}  

\newcommand{\etalchar}[1]{$^{#1}$}


\end{document}